\documentclass[11pt]{amsart}

\title[Asymptotics of radiation fields]{Asymptotics of radiation fields in asymptotically Minkowski space}
\author{Dean Baskin}
\author{Andr\'as Vasy}
\author{Jared Wunsch}
\address{Department of Mathematics, Northwestern University}
\curraddr{Department of Mathematics, Texas A\&M University}
\address{Department of Mathematics, Stanford University}
\address{Department of Mathematics, Northwestern University}
\email{dbaskin@math.northwestern.edu}
\email{andras@math.stanford.edu}
\email{jwunsch@math.northwestern.edu}

\usepackage{amsmath, amssymb, verbatim 
}
\usepackage{graphicx}

\newcommand{\abs}[1]{{\left\lvert{#1}\right\rvert}}
\newcommand{\smallabs}[1]{{\lvert{#1}\rvert}}
\newcommand{\norm}[1]{{\left\lVert{#1}\right\rVert}}
\newcommand{\Norm}[2][]{\left\|{#2}\right\|_{#1}}
\newcommand{\ang}[1]{{\left\langle{#1}\right\rangle}}
\newcommand{\smallang}[1]{{\langle{#1}\rangle}}

\newcommand{\av}[1]{}
\newcommand{\jw}[1]{}

\renewcommand{\Im}{\operatorname{Im}}
\renewcommand{\Re}{\operatorname{Re}}
\newcommand{\im}{\operatorname{Im}}
\newcommand{\re}{\operatorname{Re}}
\newcommand{\Ker}{\operatorname{Ker}}

\newcommand{\tP}{\widetilde{P}}
\newcommand{\tu}{\tilde{u}}

\newcommand{\tg}{\tilde{g}}

\newcommand{\tow}{\mathrm{ftr}}
\newcommand{\away}{\mathrm{past}}

\newcommand{\infbdf}{\rho_{\infty}}
\newcommand{\ubdry}{U}

\DeclareMathOperator{\nbd}{nbhd}
\DeclareMathOperator{\supp}{supp}

\newcommand{\CI}{\mathcal{C}^\infty}
\newcommand{\dCI}{\dot{\mathcal{C}}^\infty}

\newcommand{\pa}{{\partial}}
\newcommand{\ep}{{\epsilon}}

\DeclareMathOperator{\id}{I}
\DeclareMathOperator{\Id}{I}
\DeclareMathOperator{\Ann}{Ann}
\newcommand{\RR}{\mathbb{R}}
\newcommand{\reals}{\RR}

\newcommand{\NN}{\mathbb{N}}

\newcommand{\CC}{\mathbb{C}}

\renewcommand\Re{\operatorname{Re}}

\newtheorem{theorem}{Theorem}[section]
\newtheorem{lemma}[theorem]{Lemma}
\newtheorem{proposition}[theorem]{Proposition}
\newtheorem{corollary}[theorem]{Corollary}

\numberwithin{equation}{section}

\newtheorem{non-theorem}[theorem]{Non-Theorem}

\theoremstyle{remark}
\newtheorem{definition}[theorem]{Definition}
\newtheorem{remark}[theorem]{Remark}

\newcommand{\Lap}{\Delta}

\DeclareMathOperator{\WF}{WF}
\DeclareMathOperator{\Op}{Op}
\newcommand\WFb{\WF_{\bl}}

\newcommand{\Tscstar}{{}^{\scl}T^*}

\newcommand{\Tbstar}{{}^{b}T^*}
\newcommand\scl{{\mathrm{sc}}}
\newcommand\bl{{\mathrm{b}}}
\newcommand\Tsc{{}^{\scl} T}
\newcommand\Tb{{}^{\bl} T}
\newcommand\Sb{{}^{\bl} S}
\newcommand\Vf{\mathcal{V}}
\newcommand\Vb{\mathcal{V}_\bl}

\newcommand\Diffsc{\mathrm{Diff}_\scl}
\newcommand\Diffb{\mathrm{Diff}_\bl}
\newcommand\Psib{\Psi_\bl}

\newcommand\Hb{H_\bl}

\newcommand{\sigmab}{\sigma_{b}}
\newcommand{\hamvf}{\mathsf{H}}
\newcommand{\sH}{\hamvf}
\newcommand{\restricted}[1]{\big\rvert_{#1}}
\newcommand{\Ell}{\mathrm{Ell}}
\newcommand{\semi}{\hbar}
\newcommand{\esssupp}{\mathrm{esssupp}}

\newcommand{\cl}{\mathrm{cl}}

\newcommand\cM{\mathcal{M}}
\newcommand\cF{\mathcal{F}}
\newcommand{\cL}{\mathcal{L}}
\newcommand\cS{\mathcal{S}}
\newcommand{\cX}{\mathcal{X}}
\newcommand{\cY}{\mathcal{Y}}
\newcommand{\cR}{\mathcal{R}}
\newcommand{\cI}{\mathcal{I}}
\newcommand{\cU}{\mathcal{U}}

\newcommand{\Mbar}{\overline{M}}

\newcommand{\pd}[1][]{\partial_{#1}}
\newcommand{\sphere}{\mathbb{S}}
\DeclareMathOperator{\Diff}{Diff}

\newcommand{\hol}{\mathcal{H}}
\newcommand{\holcon}[3][-\infty]{\hol(\mathbb{C}_{#2}) \cap
  \smallang{\sigma}^{#1}L^{\infty}L^2(\mathbb{R} ; I^{(#3)}(\Lambda^+))}
\newcommand{\YS}{S_+}

\thanks{The authors acknowledge partial support from NSF grants
  DMS-1068742 (AV) and DMS-1001463 (JW), and the support of NSF
  Postdoctoral Fellowship DMS-1103436 (DB).  We also thank an
  anonymous referee for some helpful remarks on the manuscript.}

\date{December 20, 2012, revised on June 26, 2014.}


\begin{document}

\begin{abstract}
  We consider a non-trapping $n$-dimensional Lorentzian manifold
  endowed with an end structure modeled on the radial compactification
  of Minkowski space.  We find a full asymptotic expansion for
  tempered forward solutions of the wave equation in all asymptotic
  regimes.  The rates of decay seen in the asymptotic expansion are
  related to the resonances of a natural asymptotically hyperbolic
  problem on the ``northern cap'' of the compactification.  For small
  perturbations of Minkowski space that fit into our framework, our
  asymptotic expansions yield
  a rate of decay that improves on the Klainerman--Sobolev
  estimates.
\end{abstract}

\maketitle

\section{Introduction}

\label{sec:introduction}

In this paper we consider the asymptotics of solutions to the wave equation
on a class of spacetimes including asymptotically Minkowski space, as well
as more general spacetimes that have compactifications similar to the
radial compactification of Minkowski space.  Subject to a condition of
non-trapping of null geodesics, namely that they tend toward null-infinity, we find a
full asymptotic expansion for tempered forward solutions of the wave
equation in all asymptotic regimes.  Most notably, we find compound
asymptotics for the solution near null infinity.  The rates of decay seen
in the asymptotic expansion are related to the resonances of a natural
asymptotically hyperbolic problem on the ``northern cap'' of the
compactification.  (This cap corresponds to the interior of the future light cone in
Minkowski space.)  In the special case of small perturbations of Minkowski
space, these expansions imply a rate of decay that
improves on the Klainerman--Sobolev estimates.

\begin{figure}[htp]
  \centering
  \includegraphics{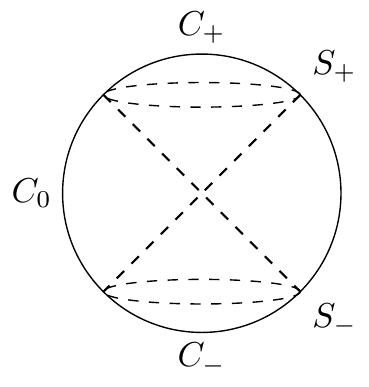}
  \caption{The polar and equatorial regions in Minkowski space}
  \label{fig:Minkowski}
\end{figure}

More specifically, we consider an $n$-dimensional Lorentzian manifold
endowed with the end structure of a ``scattering manifold'' motivated
by the analogous definition for Riemannian manifolds given by Melrose
\cite{RBMSpec}.  Our manifolds come equipped with
compactifications to smooth manifolds-with-boundary, i.e., we will consider
the Lorentzian manifold $(M^\circ,g)$ where $M$ is a manifold with
boundary denoted $X=\pa M.$  The key example is the radial compactification of
Minkowski space $\RR^{1,n-1}_{(t,x)},$ where $X$ is a ``sphere at
infinity,'' with boundary defining function
$\rho=(\abs{x}^2+t^2+1)^{-1/2}.$ 
On $M$ the forward and backward light
cones emanating from any point $q\in M^{\circ}$ terminate at $\pa M$
in manifolds $S_\pm$ independent of the choice of $q;$ we call $S_\pm$
the future and past light cones at infinity, and they bound
submanifolds (which are open subsets) $C_{\pm}\subset X$, which we
call future ($C_+$) and past infinity ($C_-$). In the case of
Minkowski space $C_+$ and $C_-$ are the ``north'' and ``south'' polar
regions (or caps) on the sphere at infinity (see Figure~\ref{fig:Minkowski}). Further, there is an
intermediate region $C_0$ (``equatorial'' on the sphere at infinity in the case of
Minkowski space) which has as its two boundary
hypersurfaces $S_+$ and $S_-$.  We assume that the
metric $g$ is \emph{non-trapping} in the sense that all maximally
extended null-geodesics approach $S_-$ at one end and $S_+$ at the
other.  The full set of geometric hypotheses is described in detail in
Section~\ref{sec:general-hypotheses}.

We consider solutions $w$ to
the wave equation
$$
\Box w=f\in\dCI(M)
$$
on such a manifold so that $w$ is tempered and vanishes near the
``past infinity'' $\overline{C}_{-}$ (thus $f$ also vanishes near
$\overline{C}_{-}$); here $\dCI(M)$ denotes $\CI$ functions on $M$
vanishing at $\pa M$ with all derivatives -- in the case of Minkowski
space this amounts to the set of Schwartz functions. In
\cite[Section~5]{Vasy-Dyatlov:Microlocal-Kerr} the asymptotic behavior
of the solution of the wave equation was analyzed in $C_+$ on
Minkowski space in a manner that extends to our more general setting
in a straightforward manner, giving a polyhomogeneous asymptotic
expansion in the boundary defining function $\rho;$ the exponents
arising in this expansion are related to the {\em resonances} of the
Laplace operator associated to a certain natural {\em asymptotically
  hyperbolic} Riemannian metric on $C_+$.

The main result of this paper is to obtain the {\em precise asymptotic
  behavior} of the solution $w$ near the light cone at infinity,
$S_+=\pa C_+$, performing a uniform (indeed, conormal, on an
appropriately resolved space) analysis as $S_+$ is approached in
different ways.  This amounts to a {\em blow-up} of $S_+$ in $M$.  In
Minkowski space $(t,x) \in \RR^{1+3},$ locally near the interior of
this front face (denoted $\operatorname{ff}$), the blow up amounts to
introducing new coordinates $\rho=(\smallabs{x}^2+t^2+1)^{-1/2}$, $s=t-\smallabs{x},$
$y=x/\smallabs{x}$, the front face itself being given by $\rho=0$, so $s=t-\smallabs{x},$
$y=x/\smallabs{x}$ are the coordinates on the front face.  More generally, if
$\rho$ is a defining function for the boundary at infinity of $M$ and
$v$ is a defining function for $S_+\subset X$ with $(v,y)$ a coordinate system
on $X,$ we can let $s=v/\rho$ and use $s,y$ as coordinates on the
interior of the front face of the blow-up. Thus, $s$ measures the
angle of approach to $S_+$, with $s\to+\infty$ corresponding to
approach from $C_+$, while $s\to-\infty$ corresponds to approach
from $C_0$.
\begin{figure}[htp]
  \centering
  \includegraphics{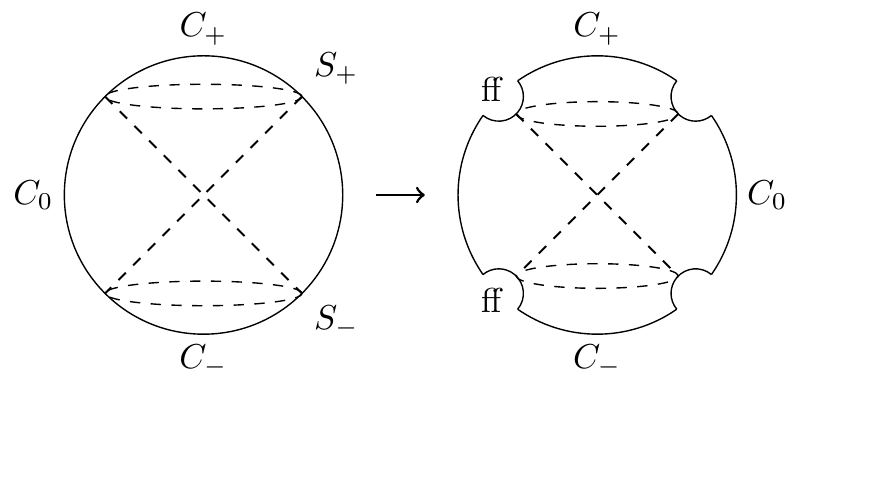}
  \caption{The radiation field blow-up of Minkowski space}
  \label{fig:blow-up}
\end{figure}

In order to make a statement without compound asymptotics, we consider
the so-called radiation fields. Thus, in this paper we show the
existence of the \emph{Friedlander radiation field}, which in the
Minkowski setting is given in the coordinates introduced above by
$$
\mathcal{R}_{+}[w](s,y)=\pa_s\rho^{-(n-2)/2} w(\rho,s,y) |_{\rho=0},
$$
i.e., by restricting an appropriate rescaling of the derivative of $w$
to the new face obtained by blowing up the future light cone at
infinity $S_+$.  (See Section~\ref{sec:asympt-radi-field} for further discussion.)
The function $\mathcal{R}_{+}$ thus measures the
radiation pattern seen by an observer far from an interaction region;
in the case of static metrics, it is known to be an explicit
realization of the Lax-Phillips translation representation as well as
a geometric generalization of the Radon transform
\cite{Friedlander:1980}.

Our main theorem concerns the asymptotics of the radiation field as
$s,$ the ``time-delay'' parameter, tends toward infinity, and more
generally the multiple asymptotics of the solution near the forward light cone.
\begin{theorem}
  \label{thm:main-theorem}
  If $(M,g)$ is a compact manifold with boundary with a
non-trapping Lorentzian scattering metric as defined in
Section~\ref{sec:general-hypotheses} and $w$ is a tempered
  solution of $\Box_{g} w = f \in
  \dCI(M)$ vanishing in a neighborhood of
  $\overline{C}_{-}$,
  then the radiation field of $w$ has an asymptotic
  expansion of the following form as $s\to \infty$:
  \begin{equation*}
    \mathcal{R}_{+}[w](s,y) \sim \sum_{j} \sum_{\kappa \leq m_{j}}
    a_{j\kappa} s^{-\imath\sigma_{j}-1} (\log s)^{\kappa}
  \end{equation*}
Moreover, $w$ has a full asymptotic expansion at all
  boundary faces with the compound asymptotics given by:
\begin{equation}\label{asymp}  \begin{aligned}
   w \sim \rho^{\frac{n-2}{2}}\sum_{j}\sum_{\ell =0}^{N}\sum_{\kappa +\alpha \leq
    \tilde{m}_{j\ell}}\rho^{\ell}\smallabs{\log\rho}^{\kappa}\left( \log\rho +
    \log s\right)^{\alpha} a_{j\ell\kappa\alpha}s^{-\imath\sigma_{j}}
 \end{aligned}\end{equation} 
\end{theorem}

\begin{remark}
  \label{rem:introduction-logs}
  Although it may appear in \eqref{asymp} that
  the log terms may obstruct the restriction to $\rho=0$ and hence the
  definition of the radiation field, we show
  in Section~\ref{sec:logs} that the log terms cancel in the $\ell=0$ term,
  enabling this restriction.

  We also remark that the power of $s$ in the second formula differs
  from the previous one by $1$ due to a derivative in the definition
  of the radiation field.
\end{remark}

\begin{remark}
  In Minkowski space, the requirement that $w$ vanish in a
  neighborhood of $\overline{C}_{-}$ implies that $w$ is the
  \emph{forward} solution of $\Box_{g}w = f$.  One should then think
  of the vanishing requirement as analogous to taking the
  forward solution of $\Box_{g} w = f$.
\end{remark}

\begin{remark}
In Minkowski-like settings, namely when the forward and backward
problems are well-posed in the sense that given an element $f$ of
$\dCI(M)$ supported away from $\overline{C_+}\cup\overline{C_-}$ there
are unique tempered $w_+$, resp.\ $w_-$, with $\Box_g w_\pm=f$ and with $w_+$,
resp.\ $w_-$ vanishing near $\overline{C_-}$, resp.\ $\overline{C_+}$,
one can turn the solution of the Cauchy problem for appropriate space-like hypersurfaces
transversal to $C_0$ (and intersecting $\pa M$ in $C_0$ only) with
Schwartz initial data into the sum of a ``forward solution'' to which our theorem
applies (supported away from $\overline{C_-}$), and a similar ``backward
solution'' (supported away from $\overline{C_+}$), to which the
analogue of our theorem (interchanging $C_+$ with $C_-$)
applies. Thus, in particular, the asymptotic behavior of solutions of the
Cauchy problem with Schwartz initial data is given by our theorem.
A statement related to, but
slightly weaker than, this forward and backward well-posedness follows
from our assumptions (see Remark~\ref{remark:solvability}), but
improving on this, possibly under some additional hypotheses, requires addressing issues beyond the scope of our paper,
and thus will be taken up elsewhere.
\end{remark}

\begin{remark}
It follows from our arguments that the expansion depends continuously
on $w$ and $f$ satisfying a fixed support condition (support in a
fixed compact set disjoint from
$\overline{C_-}$) in the tempered, respectively $\dCI(M)$,
topologies. In particular, finite expansions follow from imposing
finite regularity assumptions, but we do not attempt to make optimal
statements in this paper regarding the required regularity assumptions.
\end{remark}

A question of considerable interest is, of course, whether the
radiation field actually \emph{decays} as $s \to +\infty$ and, more
generally, the description of the exponents $\sigma_j.$ Remarkably,
these are the resonance poles of a naturally-defined asymptotically
hyperbolic metric defined on $C_+.$ (More precisely, the poles we are
interested in are those of the inverses of a family of operators that
looks to leading order like an asymptotically hyperbolic Laplacian.
It is not in general a \emph{spectral family} of the form $P-\sigma^2$
however: the $\sigma$-dependence is more complex.)  We denote the
family of asymptotically hyperbolic operators by $L_{\sigma,+},$ and
record the following corollary:

\begin{corollary}
If there exists $C>0$ such that $L_{\sigma,+}^{-1}$ has no poles at $\sigma$ with
$\Im \sigma  > -C $ then for all $\ep>0$, the radiation field decays at a rate $O(s^{-C-1+\ep})$.
\end{corollary}

One class of spacetimes to which our theorem (and corollary) applies
is that of {\em normally short-range} perturbations of Minkowski
space, i.e., perturbations of the metric which are, relative to the
original metric, $O(\rho^2)$ in the normal-to-the boundary component,
${d\rho^2}/{\rho^4}$, $O(1)$ in the tangential-to-the-boundary
components, ${dv^2}/{\rho^2}$, ${dy^2}/{\rho^2}$ and ${dv\,
  dy}/{\rho^2}$, and $O(\rho)$ in the mixed components.  In
particular, note that we are permitted make \emph{large} perturbations
of the spherical metric on the cap $C_+,$ hence in these
``tangential'' metric components our hypotheses allow a much wider
range of geometries than even traditional ``long-range'' perturbations
of Minkowski space.  

We note here that the long-range structure of the Schwarzschild and
Kerr spacetimes near null infinity does not fit into our class
of spacetimes.  Such spacetimes will be addressed in the follow-up to
this paper.  Also left to a future paper is the question of how to
integrate the decay estimates for geometries with mild (e.g., normally
hyperbolic) trapping into our analysis.

In the more restrictive setting of {\em ``normally very short range''}
perturbations (defined at the beginning of
Section~\ref{sec:minkowski-space-1}), we recover the same
asymptotically hyperbolic problem at infinity as in the Minkowski
case, and therefore exhibit the same order decay as seen on Minkowski
space.  In particular, in odd spatial dimensions one has rapid decay
of solutions of the wave equation away from the light cone.  Thus, we
obtain the following corollary for ``normally very short range''
perturbations of Minkowski space:
\begin{corollary}
  \label{cor:very-short-range}
  If $(M,g)$ is a normally very short range non-trapping perturbation of
  $n$-dimensional Minkowski space, $w$ vanishes near
  $\overline{C}_{-}$, and $\Box_{g}w = f\in
  \dCI(M)$, then the radiation field of $w$ has an
  asymptotic expansion of the following form:
  \begin{equation*}
    \mathcal{R}_{+}[w](s,\omega) \sim
    \begin{cases}
      O(s^{-\infty}) & n \text{ even} \\
      \sum_{j=0}^{\infty}\sum _{\kappa\leq j} s^{-\frac{n}{2}-j}(\log
      s)^{\kappa} a_{j\kappa} & n\text{ odd}
    \end{cases}
  \end{equation*}
More generally in the case of normally short-range perturbations,
given $\ep>0,$ if
the $O(1)$ metric perturbations of the tangential-to-the-boundary
metric components are sufficiently small then the radiation field
still decays as $s \to +\infty$:
\begin{equation*}
  \mathcal{R}_{+}[w](s,\omega) \lesssim s^{\alpha+\ep},\quad \alpha=-\min(2,n/2).
\end{equation*}
\end{corollary}

The polynomial decay of solutions of the wave equation may be compared
with the Klainerman--Sobolev estimates~\cite{Klainerman:1985}.  (We
refer the reader to the book of Alinhac~\cite[Chapter 7]{Alinhac:HPDE} for a more
detailed introduction to such estimates.) In $n$-dimensional
spacetimes where the isometries (and conformal isometries) of
Minkowski space (i.e., the translations, rotations, Lorentz boosts,
and scaling) are ``asymptotic isometries'' (or ``asymptotic conformal
isometries''), then solutions $w$ of the wave equation exhibit decay
of the form
\begin{equation*}
  \left| \pd w(t,r,\theta) \right| \lesssim \frac{1}{(t+r)^{(n-2)/2}(t-r)^{1/2}}.
\end{equation*}
In terms of these coordinates, the asymptotic expansion we obtain
implies that on our class of Lorentzian manifolds (in particular,
on normally short-range perturbations of Minkowski space), there is some $\alpha$ so
that solutions $w$ of the wave equation satisfy
\begin{equation*}
  \left| \pd w(t,r,\theta)\right| \lesssim \frac{1}{(t+r)^{(n-2)/2}} (t-r)^{\alpha}.
\end{equation*}
When there are no eigenvalues of the associated asymptotically
hyperbolic problem, then $\alpha \leq 0$; in particular, on normally
very short range perturbations of Minkowski space (see
Section~\ref{sec:minkowski-space-1}), $\alpha = -n/2$ if $n$ is odd
and $\alpha = -\infty$ if $n$ is even.  Further, the resonances of the
asymptotically hyperbolic problem depend continuously on perturbations
in an appropriate sense.  The operator $P_{\sigma}^{-1}$ introduced
below is stable, but may contain additional poles at certain pure
imaginary integers as compared to the asymptotically hyperbolic
problem (as is the case in even dimensional Minkowski space).
Although such poles do not contribute to the asymptotics of the
radiation field, under perturbations they may become poles of
$L_{\sigma, +}^{-1}$.  Thus, for {\em small} normally short range
perturbations of Minkowski space, $\alpha$ is close to $-\min (2,
n/2)$ (rather than $-\infty$).  (In higher dimensions, one may improve
this statement to obtain $\alpha$ close to $-n/2$ by a careful
analysis of resonant states supported exactly at the light cone.  As
the most interesting case is $n=4$, when $n/2=2$, we do not pursue
this improvement here.) 

The class of spacetimes we consider is geometrically more general than
the class of spacetimes on with the Klainerman--Sobolev estimates
hold, but we require a complete asymptotic expansion of the metric
(and thus considerably more smoothness at infinity).  The methods we
employ would, however, allow also for finite expansions when the
metric has a finite expansion, using more careful accounting.

Finally, we note that a remarkable extension of the radiation field
construction has been obtained by Wang \cite{Wang} to the nonlinear
setting of the Einstein vacuum equations on perturbations of Minkowski
space in spacetimes of dimension $5$ and higher.  This is based in
part on very strong estimates of Melrose-Wang \cite{Melrose-Wang} in
the linear setting for the special case of asymptotically Minkowski
metrics.

\subsection{A sketch of the proof of Theorem~\ref{thm:main-theorem}}
\label{sec:sketch-proof-main}

We start with a tempered solution $w$ of $\Box_{g}w = f'\in
\dCI(M)$ vanishing identically in a neighborhood of
$\overline{C}_{-}$.  We then fix $\chi \in C^{\infty}(M)$ supported
near $\pd M$ so that $\chi$ is identically $0$ near
$\overline{C}_{-}$, identically $1$ near the portion of the boundary
where $w$ is non-vanishing.  In particular, the support of $w\,d\chi$
is compact in $M^{\circ}$ and $\chi w = w$ near $\pd M$.  We then
consider the function $u = \rho^{-(n-2)/2}\chi w$ and set
$$L= \rho^{-2}\rho^{-(n-2)/2}\Box_{g}\rho^{(n-2)/2}.$$
The function $u$ then solves $L u = f$ for some other
function $f\in \dCI(M)$ vanishing near $\overline{C}_{-}$.  A propagation of
singularities argument (Section~\ref{sec:b-regularity}) shows that $u$
is conormal to $\{ \rho = v = 0\}$.

We now set $P_{\sigma} = \widehat{N}(L)$, where $\widehat{N}$
is the reduced normal operator, i.e., the family of operators on the
boundary at infinity obtained by Mellin transform in the normal
variable $\rho.$  If we set $\tilde{u}_{\sigma}$, $\tilde{f}_{\sigma}$ to be
the Mellin transforms of $u$, and $f$, respectively, then
$\tilde{u}_{\sigma}$ solves
\begin{equation*}
  P_{\sigma}\tilde{u}_{\sigma} = \tilde{f}_{\sigma} + \text{ errors},
\end{equation*}
where the additional correction terms arise because $L$
is not assumed to be dilation-invariant.  We show that the operator
$P_{\sigma}$ fits into the framework of
Vasy~\cite{Vasy-Dyatlov:Microlocal-Kerr} and modify the argument of
that paper to show that $P_{\sigma}$ is invertible on certain
variable-order Sobolev spaces (Section~\ref{sec:results-citevasy}).
The argument further shows that $P_{\sigma}^{-1}$ is a meromorphic
family of Fredholm operators with finitely many poles in each
horizontal strip.  In fact, the poles of $P_{\sigma}^{-1}$ may be
identified with resonances for an asymptotically hyperbolic problem
(Section~\ref{sec:ident-p_sigma-1}).

An argument of Haber--Vasy~\cite{Haber-Vasy:Radial} implies that the
residues at the poles of $P_{\sigma}^{-1}$ are $L^{2}$-based
conormal distributions.  In Sections~\ref{sec:conormal-reg} and
\ref{sec:logs} we show that they are in fact classical conormal
distributions and thus have an expansion in terms of $v$.  We
calculate the leading terms of the expansion somewhat explicitly.
Inverting the Mellin transform and iteratively shifting the contour of
integration in the Mellin inversion (Section~\ref{sec:an-asympt-expans})
realizes these residues as the coefficients of an asymptotic expansion
for $u$ in terms of $\rho$.

A slight complication is that not only do the
terms of the expansion become more singular as distributions on $\pa
M$ as one obtains more decay (as is indeed necessary for them to
contribute to the radiation field in the same way, i.e., letting
$\rho\to 0$ with $s=v/\rho$ fixed), but the remainder term also becomes
more singular. We use the a priori conormal
regularity, as shown in Section~\ref{sec:b-regularity}, to deal with
this issue.  The philosophy here is that since the algebra
of b-pseudodifferential operators, discussed in
Section~\ref{sec:b-geometry-mellin} with further references given
there, is not commutative to leading order in the sense of decay at
$\pa M$ (unlike, say, Melrose's scattering pseudodifferential algebra), one first
should obtain regularity in the differential sense, which is
the conormal regularity of Section~\ref{sec:b-regularity}, and then
proceed to obtain decay estimates.

Finally, rewriting the expansion in terms of the
radiation field blow-up $s = v / \rho$ yields an asymptotic expansion
at all boundary hypersurfaces.  The explicit computation of the leading
terms shows that the logarithmic terms match up and thus $u$ may be
restricted to the front face of the blow-up, yielding the radiation
field (after differentiation), and proving
Theorem~\ref{thm:main-theorem} in
Section~\ref{sec:asympt-radi-field}.

\section{b-geometry and the Mellin transform}
\label{sec:b-geometry-mellin}

\subsection{Introduction to b-geometry}
For a more thorough discussion of
b-pseudodifferential operators and b-geometry, we refer the reader to
Chapter 4 of Melrose~\cite{Melrose:APS}. 

\emph{In this section and the following, we initially take $M$ to be a
  manifold with boundary with coordinates $(\rho,y)\in [0,1)\times
  \RR^{n-1}$ yielding a product decomposition $M \supset U \sim [0,1)
  \times \pa M$ of a collar neighborhood of $\pa M.$ In particular,
  for now we lump the $v$ variable in with the other boundary
  variables as it will not play a distinguished role.}

The space of \emph{b-vector fields}, denoted $\mathcal{V}_{b}(M),$ is the
vector space of vector fields on $M$
tangent to $\pd M$.  In local coordinates $(\rho, y)$ near $\pd M$,
they are spanned over $C^{\infty}(M)$ by the vector fields
$\rho\pd[\rho]$ and $\pd[y]$. We note that $\rho\pd[\rho]$ is
well-defined, independent of choices of coordinates, modulo $\rho \mathcal{V}_{b}(M)$; one may call this the
{\em b-normal vector field} to the boundary. One easily verifies that $\mathcal{V}_{b}(M)$ forms a Lie
algebra. The set of b-differential operators, $\Diffb^{*}(M)$, is the
universal enveloping algebra of this Lie algebra:
it is the filtered algebra consisting of operators of the form
\begin{equation}\label{exampleboperator}
A=\sum_{\smallabs{\alpha}+j\leq m} a_{j,\alpha}(\rho,y) (\rho D_\rho)^j
D_y^\alpha \in \Diffb^m(M)
\end{equation}
 (locally near $\pa M$) with the coefficients $a_{j,\alpha} \in \CI(M).$
We further define a bi-filtered algebra by setting
$$
\Diff_{\bl}^{m,l}(M)\equiv\rho^{-l} \Diff_{\bl}^m(M).
$$

The
b-pseudodifferential operators $\Psib^{*}(M)$ are the ``quantization''
of this Lie algebra, formally consisting of operators of the form
$$
b(\rho,y,\rho D_\rho, D_y)
$$
with $b(\rho,y,\xi,\eta)$ a Kohn-Nirenberg symbol; likewise we let
$$
\Psib^{m,l}(M)=\rho^{-l}\Psib^m(M)
$$
and obtain a bi-graded algebra.

\begin{remark}
  The convention we use for the sign of the weight exponent $l$ is the
  opposite of that employed in some other treatments; we have chosen
  this convention as differential order and the weight order behave
  similarly: the space increases if either one of these is
  increased.
\end{remark}

The space $\mathcal{V}_{b}(M)$ is in fact the space of sections of a
smooth vector bundle over $M,$ the \emph{b-tangent bundle}, denoted
$\Tb M.$ The sections of this bundle are of course locally spanned by
the vector fields $\rho\pa_\rho,\pa_y.$ The dual bundle to $\Tb M$ is
denoted $\Tbstar M$ and has sections locally spanned over $\CI(M)$ by
the one-forms $d\rho/\rho, dy.$ We also employ the \emph{fiber
  compactification} $\overline{\Tb^{*}}M$ of $\Tbstar M$, in which we
radially compactify each fiber.  If we let
$$
\xi \frac{d\rho}{\rho}+ \eta \cdot dy
$$
denote the canonical one-form on $\Tbstar M$ then
a defining function for the boundary ``at
infinity'' of the fiber-compactification is
$$
\nu=(\xi^2+\abs{\eta}^2)^{-1/2};
$$
a set of local coordinates on each
fiber of the compactification near $\{v= \rho = 0\}$ is given by
\begin{equation*}
  \nu,\ \hat{\xi} = \nu \xi,\  \hat{\eta} = \nu \eta.
\end{equation*}

The symbols of operators in $\Psib^*(M)$ are thus Kohn-Nirenberg
symbols defined on $\Tbstar M.$ The principal symbol map, denoted
$\sigma_{\bl},$ maps (the classical subalgebra of)
$\Psib^{m,l}(M)$ to $\rho^{-l}$ times
homogeneous functions of order $m$ on $\Tbstar M.$ In the particular
case of the subalgebra $\Diff_{\bl}^{m,l}(M),$ if $A$ is given by
\eqref{exampleboperator} we have
$$
\sigma_{\bl}(\rho^{-l} A)=\rho^{-l} \sum_{\smallabs{\alpha}+j\leq m} a_{j,\alpha}(\rho,y) \xi^j
\eta^\alpha
$$
where $\xi,\eta$ are ``canonical'' fiber coordinates on $\Tbstar M$
defined by specifying that the canonical one-form be 
$$
\xi \frac{d\rho}\rho + \eta \cdot \frac{dy}{\rho}.
$$

In addition to the principal symbol, which specifies high-frequency
asymptotics of an operator, we will employ the ``normal operator'' which
measures the boundary asymptotics.  For a b-differential operator given by
\eqref{exampleboperator}, this is simply the dilation-invariant operator
given by freezing the coefficients of $\rho D_\rho$ and $D_y$ at $\rho=0,$ hence
$$
N(A)\equiv \sum_{\smallabs{\alpha}+j\leq m} a_{j,\alpha}(0,y) (\rho D_\rho)^j
D_y^\alpha \in \Diffb^m([0,\infty) \times \pa M).
$$
The Mellin conjugate (see Section \ref{section:mellin} below) of this operator is known as the ``reduced normal
operator'' and is simply the family in $\sigma$ of operators on $\pa M$
given by
$$
\widehat{N}(A)\equiv \sum_{\smallabs{\alpha}+j\leq m} a_{j,\alpha}(0,y) \sigma^j
D_y^\alpha.
$$
This construction can be extended to b-pseudodifferential operators, but 
we will only require it in the differential setting here.

Here and throughout this paper
we fix a
``b-density,'' which is to say a density which near the boundary is of the
form
$$
f(\rho,y) \abs{\frac{d\rho}\rho\wedge dy_1 \wedge\dots \wedge dy_{n-1}}
$$
with $f>0$ everywhere.
Let $L^2_\bl(M)$ denote the space of square integrable functions with respect to
the b-density.
We let $\Hb^m(M)$ denote the
Sobolev space of order $m$ relative to $L^2_\bl(M)$ corresponding to
the algebras $\Diffb^m(M)$ and $\Psib^m(M)$.  In other words, for
$m\geq 0$, fixing $A\in\Psib^m(M)$ elliptic, one has $w\in\Hb^m(M)$ if
$w\in L^2_\bl(M)$ and $Aw\in L^2_\bl(M)$; this is independent of the
choice of the elliptic $A$.  For $m$ negative, the space is defined by
dualization.  (For $m$ a positive integer, one can alternatively give a
characterization in terms of $\Diffb^m(M)$.)  Let
$\Hb^{m,l}(M)=\rho^{l}\Hb^m(M)$ denoted the corresponding weighted
spaces.

We recall also that associated to the calculus $\Psib^*(M)$ is
associated a notion of Sobolev wavefront set:
$\WFb^{m,l}(w)\subset \Sb^*M$ is defined only for $w\in \Hb^{-\infty,l}$
(since $\Psib(M)$ is not commutative to leading order in
the decay order); the definition is then $\alpha\notin\WFb^{m,l}(w)$ if
there is $Q\in\Psib^{0,0}(M)$ elliptic at $\alpha$ such that $Qw\in
\Hb^{m,l}(M)$, or equivalently if there is $Q'\in\Psib^{m,l}(M)$
elliptic at $\alpha$ such that $Q'w\in L^2_{\bl}(M)$. We refer to
\cite[Section~18.3]{Hormander:v3} for a discussion of $\WFb$ from a
more classical perspective, and \cite[Section~3]{Melrose-Vasy-Wunsch:Propagation}
for a general description of the wave front set in the setting of various pseudodifferential
algebras; \cite[Sections~2 and 3]{Vasy:corners} provide
another discussion, including on the b-wave front set relative to
spaces other than $L^2_{\bl}(M)$.

\subsection{Scattering geometry}
We now turn to the analogous concepts of ``scattering geometry'' which
will be less used in this paper but which are a useful motivation.
For a full discussion of scattering geometry, we
refer the reader to Melrose~\cite{RBMSpec}.

In analogy to the space of b-vector fields, we define \emph{scattering
  vector fields} as $\mathcal{V}_{\scl}\equiv \rho\mathcal{V}_b;$ that
is to say, the vector fields when applied to $\rho$ must return a
smooth function divisible by $\rho^2.$  They are locally spanned by
$\rho^2 \pa_\rho$ and $\rho \pa_y$ over $\CI(M).$ We note that as the b-normal
vector field $\rho\pa_\rho$ is well-defined modulo
$\rho\mathcal{V}_b$, the {\em span} of $\rho^2 \pa_\rho$ is
well-defined modulo $\rho\mathcal{V}_{\scl}$; we call vector fields
lying in this span {\em scattering normal vector fields}. The scattering
vector fields form sections of a bundle
$\Tsc M;$ the dual bundle, $\Tscstar M$ has sections locally spanned
by $d\rho/\rho^2,\ dy/\rho.$  As motivation for our discussions of the
form of the ``scattering metrics'' below, we remark that if we
radially compactify Euclidean space, the \emph{constant} vector fields
push forward to be scattering vector fields on the compactification,
hence sections of the tensor square of $\Tscstar M$ are the natural
place for asymptotically Euclidean or Minkowskian metrics to
live.

 The \emph{scattering differential
  operators} are those of the form (near $\pa M$)
$$
\sum_{\smallabs{\alpha}+j\leq m} a_{j,\alpha}(\rho,y) (\rho^2 D_\rho)^j
(\rho D_y)^\alpha \in \Diffsc^m(M).
$$
Again, this space of operators can be microlocalized by introducing
\emph{scattering pseudodifferential operators} which are formally
objects given by
$$
b(\rho,y,\rho^2 D_\rho,\rho D_y)
$$
with $b(\rho,y,\xi,\eta)$ a Kohn-Nirenberg symbol on the bundle
$\Tscstar M.$  There are of course associated scales of Sobolev
spaces, which we will not have occasion to use in this paper, as well
as wavefront sets which are described in detail in \cite{RBMSpec}.

\subsection{Mellin transform}\label{section:mellin}

We first recall the definition of the \emph{Mellin transform} on $\reals_{+}$.
For a smooth compactly supported function, or indeed a Schwartz function, $u$ on $\reals_{+}$,
$\tilde{u}_{\sigma} = \int _{0}^{\infty}\rho^{-\imath\sigma -1}u(\rho)
\,d\rho$. (Here Schwartz means that for all $k\in\RR$ all derivatives of $u$ are bounded
by a multiple of $\rho^k$ both at $0$ and at $\infty$, i.e.\ they
vanish rapidly.)
Because $u$ is compactly supported or Schwartz, $\tilde{u}_{\sigma}$ is
an entire function of $\sigma$ which decays rapidly along each line of
constant $\Im \sigma$.   We will also use the notation
$$
\mathcal{M} u =\tilde{u}
$$
for the Mellin transform.

The Mellin transform on $\RR_+$ is equivalent to the Fourier transform by the
substitution $x = \log \rho$; note that the Schwartz behavior amounts
to superexponential decay at $\pm\infty$ in terms of $x$.  In particular, the Plancherel theorem
guarantees that it extends to an isomorphism of Hilbert spaces
\begin{equation*}
  L^{2}(\reals_{+}; d\rho/\rho) \to L^{2}(\RR),
\end{equation*}
and, more generally, to an isomorphism with a weighted space,
\begin{equation*}
  \rho^{\delta}L^{2}(\reals_{+};d\rho/\rho) \to L^{2} ( \{ \Im \sigma = -\delta\}).
\end{equation*}
Moreover, the Mellin transform intertwines $\rho\pd[\rho]$ with
multiplication by $\imath\sigma$:
\begin{equation*}
  \widetilde{(\rho\pd[\rho]u)}_{\sigma} = \imath\sigma \tilde{u}_{\sigma}
\end{equation*}
The inverse Mellin transform is given by integrating
$\tilde{u}_{\sigma}\rho^{i\sigma}$ along a horizontal line $\{ \Im
\sigma = C\}$, provided this integral exists.

Near the boundary of $M,$ we use the boundary defining function $\rho$ to
obtain a local product decomposition: $M\supset\nbd(X) = [0, \epsilon)_{\rho} \times
\pd M$.  This local product decomposition allows us to define the
Mellin transform for functions supported near $\pd M$ via cut-off
functions that are identically $1$ for $\rho \leq \epsilon / 2$.  In
what follows, this definition suffices, as we may always cut off the
functions in which we are interested away from the boundary.  Note
that this definition of the Mellin transform depends both on the boundary
defining function $\rho$ and on the cut-off functions chosen, but this
dependence will not make a difference in the sequel.

We additionally recall the space of $L^{2}$-based conormal
distributions $I^{(s)}$ on the boundary $X=\pa M.$  Here we finally split the
boundary coordinates locally into $(v,y) \in \RR\times \RR^{n-2}$
rather than using $y$ to denote all of them.   For the hypersurface $Y = \{
v=0\} \subset X$, $u \in I^{(s)}(N^{*}Y)$ means that $u \in
H^{s}(X)$ and $A_{1}\ldots A_{k}u \in H^{s}$ for all $k$ and for
all $A_{j} \in \Psi^{1}(X)$ with principal symbol vanishing on
$N^{*}Y$.

We now record some additional mapping properties of the Mellin
transform:
\begin{definition}
  Let $\CC_\nu$ denote the halfspace $\Im \sigma>-\nu$ and let
  $\hol(\CC_\nu)$ denote holomorphic functions on this space.  For a
  Fr{\'e}chet space $\mathcal{F}$, let $$\hol(\CC_\nu)\cap \langle
  \sigma\rangle^{-k} L^\infty L^{2}(\reals ;\mathcal{F})$$ denote the
  space of $g_\sigma$ holomorphic in $\sigma \in \CC_\nu$ taking values in $\mathcal{F}$
  such that each seminorm
  $$
  \int_{-\infty}^\infty \norm{g_{\mu+\imath\nu'}}_\mathcal{\bullet}^2 \langle \mu
  \rangle^{2k} \, d\mu
  $$
  is uniformly bounded in $\nu'>-\nu$.
\end{definition}
Note the choice of signs: as $\nu$ increases, the halfspace gets
larger.

We will further allow elements of $\hol (\CC_{\nu})$ to take values in
\emph{$\sigma$-dependent} Sobolev spaces, or rather Sobolev spaces
with $\sigma$-dependent norms.  In particular, we allow
values in the standard semiclassical Sobolev spaces $H^{m}_{h}$ on a
compact manifold (without boundary), with semiclassical parameter $h =
|\sigma|^{-1}$. Recall (see \cite[Section 8.3]{Zworski}) that these are the standard Sobolev spaces
and up to the equivalence of norms, for $h$ in compact subsets of
$(0,\infty)$, the norm is just the standard $H^m$ norm, but
the norm is $h$-dependent: for non-negative integers $m$, in
coordinates $y_j$, locally the norm $\|g\|_{H^m_h}$ is equivalent to
$\sqrt{\sum_{\smallabs{\alpha}\leq m}
  \|(hD_{y_j})^\alpha g\|^2_{L^2}}$.

We will require some more detailed information about mapping properties of the
Mellin transform acting on b-Sobolev spaces.
\begin{lemma}\label{lemma:mellinspaces}
Let $u \in \Hb^{m,l}(M)$ be supported in a collar neighborhood
$[0,\epsilon)_\rho\times X$ of $\pa
M$  Then $$\mathcal{M} u \in \hol(\CC_l) \cap 
\ang{\sigma}^{\max(0,-m)}L^\infty L^2(\RR; H^m(X)).$$  

If $u \in \Hb^{m,l}(M)$ is furthermore conormal to $\rho=v=0$ then
$$
\mathcal{M} u \in \hol(\CC_l) \cap \ang{\sigma}^{-\infty} L^\infty L^2(\RR; I^{(m)}(N^{*}Y)).
$$

The inverse Mellin transform maps
$$
\hol(\CC_{\nu}) \cap
\langle\sigma\rangle^{-\infty}L^\infty L^2(\reals ;I^{(s)}(N^* Y)))
$$
into
$$
\rho^{\nu-0} H_b^\infty([0,\infty)_\rho; I^{(s)}(N^* Y))
$$
which in turn, for $s<1/2$, is contained in
$$
\rho^{\nu-0} v^{s-1/2-0} L^\infty.
$$
\end{lemma}
 \begin{proof}
For $m$ a positive integer, the first result follows since lying in
$\Hb^{m,l}$ implies that 
$$
\pa_y^\alpha\pa_v^\beta u \in \rho^{l} L^2_b
$$
for all $\abs{\alpha}+\abs{\beta}\leq m$
hence
$$
\pa_y^\alpha \pa_v^\beta \mathcal{M} u \in \hol(\CC_l) \cap 
L^\infty L^2(\RR; L^2(X));$$
the result for general $m\geq 0$ follows by interpolation. For $m<0$,
choose a positive integer $\tilde m$ such that $m+\tilde m\geq 0$;
then $u$ can be written as a finite sum of terms of the product of at
most $\tilde m$ b-vector fields applied to elements $u'$ of
$\Hb^{m+\tilde m,l}(M)$. Now, the Mellin transform of such $u'$ lies
in $\hol(\CC_l) \cap L^\infty L^2(\RR; H^{m+m'}(X))$ by the first part;
$\pd[y]$ and $\pd[v]$ act as vector fields on $X$ and thus would lead to
the conclusion that $u$ is in $\hol(\CC_l) \cap L^\infty L^2(\RR;
H^{m}(X))$ if only they appeared; however, $\rho \pd[\rho]$ Mellin
transforms to $\imath\sigma$, and thus we may obtain up to $\tilde m$
factors of $\sigma$ as well, leading to the desired weight when $m<0$
is an integer; interpolation gives the weight (without a loss) for all $m<0$.

The proof of the second and third parts is similar; here we use
Sobolev embedding, and the
fact that regularity under
$\rho \pa_\rho,$ $v \pa_v$ and $\pa_y$ intertwines under Mellin
transform with regularity under $\sigma,$ $v\pa_v,$ and $\pa_y.$
\end{proof}
We remark further that Mellin transform maps
$\Hb^{\infty,l}(M)$ into 
$$
\hol(\CC_l) \cap \ang{\sigma}^{-\infty} L^\infty L^2(\RR; H^\infty(X)).
$$
This map is not onto, as there is no iterated regularity
under $\rho \pa_v$ built into the latter space.

\section{Geometric set-up}
\label{sec:geometry}

\subsection{Minkowski metric}
\label{sec:minkowski-metric}

As a preliminary to our discussion of Lorentzian scattering metrics,
we record the asymptotic behavior of the Minkowski space on
$\RR^{n}$, endowed with the Lorentzian metric with the \emph{mostly
  minus} sign convention (here we are following the notation of
\cite{Vasy-Dyatlov:Microlocal-Kerr}).  We take coordinates $t, x_{1}, \ldots, x_{n-1}$, and
set
\begin{align*}
  t &= \rho^{-1}\cos \theta, \\
  x_{j}&= \rho^{-1}\omega_{j}\sin\theta ,
\end{align*}
with $\omega \in \sphere^{n-2}$.  The Minkowski metric is then
\begin{align*}
  dt^{2} - \sum dx^{2} &= \left( -\frac{\cos\theta
    \;d\rho}{\rho^{2}}-\sin\theta\frac{d\theta}{\rho}\right)^{2} - \sum
  \left( -\omega_{j}\sin\theta\frac{d\rho}{\rho^{2}} + \omega_{j}
  \cos\theta\frac{d\theta}{\rho} + \sin\theta
  \frac{d\omega_{j}}{\rho}\right)^{2} \\
  &= \cos 2\theta \frac{d\rho^{2}}{\rho^{4}} - \cos 2\theta
  \frac{d\theta^{2}}{\rho^{2}} + \sin 2\theta \left(
  \frac{d\rho}{\rho^{2}} \otimes \frac{d\theta}{\rho} +
  \frac{d\theta}{\rho}\otimes \frac{d\rho}{\rho^{2}}\right) -
  \sin^{2}\theta \frac{d\omega^{2}}{\rho^{2}}.
\end{align*}
Here $d\omega^{2}$ represents the standard round metric on the sphere.

As the function $\cos 2\theta$ clearly plays an important role here,
we set
\begin{equation*}
  v = \cos 2\theta,
\end{equation*}
replacing the $\theta$ coordinate by $v$, and write
\begin{equation}\label{Minkowskimetric}
  g = v \frac{d\rho^{2}}{\rho^{4}} - \frac{v}{4(1-v^{2})} \frac{dv^{2}}{\rho^{2}} -
  \frac{1}{2} \big( \frac{d\rho}{\rho^{2}}\otimes \frac{dv}{\rho} +
  \frac{dv}{\rho} \otimes \frac{d\rho}{\rho^{2}}\big) - \frac{1-v}{2}\frac{d\omega^{2}}{\rho^{2}}.
\end{equation} 
We remark that this form of the metric in these extremely natural
coordinates does not conform to the standard ``scattering metric''
hypotheses \cite{RBMSpec} often employed in the Riemannian
signature, in which cross terms of the form $(d\rho / \rho^{2})
\otimes (dy/\rho)$ with $y$ a general smooth function are forbidden.

\subsection{General hypotheses}
\label{sec:general-hypotheses}

Let $(M,g)$ be an $n$-dimensional manifold with boundary $X=\pa M$ equipped with a Lorentzian metric
$g$ over $M^{\circ}$ such that $g$ extends to be a nondegenerate
quadratic form on $\Tsc M$ of signature $(+,-,\ldots, -)$.

\begin{definition}
  We say that $g$ is a Lorentzian scattering metric if $g$ is a smooth,
  Lorentzian signature, symmetric bilinear form on $\Tsc M$, and there
  exist a boundary defining function $\rho$ for $M$, and a function
  $v\in\CI(M)$ such that
\begin{enumerate}
\item 
  When $V$ is a scattering normal
  vector field, $g(V,V)$ has the same sign as $v$ at $\rho=0$,
\item 
  in a neighborhood of $\{v=0,\ \rho=0\}$ we have
  \begin{equation*}
    g = v \frac{d\rho^{2}}{\rho^{4}} - \big(
    \frac{d\rho}{\rho^{2}}\otimes \frac{\alpha}{\rho} +
    \frac{\alpha}{\rho}\otimes \frac{d\rho}{\rho^{2}}\big) - \frac{\tilde{g}}{\rho^{2}}
  \end{equation*}
  with $\alpha$ a smooth 1-form on $M$
  and $\tilde{g}$ a smooth
  symmetric $2$-cotensor on $M$ so that
  \begin{equation*}
    \tilde{g} \rvert_{\Ann(d\rho, dv)} \text{ is positive definite.}
  \end{equation*}
  We further require that
  \begin{equation*}
    \alpha = \frac{1}{2}dv + O(v) + O(\rho) \text{ near } v = \rho = 0.
  \end{equation*}
\end{enumerate}
\end{definition}
\begin{remark}
  We remark that while it might be tempting to mandate also the
  vanishing of the $dv^{2}/\rho^{2}$ component at $v=0$ as we have in
  the exact Minkowski case, this condition is highly non-invariant, in
  that it requires a product decomposition of $X$.  
\end{remark}
\begin{remark}
  The function $v$ must have a non-degenerate $0$-level set when
  restricted to $X$ (and $dv,\ d\rho$ must be independent at that set), since otherwise our metric would be
  degenerate at $v=0$.  

  We further remark that our hypotheses imply that for {\em
    non-trapping} Lorentzian scattering metrics, even if the boundary
  is disconnected, $v$ vanishes on each component of the boundary, see
  Remark~\ref{remark:v-vanishes}. Without the non-trapping assumption
  this need not be the case: consider $\tilde X\times\overline{\RR}$,
  with $\tilde X$ compact without boundary, $\overline{\RR}$ the
  radial compactification of $\RR_r$ (so $\rho=r^{-1}$ works for $r\gg
  1$). Then $dr^2-(1+r^2)h$, $h$ a metric on $\tilde X$, is a
  Lorentzian scattering metric if one chooses $v\equiv 1$; $X$ is then
  the disjoint union of two copies of $\tilde X$.
\end{remark}
\begin{remark}
Note that near $v=0$, $V=\rho^2\pa_\rho$ gives $g(V,V)=v$, which has
the same sign as $v$, so the first and second parts of the definition
are consistent, with the second part refining the first near $v=0.$
\end{remark}

Before proceeding, note that the rescaled, or scattering, Hamilton
vector field of the metric function on $\Tscstar M\setminus o$ is a $\CI$ vector
field, tangent to the boundary. The integral curves of this Hamilton
vector field within the zero set
of the metric function (i.e., the null bicharacteristics)
over the interior of $M$ project to reparameterized null-geodesics; indeed, they are
exactly the appropriately reparameterized 
null-geodesics lifted to $T^*M^\circ$. We show later in Section~\ref{sec:b-radial-points}
that over $S=\{v=0,\rho=0\}$ the Hamilton flow has sources and sinks;
there we shift to the b-framework, and these sources and sinks are located at
the ``b-conormal bundle'' of $S$, denoted by $\cR$.

With this in mind, we make two additional global assumptions on the
structure of our spacetime:
\begin{definition}\label{definition:non-trapping}
A Lorentzian scattering metric is {\em non-trapping} if
\begin{enumerate}
\item The set $S=\{ v= 0,\rho=0\}\subset X$ splits into $S_{+}$ and $S_{-}$, each a
  disjoint union of connected components; we further assume that
  $\{v>0\}$ splits into components $C_\pm$ with $S_\pm=\pa C_\pm.$  We
  denote by $C_{0}$ the subset of $X$ where $v < 0$.
\item The projections of all null bicharacteristics on $\Tscstar M\setminus
  o$ tend to $S_{\pm}$ as their parameter tends to
  $\pm \infty$ (or vice versa).  (A discussion of the flow near
    $S_{\pm}$ is contained in
    Sections~\ref{sec:locat-radi-points}-\ref{sec:b-radial-points}.)
  \end{enumerate}
\end{definition}

In particular, this implies the time-orientability of $(M,g)$ by
specifying the future light cone as the one from which the forward (in
the sense of the Hamilton flow) bicharacteristics tend to $S_+$.

\begin{remark}\label{remark:v-vanishes}
For non-trapping Lorentzian scattering metrics, $v$ must necessarily
vanish on each component of $X$. To see this, note that on $\Tscstar
M$ a Lorentzian metric has non-trivial characteristic set over each
point, in particular over each point in $X$. Since the scattering
Hamilton vector field is a smooth vector field on $\Tscstar X$ tangent
to the boundary, the bicharacteristics through a point in a connected
component of $X$ stay in that component. Thus, for non-trapping
metrics the zero set of $v$ within each connected component must be
non-trivial.
\end{remark}

\begin{remark}\label{remark:solvability}
  $C$ stands for ``cap'' as in the Minkowski case $C_+$ is simply the
  spherical cap $\smallabs{\theta}<\pi/4.$ The assumption the $S_+$
  bounds a cap is in fact not necessary for us to prove any
  of the Fredholm properties in Section \ref{sec:results-citevasy}; however
  it is of course necessary to recognize the poles of the resulting
  operator as resonance poles on a cap, and hence in order to know
  that there are finitely many resonances in any horizontal strip in
  $\CC,$ which is crucial to the development of our asymptotic
  expansion.

Given Schwartz $f$ supported away from $\overline{C_-},$ it is natural to
consider ``forward tempered'' solutions of $\Box_g u=f$ with $u=0$
  near $\overline{C_-}.$  It is unclear whether our hypotheses guarantee
  that these always exist, though it is not hard to show
  that there are only finite dimensional obstructions to solvability
  and uniqueness in fixed weighted spaces.
\end{remark}

Near $v=0$, which is away from the critical points of $v|_{X}$, we may
choose $y_{1},\ldots ,y_{n-2}\in \CI(M)$ so that $(v,y)$ constitute a
coordinate system on $X=\pa M$ and $(\rho, v, y)$ thus give coordinates
on $M$ in a neighborhood of $X$.  Moreover, $(\rho, v, y)$ also
provide a product decomposition of that neighborhood into
$[0,\epsilon)_{\rho}\times X$.  In the frame
\begin{equation*}
  \rho^{2}\pd[\rho],\, \rho \pd[v], \, \rho\pd[y],
\end{equation*}
associated to these coordinates, the metric (when restricted to the
boundary $\{ \rho = 0\}$) thus has the block form
\begin{equation}
  \label{eq:metric-at-bdry}
  G_{0} =
  \begin{pmatrix}
    v & -\frac{1}{2} + a_{0}v & a_{1}v & \dots & a_{n-2}v \\
    -\frac{1}{2} + a_{0}v & b & c_{1} & \dots & c_{n-2} \\
    a_{1}v & c_{1} & -h_{1,1} & \dots & -h_{n-2,1} \\
    \vdots & \vdots & \vdots & \ddots & \vdots \\
    a_{n-2}v & c_{n-2} & -h_{1,n-2} & \dots & -h_{n-2,n-2}
  \end{pmatrix},
\end{equation}
with the lower $(n-1)\times (n-1)$ block negative definite, hence
$h_{ij}$ is positive definite.

Blockwise inversion shows that in the frame
\begin{equation*}
  \frac{d\rho}{\rho^{2}}, \, \frac{dv}{\rho}, \, \frac{dy}{\rho},
\end{equation*}
the inverse metric when restricted to the boundary has the block
form (the $\alpha$ here is a function and should not be
  confused with the 1-form $\alpha$ in the definition of the
  metric)
\begin{equation*}
  G_{0}^{-1} =
  \begin{pmatrix}
    -q & -2 + \alpha v & -\frac{1}{2}\Upsilon^{\mathrm{T}} + O(v)\\
    -2 + \alpha v & -4 v + \beta v^{2} & -v
    \Upsilon^{\mathrm{T}} + O(v^{2})\\
    -\frac{1}{2}\Upsilon + O(v) & -v\Upsilon + O(v^{2}) & -h^{-1} + O(v)
  \end{pmatrix}.
\end{equation*}
In the above, $h^{-1}=h^{ij}$ is the inverse matrix of $h_{ij}$, $q$,
$\alpha$, $\beta$, and $\Upsilon_{j}$ are smooth near $v
= \rho = 0$, and $A^{\mathrm{T}}$ denotes the transpose of the matrix
$A$.  

In a neighborhood of the boundary, i.e., at $\rho \neq 0$, there are
further correction terms in the inverse metric as the actual metric
is given by
\begin{align*}
  G &= G_{0} + H, \\
  H &=
  \begin{pmatrix}
    O(\rho^{2}) & O(\rho) & O(\rho) \\
    O(\rho) & O(\rho) & O(\rho) \\
    O(\rho) & O(\rho) & O(\rho)
  \end{pmatrix}.
\end{align*}
Thus in the inverse frame above, 
\begin{equation}
  \label{eq:inverse-perturbation}
  G^{-1} = G_{0}^{-1} +
  \begin{pmatrix}
    O(\rho) & O(\rho) & O(\rho) \\
    O(\rho) & O(\rho^{2}) + O(\rho v) & O(\rho) \\
    O(\rho) & O(\rho) & O(\rho) 
  \end{pmatrix}.
\end{equation}

Thus in the coordinate frame $\pd[\rho]$,
$\pd[v]$, $\pd[y]$, the dual metric becomes
\begin{equation}\label{dualmetric}
  \begin{pmatrix}
    g^{\rho\rho}\rho^{4}  + O(\rho^{5}) & g^{\rho v} \rho^{3} +
    O(\rho^{4}) & g^{\rho y} \rho^{3} + O(\rho^{4}) \\
    g^{\rho v} \rho^{3} + O(\rho^{4}) & g^{vv}\rho^{2}+ O(\rho^{4}) +
    O(\rho^{3}v) & g^{v y}\rho^{2} + O(\rho^{3}) \\
    g^{\rho y}\rho^{3} + O(\rho^{4}) & g^{vy}\rho^{2} + O(\rho^{3}) &
    g^{yy}\rho^{2} + O(\rho^{3})
  \end{pmatrix},
\end{equation}
where $g^{\bullet\bullet}$ are given by:
\begin{equation}\label{inversemetriccomponents}
\begin{aligned}
  g^{\rho\rho} &= -q & g^{\rho v} &= -2 + \alpha v & g^{\rho y} &=
  -\frac{1}{2}\Upsilon + O(v)
  \\
  g^{vv} &= -4v + \beta v^{2} & g^{vy} &=
  -v\Upsilon+O(v^{2}) & g^{yy} &= -h^{-1} -O(v)
\end{aligned}
\end{equation}
Again all terms are smooth.

Cofactor expansion of equation~\eqref{eq:metric-at-bdry} scaled to the
frame $\pd[\rho]$, $\pd[v]$, $\pd[y]$ shows that the determinant of the
metric is
\begin{equation*}
  |g| = \rho^{-2(n+1)}|G| = \rho^{-2(n+1)}\left( (f^{2}-qv)|h| +O(\rho)\right) 
\end{equation*}
In particular,
\begin{align*}
  \frac{1}{2}\pd[\rho]\log|g| &= - (n+1)\rho^{-1} + O(1) \\
  \frac{1}{2}\pd[v]\log|g| &= O(1) \\
  \frac{1}{2}\pd[y]\log|g| &= O(1).
\end{align*}

\subsubsection{Induced metrics}
\label{sec:an-induced-metric}

In this section we describe induced metrics on the ``caps'' $C_{\pm}$
(the components of $\{ v > 0\}$ bounded by $S_{\pm}$)
and on the ``side'' $C_{0}$ ($\{ v < 0\}$).

We define the metric $K$ on $T^{*}X$ via the inclusion
$r^*: T^{*}X \hookrightarrow \Tbstar _{X}M$ (which is dual to the
restriction map $r: {}^{\bl}T_{X}M \to TX$).  As $\rho^{2}g$ is
a $\bl$-metric, we define for $\omega , \eta \in T^{*}(X)$ the
dual metric $K^{-1}$ by
\begin{equation*}
  K^{-1}(\omega, \eta) = - (\rho^{2}g)^{-1}(r^*\omega, r^*\eta)|_{\rho=0}.
\end{equation*}
Observe that $K^{-1}$ is the restriction of $-(\rho^{2}g )^{-1}$ to the
annihilator of $\rho\pd[\rho]$ (the ``b-normal'' vector field) at $\rho
= 0$.  

The components of the dual metric $K^{-1}$ are given in the frame
$\pd[v], \pd[y]$ by
\begin{align*}
  \left(
    \begin{array}{cc}
      K^{vv} & K^{vy} \\ K^{vy} & K^{yy}
    \end{array}
  \right) = \left(
    \begin{array}{cc}
      -g^{vv} & -g^{vy} \\ -g^{vy} & -g^{yy}
    \end{array}
  \right),
\end{align*}
where $g^{\bullet\bullet}$ are the components of the dual metric of $g$ in
the frame $\rho^{2}\pd[\rho]$, $\rho\pd[v]$, and $\rho\pd[y]$.

Because $\rho^{2}\pd[\rho]$ is time-like near $C_{\pm}$ and $K^{-1}$
is the restriction of $-(\rho^{2}g)^{-1}$ to the annihilator of
$\rho\pd[\rho]$, $K^{-1}$ is nondegenerate, and, in fact, Riemannian
in $C_{\pm}$.  In coordinates $(v,y)$, the metric $K$ on $TX$ is given
by
\begin{equation*}
  K = \frac{1}{4v} \left( 1 + O(v)\right) dv^{2} +
  \sum_{j=1}^{n-2}O(1) \left( dv \otimes dy_{j} + dy_j
    \otimes dv\right) + \sum_{i,j=1}^{n-1}K_{ij}dy_{i}\otimes dy_{j} .
\end{equation*}
It is easy to see from the above expression that the metric
\begin{equation*}
  k_{\pm} = \frac{1}{v}K|_{C_{\pm}}
\end{equation*}
is an asymptotically hyperbolic metric (in the sense of
Vasy~\cite{Vasy-Dyatlov:Microlocal-Kerr}) on $C_{\pm}$.  Setting $v =
x^{2}$ in this region ensures that $k$ is an asymptotically hyperbolic
metric (in the sense of Mazzeo--Melrose~\cite{Mazzeo-Melrose}) which
is \emph{even} in its boundary defining function (cf.\ the work of
Guillarmou~\cite{Guillarmou:AH}).

Similarly, because $\rho^{2}\pd[\rho]$ is space-like near $C_{0}$, 
$K^{-1}|_{C_{0}}$ is Lorentzian (with the ``mostly-plus'' convention), and
\begin{equation*}
  k_{0} = \frac{1}{v}K|_{C_{0}}
\end{equation*}
is an even \emph{asymptotically de Sitter} metric (with the ``mostly-minus''
convention, as $v< 0$ here) on $C_{0}$.  Indeed, if $v = -x^{2}$, then
the metric has the form used by Vasy~\cite{Vasy:asymp-dS}.  The
non-trapping assumption~(2) above implies that the metric satisfies
the conditions in Vasy's definition of an asymptotically de Sitter metric.

The $\rho$ components of the dual metric of $g$ are also related to
the components of the dual metric of $K$.  In the
$\rho^{2}\pd[\rho], \rho\pd[v], \rho\pd[y]$ frame for $g$ and the
$\pd[v], \pd[y]$ frame for $K$, we have
\begin{equation*}
  g^{\rho\rho} = -\frac{1}{v}\left( 4q K^{vv} + O(v^{2})\right), \quad
  g^{\rho v} = -\frac{1}{2v}\left( K^{vv} + O(v^{2})\right), \quad
  g^{\rho y} = -\frac{1}{2v}\left( K^{vy} + O(v^{2})\right).
\end{equation*}

As $K^{-1}$ is the lower-right block of $-g^{-1}$ and $g_{\rho\rho} =
v$, the volume forms of $g$ and $K$ (and hence the asymptotically
hyperbolic and asymptotically de Sitter metrics $k_{\pm}, k_{0}$) are
also related:
\begin{equation*}
  \sqrt{g} = \rho^{-(n+1)}\left( v^{1/2}\sqrt{|K|} + O(\rho)\right)
  = \rho^{-(n+1)}\left( v^{n/2}\sqrt{|k_{\pm,0}|} + O(\rho) \right).
\end{equation*}

\subsection{The form of the d'Alembertian}
\label{sec:form-dalembertian}

In this section we compute the form of the operator $\Box_{g}$ and its
normal operator $\widehat{N}(\rho^{-2}\Box_{g})$.

Putting the
calculations of the metric components and the volume form in
Section~\ref{sec:general-hypotheses} together, we compute the form of
$\Box_g$ near $\rho = 0$ (here we use $\sqrt{G} = \rho^{n+1}\sqrt{g}$
and recall that $g^{\bullet\bullet}$ are given by \eqref{inversemetriccomponents}):
\begin{align*}
  \label{eq:Box1}
  -\Box_{g} &= \rho^{2} \bigg[ \left( g^{\rho\rho}+ O(\rho)\right)
  (\rho \pd[\rho])^{2} + \left( g^{\rho v}+ O(\rho)\right)
  (\rho\pd[\rho])\pd[v]
  + \left( g^{\rho y}+ O(\rho)\right)(\rho\pd[\rho])\pd[y] \\
  &\quad\quad + (2-n) \left( \left( g^{\rho\rho}+O(\rho)
    \right)\rho\pd[\rho] + \left( g^{\rho v}+O(\rho)\right)\pd[v] +
    \left( g^{vy}+O(\rho)\right)\pd[y]\right) \\
  &\quad\quad + \frac{1}{\sqrt{G}}\pd[v]\left( \left( g^{\rho v} +
      O(\rho)\right)\sqrt{G} \rho\pd[\rho] +
    \left( g^{vv} + O(\rho) \right)\sqrt{G}\pd[v] + \left( g^{vy} + O(\rho)\right)\sqrt{G}\pd[y]\right) \\
  &\quad\quad + \frac{1}{\sqrt{G}}\pd[y] \left( \left( g^{\rho y}+
      O(\rho)\right)\sqrt{G}\rho \pd[\rho] + \left( g^{vy}+
      O(\rho)\right)\sqrt{G}\pd[v] + \left( g^{yy} + O(\rho)
    \right)\sqrt{G}\pd[y]\right)\bigg]
\end{align*}
Adopting now the notation of Vasy \cite{Vasy-Dyatlov:Microlocal-Kerr}, 
\begin{align*}
  -\widetilde{P}_{\sigma} &= -\widehat{N}(\rho^{-2}\Box_{g}) \\
  &= \frac{1}{\sqrt{G}}\left[ \pd[v] \left( g^{vv}\sqrt{G}\pd[v] +
        g^{vy}\sqrt{G}\pd[y]\right) + \pd[y]\left(
        g^{vy}\sqrt{G}\pd[v] + g^{yy}\sqrt{G}\pd[y]\right)\right] \\
    &\quad\quad +  g^{\rho v}\left( 2\imath\sigma + 2-n\right)
    \pd[v] +  g^{\rho y}\left( 2\imath\sigma +
        2-n\right) \pd[y] \\
    &\quad\quad + \imath\sigma\left[ \frac{1}{\sqrt{G}}\pd[v]\left( g^{\rho
          v}\sqrt{G}\right) + \frac{1}{\sqrt{G}}\pd[y]\left( g^{\rho
          y}\sqrt{G}\right) + g^{\rho\rho}\imath\sigma\right]
\end{align*}
In particular, near $v=0$,
\begin{align*}
  -\widetilde{P}_{\sigma} &= \left( -4v + O(v^{2})\right)\pd[v]^{2} +
  O(v)\pd[v]\pd[y] -\left( h^{ij} + O(v)\right)\pd[y_{i}]\pd[y_{j}] +
  O(1)\pd[y] \\
  &\quad\quad + 2\left( n - 4 - 2\imath\sigma  + O(v)\right)\pd[v] + q(\sigma),
\end{align*}
with $q$ a smooth function in $v$ and $y$ with values in quadratic
polynomials in $\sigma$.  

In our asymptotic expansions (and in the analysis of the radiation
field), it is more convenient to deal with
\begin{equation}\label{Psigmadef}
  P_{\sigma} \equiv \widehat{N}\left( \rho ^{-(n-2)/2}\rho^{-2}\Box_{g}\rho^{(n-2)/2}\right)
\end{equation}
than with $\widetilde{P}_{\sigma}$, in part to more directly correspond to
the setting of \cite{Vasy-Dyatlov:Microlocal-Kerr}.  To this end we
note simply that
\begin{equation}\label{Psigmatilde}
  P_{\sigma} = \widetilde{P}_{\sigma -\imath (n-2)/2}
\end{equation}
hence since the $\pd[v]^{2}$ and $\pd[v]$ terms of
$\widetilde{P}_{\sigma}$ may be written
\begin{equation}
  \label{eq:Psigma}
  -4\left( (v+O(v^{2}))D_{v}^{2} + \left( \frac{\imath}{2}(n-4-2\imath
    \sigma) + O(v)\right)D_{v}\right),
\end{equation}
we have
\begin{equation}\label{Psigma1}
  P_{\sigma} = -4\left( \left(v + O(v^{2})\right)D_{v}^{2} + \left(
      (\sigma - \imath) + O(v) \right)D_{v}\right) + O(1)\pd[y]^{2} +
  O(1)\pd[y] + O(v)\pd[v]\pd[y] + O(\sigma^{2}).
\end{equation}

\subsubsection{Relationship with the induced metrics}
\label{sec:relat-with-induc}

In the regions $C_{\pm}$ and $C_{0}$ of the boundary, $P_{\sigma}$ may
be written in terms of the metrics $k_{\pm}$ and $k_{0}$.  

We first work near $C_{\pm}$.  By an explicit computation, there is a
($\sigma$-dependent) vector field $\mathcal{X}(\sigma)$ tangent to
$v=0$ and a ($\sigma$-dependent) smooth potential $V(\sigma)\in
C^{\infty}(X)$ so that:
\begin{equation}
  \label{eq:relat-with-asymp-hyp}
  v^{\frac{1}{2}}v^{\frac{n}{4}+\frac{\imath\sigma}{2}}P_{\sigma}v^{-\frac{n}{4}
    - \frac{\imath\sigma}{2}}v^{\frac{1}{2}} = -
  \Lap _{k_{\pm}} + \left( \sigma^{2} + \frac{(n-2)^{2}}{4} \right) + v\mathcal{X}(\sigma) + vV(\sigma).
\end{equation}
(In terms of the variable $x$ given by $v = x^{2}$, the vector field
$\mathcal{X}$ is in fact a $0$-vector field in the sense of
Mazzeo--Melrose~\cite{Mazzeo-Melrose}.)  Moreover, if all $a_{j}$ and
$q$ vanish identically on $X$ (as is the case in Minkowski space and
the normally very short range perturbations of Minkowski space defined
in Section~\ref{sec:minkowski-space-1}) then $\mathcal{X} = 0$ and
$V=0$.
\begin{remark}
  The result of this computation should not be too surprising, as the
  entries of the inverse metric of $k$ agree up to a factor of $v$
  with a block of the inverse metric of $g$, accounting for the
  second-order terms. Moreover, it is easy to check that the operator
  on the left side is a b-differential operator on $X.$ The remainder
  of the computation requires checking only that the b-normal
  operators of the two sides agree.  A similar computation is carried
  out in \cite[Section 5]{Vasy-Dyatlov:Microlocal-Kerr}. 
\end{remark}

We now consider $C_{0}$.  The same calculation as above implies
that
\begin{equation}
  \label{eq:connection-with-dS}
  \smallabs{v}^{\frac{1}{2} + \frac{n}{4} +
    \frac{\imath\sigma}{2}}P_{\sigma}\smallabs{v}^{\frac{1}{2} -
    \frac{n}{4} - \frac{\imath\sigma}{2}} = \Box_{k_{0}} - \left(
    \sigma^{2} + \frac{(n-2)^{2}}{4}\right) + v \mathcal{X}(\sigma) + vV(\sigma),
\end{equation}
where $\mathcal{X}$ and $V$ are as above.  \emph{In particular, note that $P_{\sigma}$
is a hyperbolic operator on $C_{0}$ and an elliptic operator on $C_{\pm}$.}

\subsection{Location of radial points}
\label{sec:locat-radi-points}

We now study the flow associated to the Hamilton vector field of
$P_{\sigma}$.  In particular, we are interested in the radial points
of the vector field, i.e., those points in the characteristic set
where it is proportional to the fiber-radial vector field.  As
$P_{\sigma}$ is hyperbolic for $v< 0$ and elliptic for $v>0$, the only
possible radial points must occur when $v= 0$.  

As $0$ is not a critical point of $v$, we may take
\begin{equation*}
  \gamma dv + \eta \cdot dy
\end{equation*}
to be the canonical one-form on $T^{*}X.$ The principal symbol of
$P_{\sigma}$ is (employing summation convention) given by
\begin{equation*}
  \sigma (P_{\sigma}) = -(4v - \beta v^{2})\gamma^{2} - (2v\Upsilon +
  O(v^{2}))\gamma \eta - (h^{ij} + O(v))\eta_{i}\eta_{j}.
\end{equation*}
Letting $\hamvf$ denote the resulting Hamilton vector field on
$T^{*}X$, we have
\begin{equation}
  \label{eq:hamvf1}
  \frac{1}{2}\hamvf = (-4v\gamma + \beta v^{2}\gamma + v \eta \cdot
  \Upsilon) \pd[v] + (v\gamma \Upsilon_{j} + g^{y_{i}y_{j}}\eta_{i})\pd[y_{j}]
  + \bullet \pd[\gamma] + \bullet \pd[\eta],
\end{equation}
with the $\bullet$ terms homogeneous of degree 2 in the fiber
variables.  We now analyze the radial points of the vector field.  The
components in the base variables are given by
\begin{equation*}
  (-4v\gamma + \beta v^{2}\gamma + v\eta\cdot \Upsilon)\pd[v] +
  (v\gamma\Upsilon_{j} + g^{y_{i}y_{j}})\pd[y_{j}].
\end{equation*}
These coefficients must vanish at the radial set, which we have
already observed to lie over $v=0$.  In particular, we must have
\begin{equation*}
  g^{y_{i}y_{j}}\eta _{i} = 0
\end{equation*}
for all $j$.  As $g^{y_{i}y_{j}}$ is nondegenerate at $v=0$, we must
have $\eta =0$ on the radial set.  

We now easily verify that indeed the vector field at points
\begin{equation*}
  v=0, \ \eta = 0
\end{equation*}
is radial; hence these are in fact precisely the radial points.

\subsection{Structure near radial points}
\label{sec:struct-near-radi}

We now verify several of the hypotheses of \cite{Vasy-Dyatlov:Microlocal-Kerr}
near the radial points.  We have
established that the radial points occur at
\begin{equation}\label{radpoints}
  \Lambda^{\ep_1}_{\ep_2} \equiv \{ v=0, \ \eta = 0, \ep_2\gamma > 0\}\cap \pi^{-1}(S_{\ep_1})
  \subset T^{*}X,\ \ep_i=\pm;
\end{equation}
thus the $\pm$ in the superscript distinguishes ``past'' from
``future'' null infinity, while that in the subscript separates the
intersections with the two
components of the characteristic set.  We will
write $$\Lambda^{\pm}=\Lambda^{\pm}_+\cup\Lambda^{\pm}_-,\quad \Lambda_{\pm}=\Lambda_{\pm}^+\cup\Lambda_{\pm}^-.$$

We must now verify the following:
\begin{enumerate}
\item For a degree $-1$ defining function $\infbdf$ of $S^{*}X$
  inside the fiber-radial-compactification of $T^{*}X$, we
  have
  \begin{equation*}
    \infbdf \hamvf\infbdf \restricted{\Lambda_{\pm}} = \mp \beta_{0}, \
    \beta_{0} \in \CI(\Lambda_{\pm}), \ \beta _{0}> 0
  \end{equation*}
  (equation (2.3) of \cite{Vasy-Dyatlov:Microlocal-Kerr}).
\item There exists a non-negative homogeneous degree $0$ function
  $\rho_{0}$ vanishing quadratically and non-degenerately exactly at
  $\Lambda_{\pm}$ and a $\beta_{1} > 0$ such that
  \begin{equation*}
    \mp \infbdf \hamvf \rho_{0} - \beta _{1} \rho _{0} \geq 0 \text{
      modulo cubic terms vanishing at }\Lambda_{\pm}
  \end{equation*}
  (equation (2.4) of \cite{Vasy-Dyatlov:Microlocal-Kerr}).
\end{enumerate}

To deal with the first property, we remark that from
\eqref{eq:hamvf1}, we have
\begin{equation*}
  \frac{1}{2}\hamvf = (2\gamma ^{2} + O(v) + O(\eta)) \pd[\gamma] +
  (O(\eta^{2}) + O(v\eta) + O(v^{2}))\pd[\eta] + (-4v\gamma + \beta
  v^{2} \gamma + v\eta\cdot \Upsilon) \pd[v] + \bullet \pd[y]
\end{equation*}
where the big-Oh terms all have the homogeneities in $\gamma, \eta$
required to make the overall vector field homogeneous of degree $1$.
Near $\eta = 0$ we may employ the homogeneous coordinates
\begin{equation*}
  \infbdf = \frac{1}{\abs{\gamma}}, N = \frac{\eta}{\abs{\gamma}}
\end{equation*}
on the radial compactification of $T^{*}X$, hence we compute
that near $\Lambda_{\pm}$
\begin{equation}
\begin{aligned}
  \label{eq:compactifiedhamvf}
  \frac{1}{2} \hamvf &= \infbdf ^{-1} \big( (\mp 2 + O(v) + O(N))\infbdf
  \pd[\infbdf]\\
&\qquad + (\mp 2 N + O(v^{2}) + O(vN) + O(N^{2}))\pd[N]\\ &\qquad+ (\mp 4
  \pm \beta v + N\cdot \Upsilon ) v\pd[v] + \bullet \pd[y]\big),
\end{aligned}
\end{equation}
hence
\begin{equation*}
  \mp \infbdf \hamvf \infbdf \restricted{\Lambda_{\pm}} = 4,
\end{equation*}
i.e., the first property holds with $$\beta_{0} = 4.$$

To verify the second property, we take, in our compactified
coordinates,
\begin{equation*}
  \rho_{0} = v^{2} + N^{2}.
\end{equation*}
Applying~\eqref{eq:compactifiedhamvf} yields
\begin{equation*}
  \infbdf\hamvf(\rho_{0}) = \mp (16v^{2} + 8N^{2}) + \text{cubic terms
    in }(v,N),
\end{equation*}
i.e.,
\begin{equation*}
  \mp \infbdf \hamvf(\rho_{0}) - 8\rho_{0} = \text{cubic terms in }(v,N),
\end{equation*}
hence the second property is satisfied with $\beta_{1} = 8$.  

We may thus compute the subprincipal symbol of $\widetilde{P}_{\sigma}$
(and hence of $P_{\sigma}$) in terms of $\beta_{0}$.  Indeed, we compute
\begin{equation}\begin{aligned}
&-(2\imath)^{-1}(\tP_\sigma-\tP_\sigma^*)\\
&\qquad=
(2\imath)^{-1}\big((8+4(n-4+2\im\sigma)+O(v))\pa_v+O(1)\pa_y\big)+O(1)\\
&\qquad = -2\imath\big((n-2+2\im\sigma)+O(v)\big)\pa_v+O(1)\pa_y+O(1)
\end{aligned}\end{equation}
and consequently
\begin{equation}\label{subprincipal1}\begin{aligned}
\sigma \big((2\imath)^{-1}(\tP_\sigma-\tP_\sigma^*)\big)   \restricted{v=0,\,
  \eta=0}  &= \pm 4\left(-\frac{n-2}{2}-\im\sigma\right)|\gamma|\\
&=\pm\beta_0\left(-\frac{n-2}{2}-\im\sigma\right)|\gamma|.
\end{aligned}\end{equation}
Note that, even apart from the shift by $(n-2)/2$,
the sign of $(2\imath)^{-1}(\tP_\sigma-\tP_\sigma^*)$ is
switched as compared to \cite{Vasy-Dyatlov:Microlocal-Kerr} (where $\pm
\beta_0\im\sigma|\gamma|$ was used with the present notation). Switching the roles of $\tP_\sigma$ and $\tP_\sigma^*$ reverses
this sign, and thus what we do here corresponds to what was discussed
in \cite{Vasy-Dyatlov:Microlocal-Kerr} for the adjoint operator in the context of the general
theory, though this reversal was pointed out there already in the Minkowski
context in Section~5 of \cite{Vasy-Dyatlov:Microlocal-Kerr}.

Returning to the operator $P_\sigma$ itself, we compute for later
reference that by \eqref{subprincipal1} and \eqref{Psigmatilde},
\begin{align}\label{hatbeta}
  \hat\beta^\pm(\sigma)&\equiv\pm\frac{\infbdf}{2\imath\beta_0}\sigma_{1}(P_\sigma-P_\sigma^*)|_{\Lambda^\pm}\\
  &= \Big( -\frac{(n-2)}2-\Im (\sigma-\imath(n-2)/2)\Big) =-\Im\sigma. \notag
\end{align}

\subsection{b-radial points}\label{sec:b-radial-points}
It is also useful to compute the full b-structure of the
radial set of
$$
L=\rho^{-(n-2)/2}\rho^{-2}\Box_g\rho^{(n-2)/2}\in\Diffb^2(M).
$$
in $\Tb^* M$. Note that the powers are chosen here so that $L$ is
formally self-adjoint with respect to the b-density
$$
\rho^n\,|dg|.
$$
The b-principal symbol of $L$ is the ``same'' as the sc-principal symbol
of $\Box_g$ under the identification of $\Tb^*M$ and $\Tsc^* M$, namely
\begin{equation}\begin{aligned}\label{eq:b-princ-symbol}
  \lambda=\sigmab (L) &= g^{\rho\rho}\xi^{2} - (4v - \beta v^{2} +
  O(\rho v) + O(\rho^{2})) \gamma^{2} - 2(2-\alpha v+ O(\rho))\xi
  \gamma \\
  &\quad + 2g^{\rho y}\cdot \eta \xi + \big(2 v \Upsilon + O(\rho)
  \big) \cdot \eta \gamma + g^{y_{i}y_{j}}\eta_{i}\eta_{j},
\end{aligned}\end{equation}
where we write b-covectors as
\begin{equation}\label{b-covectors}
\xi\frac{d\rho}{\rho}+\gamma\,dv+\eta\,dy.
\end{equation}
The b-Hamilton vector field of a symbol $\lambda$ is
\begin{equation}\label{eq:b-ham-vf-gen}
(\pa_\xi\lambda) (\rho\pa_\rho)+(\pa_\gamma\lambda)
\pa_v+(\pa_\eta\lambda)\pa_y-(\rho\pa_\rho\lambda)\pa_\xi-(\pa_v\lambda)\pa_\gamma-(\pa_y\lambda)\pa_\eta,
\end{equation}
so in our case we obtain
\begin{equation}\label{bHamvf1}\begin{aligned}
\sH_{\lambda}=&\big(2g^{\rho\rho}\xi+2g^{\rho y}\eta-2\gamma(2-\alpha
v+O(\rho))\big)(\rho\pa_\rho)\\
&-2\big((4v-\beta v^2+O(\rho
v)+O(\rho^2))\gamma\\
&\qquad\qquad\qquad+(2-\alpha
v+O(\rho))\xi+(v\Upsilon+O(\rho))\eta\big)\pa_v\\
&+2\big(g^{\rho y}\xi+(v\Upsilon+O(\rho))\gamma+g^{y_iy_j}\eta_j\big)\pa_y\\
&-(\rho\pa_\rho \lambda)\pa_\xi-(\pa_v \lambda) \pa_\gamma-(\pa_y \lambda) \pa_\eta.
\end{aligned}\end{equation}

We now investigate when this vector field has radial points when
restricted to $\rho=0.$  The symbol $\lambda$ is a
nondegenerate (Lorentzian) metric on the fibers of $\Tbstar M,$ hence
the projection to the base of $\sH_\lambda$ must be a nonvanishing
b-vector field; for $\pi(\sH_\lambda)$ to vanish over the boundary,
then, it must lie in the span of $\rho\pa_\rho.$  On the other hand,
letting $g_b$ denote our induced b-metric given by (the dual of)
$\lambda$ and letting $p
\in \Tbstar M$ we have $g_b(\pi\sH_\lambda|_{p},
\pi\sH_\lambda|_{p})=\lambda(p)=0$ assuming $p$ lies in the characteristic
set.  Thus, at a radial point over $\rho=0$, $\rho\pa_\rho$ must be a null vector
field, hence we must have $v=0.$

We further see by examining the coefficients of $\pa_v,$ $\pa_y$ that the radial set $\cR$
within $\rho=0$ is exactly $v=0$, $\eta=0$, $\xi=0$. Further, there
are no radial points in $\rho>0$, since the metric is a standard
Lorentzian metric there (and there is no distinction between b-metrics
and standard metrics in the interior).  Now, on the fiber
compactification of $\Tb^*M$ near $\cR$ we can use local coordinates,
\begin{equation}\label{projcoords}
\nu=\frac{1}{\gamma},\ \hat\xi=\frac{\xi}{\gamma},\ \hat\eta=\frac{\eta}{\gamma},
\end{equation}
to obtain the linearization of $\sH_{\lambda}$ at $\cR.$ That is,
$\nu\sH_{\lambda}\in\Vb(\overline{\Tb^*}M)$, i.e.\ is tangent to both
$\rho=0$, defining $\pa M$, and $\nu=0$, defining fiber infinity,
vanishes at $\pa\cR$ (fiber infinity of the radial set), thus maps the ideal $\cI$ of $\CI$ functions
vanishing at a point
$q\in\pa\cR$ to themselves, and thus $\cI^2$ to $\cI^2$, so it acts on
$\cI/\cI^2\cong T^*_q\overline{\Tb^*}M$. In computing this, terms of
$\nu\sH_{\lambda}$ which vanish quadratically at $\pa\cR$ can be
neglected; modulo these we have
\begin{equation*}\begin{aligned}
\nu\sH_\lambda=-4\rho\pa_\rho+(-8v-4\hat\xi)\pa_v&+2(g^{\rho
  y_i}\hat\xi+v\Upsilon+\rho c+g^{y_i y_j}\hat\eta_j)\pa_{y_i}\\
&-4(\nu\pa_\nu+\hat\xi\pa_{\hat\xi}+\hat\eta\pa_{\hat\eta})+\cI^2\Vf(\overline{\Tb^*}M),
\end{aligned}\end{equation*}
with $c$ smooth. Correspondingly, the eigenvectors and eigenvalues of
$\nu\sH_\lambda$ are
\begin{equation}\label{eigenvectors}\begin{aligned}
&dv+d\hat\xi,\ \text{with eigenvalue}\ -8,\\
&d\rho,d\nu,d\hat\xi,d\hat\eta,\ \text{with eigenvalue}\ -4,\\
&2\,dy+g^{\rho y}\,d\hat\xi+\Upsilon\,dv-c\,d\rho+g^{y_i y_j}\,d\hat\eta_j,\ \text{with eigenvalue}\ 0.
\end{aligned}\end{equation}

\subsection{Semiclassical symbol and flow}\label{sec:semicl}
In this section we record the relationship of the computations performed
above in in the b-cotangent bundle to the semiclassical results on
$P_\sigma$ that we will require below.  Fortunately, these computations
are nearly identical: as $P_\sigma$ is obtained from $L$ by Mellin
transform, if we let $\sigma_\hbar$ denote the semiclassical principal
symbol of an operator with parameter $\sigma,$ where $\pm \Re \sigma$ is the
semiclassical parameter, then we have simply
$$
\sigma_\hbar(\abs{\sigma}^{-2}P_\sigma)\cong \sigma_{\bl}(L)|_{\xi=\pm 1},\
\Re \sigma \to \pm \infty.
$$
The computation of the semiclassical Hamilton flow is similarly
simple: the vector field in \eqref{bHamvf1} is tangent to $\Tb^*_{\pa
  M}M$, with a vanishing $\pa_\xi$ component at $\rho=0$; thus the
semiclassical flow associated to $P_\sigma$ is given precisely by
\eqref{bHamvf1}, restricted to $\rho=0$ and with $\xi=\pm 1.$

We let $\Sigma_{\semi}$ denote the semi-classical characteristic
set, where $\sigma_\hbar(\abs{\sigma}^{-2}P_\sigma)=0,$ and let
$\Sigma_{\semi,\pm}$ denote its two components; by the above
discussion, these are simply the same as the characteristic sets of
the rescaling of $\Box$ viewed as a b-operator.

For later use, we also recall the \emph{semi-classical Sobolev} spaces
appropriate to our problem.  These spaces are denoted
$H^s_{\smallabs{\sigma}^{-1}}$ and when $s$ is a positive integer they
are given in local coordinates by
$$
u \in H^s_{\smallabs{\sigma}^{-1}} \Longleftrightarrow
\smallabs{\sigma}^{-\smallabs{\alpha}} D^\alpha u\in L^2,\ \text{ for
  all } \smallabs{\alpha} \leq s.
$$
The definition can be extended to non-integer $s$ via interpolation
and duality, or else by using the semiclassical pseudodifferential
calculus with parameter $h=\smallabs{\sigma}^{-1}.$ We refer the reader to \S2.8 of
\cite{Vasy-Dyatlov:Microlocal-Kerr} for details.

\subsection{The radiation field blow-up}
\label{sec:blowup}

Having described our geometric set-up, we now digress slightly to
revisit the ``usual'' construction of the radiation field in the
context at hand.  This section is not necessary in the logical
development of the paper but rather serves to situation our results
in the context of prior theorems.

In particular, although the existence of the radiation field for tempered solutions
of $\Box_{g} w = f\in \dCI(M)$ with appropriate support properties is
a consequence of our main theorem, in this section we recall the definition
of the radiation field for metrics of the form in
Section~\ref{sec:general-hypotheses}.  In order to do this we will
also need to assume an additional
support condition on the solution $w$ analogous to the one satisfied
by the forward fundamental solution in more familiar contexts.  In
particular, then, assume that $g$ is a non-trapping Lorentzian
scattering metric as described above, and further assume that
\begin{quote}The function $w$ solves $\Box_{g}w = f\in
   C^{\infty}_{c}(M^{\circ})$ and there is an $s_{0}$ so that near
   $S_{+}$, $w$ vanishes identically for $v/\rho \geq s_{0}$.
 \end{quote}

We now \emph{blow up} $S=\{ v = \rho = 0\}$ by replacing it with its inward
pointing spherical normal bundle.  (The reader may wish to
  consult \cite{Melrose:APS} for more details on the blow-up
  construction than we give here.) 
 This process replaces $M$ with a
new manifold $\Mbar = [ M ; S]$ on which polar
coordinates around the submanifold are smooth, and depends only on $S$
(not the actual functions $v$ and $\rho$). The blow-up comes
equipped with a natural blow-down map $\Mbar \to M$ which is a
diffeomorphism on the interior.  $\Mbar$ is a manifold with corners
with two boundary hypersurfaces: $\operatorname{bf}$, the closure of
the lift of $X \setminus S$ to $\Mbar$; and
$\operatorname{ff}$, the lift of $S$ to $\Mbar$. Further, the fibers
of $\operatorname{ff}$ over the base, $S$, are diffeomorphic to
intervals, and indeed, the interior of the fibers is naturally an affine
space (i.e.\ these interiors have $\RR$ acting by translations, but there is no
natural origin). Figure~\ref{fig:blow-up} depicts this blow-up construction.

Given $v$ and $\rho,$ the fibers of the interior of
$\operatorname{ff}$ in $[M;S]$ can be identified with $\RR$, via the
coordinate $s=v/\rho$. In particular, $\pd[s]$ is a well-defined vector
field on the fibers.

We define ``polar coordinates''
\begin{equation*}
  R = \left( v^{2} + \rho ^{2}\right) ^{1/2} \in [0, \infty), \quad \Theta =
  \frac{\left( \rho, v\right)}{r} \in \sphere ^{1}_{+},
\end{equation*}
which are smooth on $\Mbar$.  Near the interior of
$\operatorname{ff}$, we use the projective coordinates $\rho, s
= v/\rho$ as well as local coordinates $y$ on $S$. 
In these coordinates, a simple computation shows
that the unbounded terms of $\rho^{2}g$ cancel near $\rho = 0$ and
hence $\rho^{2}g$ is a smooth Lorentzian metric in a neighborhood of
the interior of $\operatorname{ff}$ (i.e., down to $\rho = 0$).

Given a solution $w(\rho, v, y)$ of $\Box_{g}w = f$ with $f$ smooth
and compactly supported, we define the function
\begin{equation*}
  u (\rho, s, y) = \rho ^{-\frac{n-2}{2}}w (\rho, \rho s, y).
\end{equation*}
The wave operators for the metrics $g$ and $\rho^{2}g$ are related by
the somewhat remarkable identity
\begin{align*}
  \rho^{\frac{2-n}{2}}\Box_{g}w &= \rho^{\frac{2-n}{2}}\Box_{g}\left(
    \rho^{\frac{n-2}{2}}u\right) = \rho^{2}\Box_{\rho^{2}g}u -
  \left(\rho
    ^{\frac{n+2}{2}}\Box_{\rho^{2}g}\rho^{\frac{2-n}{2}}\right)u \\
  &= \rho^{2}\Box_{\rho^{2}g}u - \rho^{2}\gamma u;
\end{align*}
we refer the reader to \cite{Friedlander:1980} for the details of this
computation.  Note that $\gamma$ is smooth on $M$ because $\rho$ is.
Moreover, $\rho^{2} g$ is a nondegenerate metric near the interior of
$\operatorname{ff}$, and so $\Box_{\rho^{2}g} - \gamma$ is a
nondegenerate hyperbolic operator near $\operatorname{ff}$.  This
calculation thus shows that if $w$ is a solution of $\Box_g w = f$
with smooth compactly supported $f$, vanishing identically for $s \leq
s_{0}$, then the argument of Friedlander~\cite[Section
1]{Friedlander:1980} shows that $w$ may be smoothly extended across
$\operatorname{ff}$.  In particular, $w$ and its derivatives may be
restricted to $\operatorname{ff}$.  Note that the condition on the
support of $w$ is analogous to the support condition satisfied by
forward solutions of the inhomogeneous equation.  (The argument
applies equally well to solutions of the homogeneous initial value
problem with the same support property on globally hyperbolic
spacetimes of this form.)

Thus,
  if $w$ is a solution of $\Box_{g}w = f$ satisfying the above support
  property, with $f$ smooth and compactly supported, we may define the
  (forward) \emph{radiation field} of $w$ by
  \begin{equation*}
    \mathcal{R}_{+}[w](s, y) = \pd[s] u(0, s, y).
  \end{equation*}
  This agrees with Friedlander's original construction.  We will later
  show that we may make this definition even without the hypothesis on
  the support of $w$ in the $s$ variable.

\begin{remark}
  \label{rem:radiation-field-blow}
  Note that the smooth expansion of $w$ across $\operatorname{ff}$
  implies that it does not have singularities at $s = 0$.
\end{remark}

\section{Propagation of b-regularity}
\label{sec:b-regularity}

In this section we prove an initial conormal estimate for tempered
solutions $w$ of $\Box_g w = f\in \dCI(M)$ vanishing near
$\overline{C}_-$.  This estimate is used to begin the iterative scheme in
Section~\ref{sec:an-asympt-expans}.  Our goal is conormal regularity at $\Lambda^+.$

The basic background in this section is the propagation of
b-regularity away from radial points (see, e.g., \cite{Vasy:corners}),
which we briefly recall here.  Let $L\in\Psib^{s,r}(M)$, and let
$\Sigma\subset\Sb^*M$ denote the characteristic set of $L$, $\lambda$
denote the principal symbol of $L$ in $\Psib^{s,r}(M)$.
\begin{proposition}\label{proposition:bpropagation} Suppose $w \in \Hb^{-\infty,l}(M).$
Then
\begin{enumerate}\item
Elliptic regularity holds away from $\Sigma$, i.e.,
$$
\WFb^{m,l}(w)\subset\WFb^{m-s,l-r}(Lw)\cup\Sigma,
$$
\item In $\Sigma$, $\WFb^{m,l}(w)\setminus\WFb^{m-s+1,r-l}(Lw)$ is a
union of maximally extended bicharacteristics, i.e., integral curves
of $\sH_\lambda$.
\end{enumerate}
\end{proposition}
Note that the order in $\WFb^{m-s+1,r-l}(Lw)$ is shifted by $1$
relative to the elliptic estimates, corresponding to the usual
hyperbolic loss.  This arises naturally in the positive commutator
estimates used to prove such hyperbolic estimates: commutators in
$\Psib(M)$ are one order lower than products in the differentiability
sense (the first index),  but not in the decay order (the second
index); hence the change in the first order relative to elliptic estimates
but not in the second.

We now turn to the radial set $\cR$, where
Proposition~\ref{proposition:bpropagation} does not yield any
interesting statements, and more refined arguments are needed.
\begin{definition}
Let $\cM\subset \Psib^1(M)$ denote the $\Psib^0(M)$-module of pseudodifferential
operators with principal symbol vanishing on the radial set
$\cR=\{\rho=0,\ v=0,\ \xi=0,\ \eta=0\}$.
\end{definition}
Note that a set of generators
for $\cM$ over $\Psib^0(M)$ is $\rho\pd[\rho],\rho
    \pd[v],v\pd[v],\pd[y]$ (with symbols
    $\xi,\rho\gamma,v\gamma,\eta$; $\gamma$ enters to convert $\rho$ and
    $v$ to first order operators) and $\Id$.  (Recall that these fiber
    variables are defined by the canonical one-form \eqref{b-covectors}.)

\begin{lemma}
The module $\cM$ is closed under commutators.
\end{lemma}
\begin{proof}
  While this can be checked directly from \eqref{eq:b-ham-vf-gen}, a
  more conceptual proof is as follows.

  We observe from the formula \eqref{eq:b-ham-vf-gen} that whenever
  $f\in\CI(\Tbstar M)$ (or just defined on an open subset, such as
  $\Tbstar M\setminus o$), the b-Hamilton vector field $\sH_f$ is
  tangent to the submanifold $\{\rho=\xi=0\} \subset \Tbstar M$. This
  submanifold is the image of $T^*_XM$ in $\Tbstar M$ via the
  canonical map dual to the inclusion of $\Vb(M)$ in $\Vf(M)$, and is
  thus canonically identified with $T^*X$; we tacitly use this
  identification from now on. Correspondingly, it is a symplectic
  manifold with symplectic form $d\gamma \wedge dv + d\eta\wedge dy$,
  and further, the restriction of the b-Hamilton vector field $\sH_f$
  of $f\in\CI(\Tbstar M)$ as above is equal to the standard Hamilton
  vector field of $f|_{T^*X}\in\CI(T^*X)$.  Since $\cR$ is a
  \emph{Lagrangian} submanifold of $T^*X,$ if $f|_{T^*X}$ vanishes on
  $\cR$, the Hamilton vector field of $f$ is tangent to $\cR.$
  Correspondingly the set of $\CI$ functions vanishing on $\cR$ is
  closed under Poisson brackets; taking into account that the Poisson
  bracket of homogeneous degree one functions on $\Tbstar M\setminus
  o$ is also such and that the principal symbol of the commutator of
  two b-pseudodifferential operators is given by Poisson brackets, we
  immediately conclude that $\cM$ is closed under commutators.
\end{proof}

\begin{proposition}\label{prop:radial-b-estimate}
Let $L=\rho^{-(n-2)/2}\rho^{-2}\Box_g\rho^{(n-2)/2}\in\Diffb^2(M)$.
If $w\in \Hb^{-\infty,l}(M)$ for some $l$, $L w\in
\Hb^{m-1,l}$ and $w\in \Hb^{m,l}$ on a punctured neighborhood $U\setminus\pa\cR$ of
$\pa\cR$ in $\Sb^*M$ (i.e.\ $\WFb^{m,l}(w)\cap
(U\setminus\pa\cR)=\emptyset$)
then for $m'\leq m$ with $m'+l<1/2$,
$w\in \Hb^{m',l}(M)$ at $\pa\cR$ (i.e.\
$\WFb^{m',l}(w)\cap\pa\cR=\emptyset$),
and for $N\in\NN$ with $m'+N\leq m$ and for
$A\in\cM^N$, $Aw$ is in $\Hb^{m',l}(M)$ at $\pa\cR$ (i.e.\
$\WFb^{m',l}(Aw)\cap\pa\cR=\emptyset$).
\end{proposition}

\begin{remark}
In the situation that we care about, $\cR=\cR_+\cup\cR_-$ splits into two components
(``future'' and ``past'') and we note that the proof in fact shows
that the result holds at each component separately.

This result is analogous to \cite{Haber-Vasy:Radial}, except $\rho=0$
produces an extra boundary (so we are in codimension 2), and $\cR$ is
not Lagrangian ($\Tb^*M$ is not symplectic at the boundary).  The
relevant input of the Lagrangian nature in
\cite{Haber-Vasy:Radial} is the eigenvectors and eigenvalues of
the linearization, hence much the same proof goes through. This is also analogous to the ``easy'' part,
Section~11, of \cite{Melrose-Vasy-Wunsch:Propagation}, describing the
propagation of edge singularities, except here we have a source/sink
rather than a saddle point, and thus the treatment is simpler.
\end{remark}

\begin{proof}
First we ignore the module. We will inductively show that 
$\WFb^{\tilde m,l}(w)\cap\pa\cR=\emptyset$ assuming that we
already has shown $\WFb^{m'',l}(w)\cap\pa\cR=\emptyset$ with
$m''=\tilde m-1/2$. We may start with $\tilde m=\min(m_0+1/2,m')$,
and increasing $\tilde m$ by $\leq 1/2$, we reach $\tilde m=m'$ in
finitely many steps.

Thus, consider $A\in\Psib^{s,r}(M)=\rho^{-r}\Psib^s(M).$  Then
$$
\imath [L,A]\in\Psib^{s+1,r}(M),\ \sigma_\bl(\imath [L,A])=\sH_\lambda
a,\ a=\sigma_{\bl,r,s}(A).
$$
We choose
$$
a=\rho^{-r}\nu^{-s}\phi^2,
$$
where $\phi\geq 0$, $\phi\equiv 1$ near $\cR$, supported in $U$
($\supp \phi$ will be further constrained below).  By
\eqref{eigenvectors}, $\nu\sH_\lambda a=(4(r+s)+c)a+e$, where $c$
vanishes at $v=0$, and $e$ is supported in $\supp d\phi$. We take
$r+s<0$, and we choose the support of $\phi$ so that $|c|<|r+s|$ on
the support of $\phi$. Note
that $r+s<0$ means $\nu\sH_\lambda a$ necessarily has negative sign
at least in some place on $\supp d\phi$, since $\phi$ has to increase
along the flow as it approaches $\cR$. Then we have $\nu\sH_\lambda
a=-b^2+e$, with $b$ elliptic near $\cR$. Then with
$B\in\Psib^{(s+1)/2,r/2}(M)$ with principal symbol $b$ and
$\WFb'(B)\subset \supp b\cap \Sb^* M$ (so for instance $B$ can be a
quantization of $b$), and similarly with $E\in\Psib^{s+1,r}(M)$,
$$
\imath[L,A]=-B^*B+E+F,\ F\in\Psib^{s,r}(M).
$$
This gives an estimate
\begin{equation}\label{foo.Feb21}
\|Bw\|^2\leq |\langle Ew,w\rangle|+|\langle Fw,w\rangle|+2|\langle Lw,Aw\rangle|
\end{equation}
when $w$ is a priori sufficiently regular. Given $\tilde m,l$, we now take
$s=2\tilde m-1$, $r=2l$, so $s+r<0$ indeed. Note that $F$ has order
$\leq 2m''$, so the inductive assumption gives a bound for $|\langle
Fw,w\rangle|$.
A standard regularization argument can be used to complete the proof
by allowing us to apply \eqref{foo.Feb21} to any $w$ for which the
right-hand-side is a priori finite:
for instance one can use a regularizer
$\psi_\ep(\nu)=(1+\ep\nu^{-1})^{-1}=\frac{\nu}{\nu+\ep}$, $\ep>0$,
which is in $S^{-1}$ for $\ep>0$ and is uniformly bounded in $S^0$ for
$\ep\in (0,1]$. Thus, one lets
$$
a_\ep=a\psi_\ep(\nu)^2;
$$
then $\nu\sH_\lambda
\psi_\ep=\ep\nu^{-2}\psi_\ep^2(\nu\sH_\lambda\nu)$ shows that the
contribution of the regularizer to the principal symbol of the
commutator is the negative of a square, provided again that $\phi$ has
sufficiently small support, i.e.\ adds another ``good term'' beside
$-b^2$. One can drop the corresponding term in the inequality given by
quantized version,
$$
\|B_\ep w\|^2\leq |\langle E_\ep w,w\rangle|+|\langle F_\ep
w,w\rangle|+2|\langle Lw,A_\ep w\rangle|,
$$
where the calculation (involving the pairing) now makes sense for $\ep>0$.
Now letting $\ep\to 0$ the right hand side remains bounded, while
$B_\ep \to B$ strongly in $\cL(L^2_\bl(M))$, so one concludes $Bw\in
L^2_\bl(M)$ and obtains the desired inequality.  This completes the proof in the case when $N=0,$
i.e., when we have included no factors from the module $\cM$ in the
test operator.

In the general case $N\geq 0$ one employs the methods developed by Hassell,
Melrose and Vasy \cite{Hassell-Melrose-Vasy:Spectral,
  Hassell-Melrose-Vasy:Microlocal}, adapted to a similar, but
different (edge), setting by Melrose, Vasy and Wunsch in
the appendix of \cite{Melrose-Vasy-Wunsch:Propagation}.
For this purpose one uses generators of the module, denoted by
$G_0=\Id,G_1,\ldots,G_{n},G_{n+1}=\Lambda L$, where
$\Lambda\in\Psib^{-1}$ is elliptic near $\cR$ and $G_1,\dots, G_n \in
\Psib^1(M).$  A sufficient condition
for these methods is that for $i=1,\ldots,n$,
\begin{equation}
\imath\Lambda [G_i,L]=\sum_j C_{ij} G_j,
\end{equation}
where \begin{equation}\label{vanishingsymbol}\sigma_{\bl,0,0}(C_{ij})|_{\cR}=0.\end{equation}
In our case this sufficient condition is satisfied by choosing $dg_i$
to be an eigenvector of $\nu\sH_{\lambda}$ at $\cR$, with eigenvalue
$-4$ (cf.\ \eqref{eigenvectors}), where $G_i$ has principal symbol
$\nu^{-1}g_i$. For instance, we may take the $g_i$ to be
$\rho,\hat\xi, \hat \eta$ since these, together with $\lambda \nu^2$
cut out $\cR$ in the cosphere bundle.  Since $d\nu$ and $dg_i$ have
equal eigenvalues then, the conclusion for $C_{ij}$ follows.  (We note
that strictly speaking, because $\hat \eta$ is not globally defined,
we must include additional generators to account for different
coordinate charts in the tangential variables.  Including additional
generators does not cause any problem, as we do not require the
generating set to be independent.)  

We thus prove
iterative regularity under $\cM$ inductively in the power of the
module as follows: we repeat the previous commutator argument, but with the
commutant $A$ replaced by 
$$
\Op(\sqrt{a})^* (G^\alpha)^* (G^\alpha) \Op(\sqrt{a})
$$
where $G^\alpha=G_1^{\alpha_1}\dots G_{n+1}^{\alpha_{n+1}}$ denotes a product of powers of
the generators of $\cM,$ hence $G^\alpha \in \cM^{\smallabs{\alpha}}.$
Considering all of these commutators at once, as $G^\alpha$ runs over a
basis of $\cM^N/\cM^{N-1},$ we then follow the same argument as used when
$N=0$ but now with systems of operators, taking values in $\CC^d$ with $d=\dim
\cM^N/\cM^{N-1}.$  The main term in the commutator, arising from the
commutators $[L,\Op(\sqrt{a})],$ is diagonal and positive, just as before
(again, because the factor $4(r+s)+c$ is negative).
Moreover the condition \eqref{vanishingsymbol} permits us
to absorb into this positive term those new terms that arise from commutators of $L$
with $G^\alpha$ and that have the maximum number of module factors.  Thus
we are in the end able to estimate the terms $\norm{B G^\alpha w}^2$
(with $B$ as before) where $\smallabs{\alpha}=N$ by terms microsupported
away from $\cR$ and by terms involving $G^\beta w$ with
$\smallabs{\beta}\leq N-1,$ thus proving the result inductively.
\end{proof}

Putting Propositions~\ref{prop:radial-b-estimate} and
\ref{proposition:bpropagation} together then yields the following:

\begin{corollary}\label{cor:radial-b-estimate}
Let $L=\rho^{-(n-2)/2}\rho^{-2}\Box_g\rho^{(n-2)/2}\in\Diffb^2(M)$, and
let $\pi:\Tb^*M\to M$ be the projection.
Suppose $w\in \Hb^{-\infty,l}(M)$ for some $l,$ $L w\in
\Hb^{m-1,l}(M)$. Suppose $\cU$ is a neighborhood of $\pi(\pa\cR)$
and that all bicharacteristics (in $\Sigma$) of $L$
that enter $\cU$, other than those in $\cR$, possess a point disjoint from $\WFb^{m,l}(w)$.
Then for $m'\leq m$ with $m'+l<1/2$, $w$ is in $\Hb^{m',l}(M)$ on $\cU$
and for $N\in\NN$ with $m'+N\leq m$ and for
$A\in\cM^N$, $Aw$ is in $\Hb^{m',l}(M)$ on $\cU$.
\end{corollary}

Note that the hypotheses of the corollary at the future radial set
hold automatically if $L$ is non-trapping, i.e.\ all bicharacteristics
tend to the future and past radial sets in the two directions of flow,
and if $w$ vanishes near $S_{-}$.  

\begin{remark}
  \label{rem:blow-up-conormal}
Corollary~\ref{cor:radial-b-estimate} implies that when $L w$ is
  in $\Hb^{\infty,l}(M)$ (so in particular if $Lw\in\dCI(M)$) and $w$ vanishes near $\overline{C_-}$ then
$w$ is in fact
  \emph{conormal} to the front face of the blow-up defined in
  Section~\ref{sec:blowup} since we obtain $\Hb^{\infty,l}$
  regularity sufficiently far along all bicharacteristics (indeed, the solution
  vanishes in a neighborhood of the boundary there by hypothesis).  In particular, this implies that the
  coefficients in the asymptotic expansion of
  Theorem~\ref{thm:main-theorem} may be taken to be smooth.
\end{remark}

\section{The mapping properties of $P_{\sigma}$ }
\label{sec:results-citevasy}

Having verified that the operator $P_{\sigma}$ satisfies many of the
hypotheses of the theorem of Vasy~\cite{Vasy-Dyatlov:Microlocal-Kerr},
we now show that $P_{\sigma}$ is Fredholm on appropriate function
spaces.  In this section we modify the argument of
\cite{Vasy-Dyatlov:Microlocal-Kerr} to our current setting.

Recall that under our global assumptions, the characteristic set of
$P_\sigma$ in $S^*X$ has two parts $\Sigma_\pm$ (each of
which is a union of connected components) such that the integral
curves of the Hamilton flow in $\Sigma_{\pm}$ tend to $S_{\pm}$ as the
parameter tends to $+\infty$.  Writing the radial
sets  at future, resp.\ past,
infinity as $\Lambda^+$, resp.\ $\Lambda^-$ (and the components of each as
$\Lambda^\pm_\pm=\Lambda^\pm\cap\Sigma_\pm$), one is interested in the
following two kinds of Fredholm problems, in which one requires a
relatively high degree of regularity at $\Lambda^+$, resp.\
$\Lambda^-$, but allows very low regularity at the other radial set,
$\Lambda^-$, resp.\ $\Lambda^+$.

To make this into a Fredholm problem it is convenient to introduce
variable order Sobolev spaces and variable order pseudodifferential
operators. This was originally done by Visik, Eskin
\cite{Visik-Eskin:Sobolev}, Unterberger \cite{Unterberger:Resolution}
and Duistermaat \cite{Duistermaat:Carleman}, and we recall this theory
in Appendix~\ref{section:variableorder}.  More recently,
Faure--Sj{\"o}strand~\cite{Faure-Sjostrand} used variable-order
Sobolev spaces in a manner similar to ours in their work on Ruelle
resonances for Anosov flows.  The main result that we
use is Proposition~\ref{proposition:variablepropagation}, which shows
that standard propagation of singularities arguments along forward
null-bicharacteristics hold with respect to the spaces $H^s$ with $s
\in \CI(S^*X)$ defining the variable order, provided $s$ is
non-decreasing along the Hamilton flow.
\begin{remark}
  We recall that in \cite{Vasy-Dyatlov:Microlocal-Kerr} such issues
  were avoided by using complex absorption arranged so that the
  resulting operator is elliptic at one of the radial sets, say
  $\Lambda^-$, but is unchanged near $\Lambda^+$. Thus, each
  bicharacteristic enters the complex absorption region in either the
  forward or backward direction; in this region the operator becomes
  elliptic due to the imaginary part of its principal symbol, hence
  only $\Lambda^+$ acts as a radial set for the operator with complex
  absorption added, and one could use standard Sobolev spaces as one
  did not have to deal with different regularity thresholds at
  $\Lambda^+$ and $\Lambda^-$.
\end{remark}

Now we recall, as computed in \eqref{hatbeta}, that the quantity
$$
\hat\beta^\pm(\sigma)=\pm\frac{\infbdf}{2\imath\beta_0}\sigma_{1}(P_\sigma-P_\sigma^*)|_{\Lambda^\pm}
$$
is given the ``constant'' value $-\Im \sigma:$ it is independent of the
point in $\Lambda^\pm$.  \emph{Here the $\pm$ at the front of the right hand
side corresponds to $\Sigma_\pm$, i.e.\ the subscript of
$\Lambda^\pm_\pm$.}  Let
$$
\bar s^\pm(\sigma)=\frac 12 -\hat\beta^\pm(\sigma)=\frac 12+\Im\sigma
$$
denote the threshold Sobolev exponents at
$\Lambda^\pm$, i.e.\ at the future and past radial sets. Thus,
$$
\bar s^+(\sigma)=\bar s^-(\sigma),
$$
but this is actually not important below.
Let $s_{\tow}$ be a function on $S^*X$, such that
\begin{enumerate}
\item
$s_{\tow}$ is constant near $\Lambda^\pm$,
\item
$s_{\tow}$ is decreasing along the $\hamvf_p$-flow on $\Sigma_+$,
increasing on
$\Sigma_-$,
\item
$s_{\tow}$ is less than the threshold exponents at $\Lambda^+$,
towards which we propagate our estimates, i.e.\ 
$s_{\tow}|_{\Lambda^+}<\bar s^+(\sigma)$,
\item
$s_{\tow}$ is greater than the threshold value at
$\Lambda^-$, away from which we propagate our estimates,
i.e.\ $s_{\tow}|_{\Lambda^-}>\bar s^-(\sigma)$.
\end{enumerate}
Since $\ubdry \in H^{s_\tow}$ near $\Lambda^-$ a priori (indeed it is
residual there), one can propagate regularity and estimates from $\Lambda^-$ to
$\Lambda^+$ as in \cite[Section 2.4]{Vasy-Dyatlov:Microlocal-Kerr},
and for all $N$ (in practice taken very large) obtain estimates for
such  $U$
\begin{equation}\label{eq:fd-est-op}
\|\ubdry\|_{H^{s_{\tow}}}\leq C(\|P_\sigma \ubdry\|_{H^{s_{\tow}-1}}+\|\ubdry\|_{H^{-N}}).
\end{equation}
(More generally, the Sobolev exponent on the first term on right hand side
would be $s_{\tow}-m+1$ where $m$ is the order of $P_\sigma;$ here $m=2.$)

On the other hand, if
$s_{\away}$ is a function
on $S^*X$, such that
\begin{enumerate}
\item
$s_{\away}$ is constant near $\Lambda^\pm$,
\item
$s_{\away}$ is increasing along the $\hamvf_p$-flow on $\Sigma_+$, decreasing on
$\Sigma_-$,
\item
$s_{\away}$ is less than the threshold exponents at $\Lambda^-$,
towards which we propagate our estimates, i.e.\ 
$s_{\away}|_{\Lambda^-}<\bar s^-(\sigma)$,
\item
$s_{\away}$ is greater than the threshold value at
$\Lambda^+$, away from which we propagate our estimates,
i.e.\ $s_{\away}|_{\Lambda^+}>\bar s^+(\sigma)$,
\end{enumerate}
then one can propagate regularity and estimates from
$\Lambda^+$ to $\Lambda^-$, and for all $N$ obtain
estimates
$$
\|\ubdry\|_{H^{s_{\away}}}\leq C(\|P_\sigma \ubdry\|_{H^{s_{\away}-1}}+\|\ubdry\|_{H^{-N}}).
$$

With $\bar s^{\pm,*}(\sigma)$ denoting the threshold Sobolev exponents
for $P_\sigma^*$, the same considerations apply to $P_\sigma^*$, i.e.,
if $s^*_{\away}$ is a function on $S^*X$ such that
\begin{enumerate}
\item
$s^*_{\away}$ is constant near $\Lambda^\pm$,
\item
$s^*_{\away}$ is increasing along the $\hamvf_p$-flow on $\Sigma_+$, decreasing on
$\Sigma_-$,
\item
$s^*_{\away}$ is less than the threshold exponents at $\Lambda^-$,
towards which we propagate our estimates, i.e.\ 
$s^*_{\away}|_{\Lambda^-}<\bar s^{-,*}(\sigma)$,
\item
$s^*_{\away}$ is greater than the threshold value at
$\Lambda^+$, away from which we propagate our estimates,
i.e.\ $s^*_{\away}|_{\Lambda^+}>\bar s^{+,*}(\sigma)$,
\end{enumerate}
then one can propagate regularity and estimates from
$\Lambda^+$ to $\Lambda^-$, and for all $N$ obtain
estimates
\begin{equation}\label{eq:back-est-adj}
\|\ubdry\|_{H^{s^*_{\away}}}\leq C(\|P_\sigma^* \ubdry\|_{H^{s^*_{\away}-1}}+\|\ubdry\|_{H^{-N}}),
\end{equation}
with analogous results for $s^*_{\tow}$.

Now, as $\bar s^{\pm,*}(\sigma)=-\bar s^{\pm}(\sigma)+1$, if one chooses
$s_{\tow}$ as above, then one can take $s^*_{\away}=-s_{\tow}+1$:
with this choice,
$$
(H^{s_{\tow}})^*=H^{s^*_{\away}-1},\ (H^{s_{\tow}-1})^*=H^{s^*_{\away}},
$$
i.e., the space on the left hand side of
\eqref{eq:fd-est-op} is dual to the (non-residual) space on the right
hand side of \eqref{eq:back-est-adj}, and (non-residual) the space on the right hand
side of \eqref{eq:fd-est-op} is dual to the space on the left
hand side of \eqref{eq:back-est-adj}. Taking $N$ sufficiently large
such that the inclusions of the spaces on the left hand side of
\eqref{eq:fd-est-op}, resp.\ \eqref{eq:back-est-adj}, into $H^{-N}$ are
compact, this implies Fredholm properties
at once for $P_\sigma$ and $P_\sigma^*$, with a slight change in the
spaces as follows. Let
$$
\cY^{s_{\tow}-1}=H^{s_{\tow}-1},\ \cX^{s_{\tow}}=\{\ubdry\in H^{s_\tow}:\ P_\sigma \ubdry\in \cY^{s_{\tow}-1}\}
$$
(note that the last statement in the definition of $\cX^{s_{\tow}}$ depends on
the principal symbol of $P_\sigma$ only, which is independent of
$\sigma$).

Thus, we finally have the following, which follows from
Propositions~2.3 and 2.4 of \cite{Vasy-Dyatlov:Microlocal-Kerr}
together with the propagation of singularities in variable order
Sobolev spaces away from radial points
(Proposition~\ref{proposition:variablepropagation} in the appendix).

\begin{proposition}\label{proposition:fredholm} The family of maps
  $P_\sigma$ enjoys the following properties:
\begin{enumerate}
\item
$$
P_\sigma:\cX^{s_{\tow}}\to\cY^{s_{\tow}-1},\ P_\sigma^*:\cX^{s_{\away}^*}\to\cY^{s_{\away}^*-1}
$$
are Fredholm.

\item
$P_\sigma$ is a holomorphic Fredholm family on these spaces
in
\begin{equation}
\CC_{s_+,s_-}=\{\sigma\in\CC:\ s_+<\bar s^+(\sigma),\ s_->\bar s^-(\sigma)\},\label{eq:P-sigma-holomorphic}
\end{equation}
with
$$
s_{\tow}|_{\Lambda^\pm}=s_\pm.
$$
$P_\sigma^*$ is antiholomorphic in the same region.

\item If $P_\sigma$ is invertible (or if simply
$u\in\cX^{s_\tow}$, $f\in\cY^{s_{\tow}-1}$, $P_\sigma u=f$), and $\WF(f)\cap\Lambda^-=\emptyset,$ then
$$\WF(P_\sigma^{-1}f)\cap\Lambda^-=\emptyset.$$

\item
If $f$ is $\CI$, then $$\WF(P_\sigma^{-1}f)\subset\Lambda^+.$$ For the adjoint, corresponding to
propagation in the opposite direction, $f \in \CI$ yields
$$
\WF((P^*_\sigma)^{-1}f)\subset \Lambda^-.
$$
\end{enumerate}
\end{proposition}

For the semiclassical problem, with $h^{-1}\equiv\smallabs{\Re \sigma} \to \infty,$ a natural assumption is {\em
  non-trapping}, i.e.\ all semiclassical bicharacteristics in
$\Sigma_{\semi,\pm}$ apart from those in the radial sets
are required to tend to $L^+$ in the forward
direction and $L^-$ in the backward direction in $\Sigma_+$, while the
directions are reversed in $\Sigma_-;$ here $L^\pm$ denotes the image
of $\Lambda^\pm$ in $S^*X$ under the quotient map, and one
considers $S^* X$ as the boundary of the radial compactification of
the fibers of $T^*X$. In particular, the non-trapping assumptions on
$M$ made in Section~\ref{sec:general-hypotheses} imply that the
operator $P_{\sigma}$ is semiclassically non-trapping.

Under this assumption, one has non-trapping semiclassical estimates
(analogues of hyperbolic estimates, i.e.\ with a loss of $h$ relative
to elliptic estimates), which, in the non-semi-classical language
employed here, corresponds to an understanding of asymptotics as $\smallabs{\Re
\sigma} \to \infty.$  The following is proved in the same way as
Theorem~2.15 of \cite{Vasy-Dyatlov:Microlocal-Kerr}.

\begin{proposition}\label{prop:nontrappingestimates}
If the non-trapping hypothesis holds, then:
\begin{enumerate}
\item
$P_{\sigma}^{-1}$ has finitely many poles in each strip
$$
a<\Im \sigma<b.
$$
\item
For all $a,b,$ there exists $C$ such that
$$
\norm{P_{\sigma}^{-1}}_{\cY_{\smallabs{\sigma}^{-1}}^{s_{\tow}-1}\to \cX_{\smallabs{\sigma}^{-1}}^{s_{\tow}}} \leq C
\ang{\Re \sigma}^{-1}
$$
on
$$
a<\Im \sigma<b,\ \smallabs{\Re\sigma}>C.
$$
\end{enumerate}
\end{proposition}

Here the spaces with $\smallabs{\sigma}^{-1}$ subscripts refer to the
variable-order versions of the semiclassical Sobolev spaces discussed
in Section~\ref{sec:semicl}.

\section{Conormality of coefficients}
\label{sec:conormal-reg}

In this section we show that the coefficients in the asymptotic
expansion which will appear in the sequel are in fact classical
conormal distributions with a very explicit singular structure.

\subsection{Conormal and homogeneous distributions}

We begin with some preliminaries on conormal and homogeneous distributions.
For $Y$ a connected component of $\{ v = 0\}$ in $X$ (such as $S_+$),
let $\cM_{\pa}$ denote the module of first order pseudodifferential
operators on $X=\pa M$ with principal symbol vanishing on $N^{*}Y$.
In particular, $\Psi^{0}(X)\subset \cM_{\pa}$.  Note that vector
fields on $X$ tangent to $Y$ lie in $\cM_{\pa}$, and indeed if $A\in
\cM_{\pa}$, then because $N^{*}Y$ is locally defined by $v=0$, $\eta =
0$, and $\gamma$ is elliptic on it, $\sigma_{1}(A) = a_{0}v\gamma +
\sum a_{j}\eta_{j}$, where $a_{j}\in S^{0}$.  In particular, then, $A
= A_{0}(v D_{v}) + \sum A_{j}D_{y_{j}} + A'$, where $A_{j},A'\in
\Psi^{0}(X)$, and so $vD_{v}, D_{y_{j}}$, and $\Id$ generate
$\cM_{\pa}$ as a $\Psi^{0}(X)$-module.

Below we work with the $L^{2}$-based conormal spaces $I^{(s)}(X)$ defined
in Section~\ref{section:mellin} above.
Recall that $u\in I^{(s)}(X)$ means that $u\in H^{s}(X)$
and $A_{1}\ldots A_{k}u\in H^{s}(X)$ for all $k \in \NN$ and $A_{j}\in
\cM_{\pa}$.  Thus $I^{(s)}$ is preserved by elements of $\cM_{\pa}$,
while elements of $\Psi^{k}(X)$ map $I^{(s)}$ to $I^{(s-k)}$.  In
particular, when restricted to a product neighborhood of $Y$, elements
of $I^{(s)}$ can be considered as $\CI$ functions on $Y$ with values
in distributions on $(-\delta , \delta)$ which are conormal to $\{ v=
0\}$, i.e., $I^{(s)}(N^{*}Y) = \CI(Y; I^{(s)}(N^{*}\{0\}))$.  We will
also use the notation
$$
I^{(-\infty)}(N^*Y)\equiv \bigcup_s I^{(s)}(N^*Y).
$$

We also recall the standard conormal spaces, defined using the
$L^{\infty}$-based symbol spaces: $a\in S^{k}(Y\times (-\delta,
\delta); \RR)$ if $a$ is a compactly supported (in the $(y,v)$
variables) and smooth (in all variables) and satisfies the estimates
\begin{equation*}
  \left| D_{y}^{\alpha}D_{v}^{\ell}D_{\gamma}^{N}a\right| \leq C_{\alpha\ell
    N}\langle \gamma\rangle^{k-N}.
\end{equation*}
Elements of $I^{r}(N^{*}Y)$ are defined as certain oscillatory integrals
(which in this case are essentially partial Fourier transforms): $u\in
I^{r}(N^{*}Y)$ if and only if
$$
u = \int e^{\imath v\gamma}a(v,y,\gamma)\,d\gamma\ \text{with}\ a\in
S^{r+(n-3)/4},\ \text{modulo}\ \CI.
$$
Thus $a\in S^{k}$ corresponds to
$u\in I^{k-(n-3)/4}(N^{*}Y)$.  Since $a\in S^{k}$ corresponds to $a$
lying in the weighted $L^{2}$ space $\langle \gamma\rangle^{k+1/2 +
  \epsilon}L^{2}$ for $\epsilon > 0$,
\begin{equation*}
  I^{k-(n-3)/4}(N^{*}Y) \subset \bigcap_{\epsilon > 0}I^{(-k-1/2-\epsilon)}(N^{*}Y).
\end{equation*}
Note that $N^{*}Y$ corresponds to $v=0$ in this parameterization, and
so the principal symbol is identified with an elliptic multiple of
$a|_{v=0}$.

Now, if $a$ is homogeneous outside a compact set in
$\gamma$ (and $a$ is independent of $v$ near $v=0$), one regards it
for convenience as a homogeneous function on $Y\times
(\RR\setminus\{0\})$, and then a basis of such functions of degree
$\kappa$ over $\CI(Y)$ is given by $\gamma ^{\kappa}$ times the
characteristic function of $(0,\infty)_{\gamma}$, resp.\ $(-\infty,
0)_{\gamma},$ which we denote $\gamma_\pm^\kappa.$  If $\kappa$ is not a negative integer, one can go
further, and consider
the homogeneous distributions $\chi_{\pm}^\kappa(\gamma)$ on $\RR$ (or
$Y\times \RR$ in our setting) defined by (the
analytic continuation in $\kappa$, from $\kappa>-1$, when they are
locally $L^1$, of)
$$
\chi_{\pm}^\kappa(\gamma)=\frac{\gamma_\pm^\kappa}{\Gamma(1+\kappa)}.
$$
The inverse Fourier transform of these distributions are
elliptic multiples of
\begin{equation*}
v_{\pm \imath 0}^{-1-\kappa} \equiv (v\pm \imath 0)^{-1-\kappa}
\end{equation*}
 (see
Section~7.1 of \cite{Hormander:v1});
these are thus a basis for $I^{k}_{\cl}(N^{*}Y)/I^{k-1}_{\cl}(N^{*}Y)$ for
$k = \kappa - (n-3)/4$ over $\CI(Y)$.  (The ``cl'' subscript stands for
``classical'' and refers to conormal distributions whose symbols have
polyhomogeneous asymptotic expansions.)  For negative integers $\kappa=-k$,
one must be more careful in describing a basis, as $\chi_\pm^{-k}$
is then supported at the origin.  We instead simply consider directly the inverse
Fourier transform of $\psi(\pm \gamma) \gamma^{-k}$ where $\psi$ is a
smooth function equal to $0$ for $\gamma<1$ and $1$ for $\gamma>2.$
The result is a $k$'th antiderivative of $\psi(\pm\gamma),$ whose
inverse Fourier transform differs by a smooth function from a multiple
of $(v\pm \imath 0)^{-1},$ hence differs by a smooth function from a
multiple of
\begin{equation}\label{homogdist}
  v_{\pm \imath 0}^{-1-\kappa}\equiv (v \pm \imath 0)^{-1-\kappa}\log (v
  \pm \imath 0).
\end{equation}
Note that these are no longer homogeneous distributions.
(We also remark that when $\kappa$ is a negative integer of course we may
also write more simply
$$
v_{\pm \imath 0}^{-1-\kappa}=v^{-1-\kappa} \left(\log\smallabs{v} \pm \imath
\pi H(-v)\right)
$$
with $H$ the Heaviside function; however it is more convenient to
stick with the consistent notation offered by the expression as
\eqref{homogdist}.)

\subsection{Spaces of solutions}

We now turn to the solution spaces of a class of operators 
including our $P_\sigma.$.  We consider a general operator of
the form
\begin{equation}\label{modelform}
  P = vD_{v}^{2} + \alpha D_{v} + Q,\ Q\in \cM_{\pa}^{2},\ \alpha \in \CI(Y);
\end{equation}
since $D_{v}$ is elliptic on $N^{*}Y$, we in particular have $P\in
\cM_{\pa}\cdot \Psi^{1}(\pd
M)$.
Note in particular that the operator family $P_\sigma$ defined by \eqref{Psigmadef} has
the form \eqref{modelform} by \eqref{Psigma1}, hence the results here apply
if $P_\sigma u \in \CI.$

\begin{lemma}
  \label{lemma:asymp-leading}
  If $Pu = f\in I^{(s)}(N^{*}Y)$, $u\in I^{(s)}(N^{*}Y)$, then
  \begin{equation}
    \label{eq:conormal-expansion1}
    u = g_{+}v_{+\imath 0}^{1-\imath \alpha} + g_{-}v_{-\imath 0}^{1-\imath \alpha} + \tilde{u},
  \end{equation}
  with $g_{\pm}\in \CI(Y)$ (pulled back via a local product
  decomposition) and $\tilde{u} \in I^{(s+1-\epsilon)}(N^{*}Y)$ for
  all $\epsilon > 0$.
\end{lemma}
\begin{remark}
  If $v_{\pm\imath 0}^{1-\imath\alpha}\in I^{(s+1-\epsilon)}$ for all
  $\epsilon > 0$, then the conclusion is simply $u\in
  I^{(s+1-\epsilon)}(N^{*}Y)$.  On the other hand, if $v_{\pm \imath
    0}^{1-\imath\alpha}\notin I^{(s)}$, then the conclusion is
  $g_{\pm}=0$, and thus again $u\in I^{(s+1-\epsilon)}(N^{*}Y)$.

  If $f \in \CI(\pa M),$ iterative use of the lemma will yield a full
  expansion of $u$, provided we replace $g_{\pm}$ by appropriate
  functions $\tilde{g}_{\pm}$ with $P(\tilde{g}_{\pm}v_{\pm\imath 0}^{1-\imath\alpha}) \in \CI(X)$ (see
  Lemma~\ref{lemma:approx-soln-expn} below).
\end{remark}

\begin{proof}
  We may assume that $u$ is supported in a product neighborhood of
  $Y$, identified as $(-\delta, \delta)_{v}\times Y$, since if $\chi
  \in \CI(\pa M)$ is compactly supported in such a neighborhood and is
  identically $1$ near $Y$, then $\WF'\left([\chi, P]\right) \cap
  N^{*}Y = \emptyset$, so $[\chi, P]u\in \CI(X)$ and thus $P(\chi
  u) \in I^{(s)}$ as well. 

  Note that $vD_{v}^{2} = D_{v}vD_{v} + \imath D_{v}$.  Thus, if $G\in
  \Psi^{-1}(X)$ is a parametrix for $D_{v}$ near $N^{*}Y$ (where
  $D_{v}$ is elliptic), applying $G$ to $Pu$ yields
  \begin{equation*}
    (vD_{v} + (\imath+\alpha) + GQ)u\in I^{(s+1)}(N^{*}Y).
  \end{equation*}
  Since $u\in I^{(s)},$  $Qu\in I^{(s)}$ and thus $GQu\in
  I^{(s+1)}$, so we have
  \begin{equation*}
    (vD_{v} + \imath + \alpha)u \in I^{(s+1)} = \CI (Y; I^{(s+1)} (N^{*}\{0\})).
  \end{equation*}
  With $J$ a compact interval, let $I^{(\ell)}_{\cS}(N^{*}\{0\})$ denote
  the sum of elements of $I^{(\ell)}(N^{*}\{0\})$ supported in $J$ and
  Schwartz functions on $\RR$.  Then the Fourier transform
  on $\RR$ maps elements of $I^{(\ell)}_{\cS}(N^{*}\{0\})$ to
  $L^{2}$-based symbols.  More precisely, if $S^{(\ell)}$ is the
  set of smooth functions $\phi$ on $\RR_{\gamma}$ such
  that
 $$(\gamma
 D_{\gamma})^{N}\phi\in L^{2,\ell} \equiv \langle \gamma
 \rangle^{-\ell}L^{2}$$ for all $N\in \NN$, then the Fourier transform
 is an isomorphism $I^{(\ell)}_{\cS}(N^{*}\{0\})\to S^{(\ell)}$.  (We
 remark that while our conventions yield a nice correspondence between
 conormal Sobolev orders and symbol orders as discussed here, they
 have the unfortunate result that $S^{\ell} \subset S^{\ell'}$ if
 $\ell\geq \ell',$ contrary to the usual sign convention for symbol
 spaces.) 

  Taking the partial Fourier transform, $\widetilde{\cF}$, in the interval
  variable, $v$, yields
  \begin{equation*}
    (-\gamma D_{\gamma} + 2\imath + \alpha)\widetilde{\cF} u =
    (-D_{\gamma}\gamma + \imath + \alpha)\widetilde{\cF}u \in \CI(Y; S^{(s+1)}).
  \end{equation*}
  Now, to analyze the behavior of $\widetilde{\cF}u$ at infinity, we
  conjugate the differential operator by $\gamma ^{-2+\imath\alpha}$ on
  $\RR\setminus\{0\}$, where
  \begin{equation*}
    \gamma^{2-\imath\alpha}(-\gamma D_{\gamma} + 2\imath +
    \alpha)\gamma^{-2 + \imath \alpha} = -\gamma D_{\gamma},
  \end{equation*}
  so one has
  \begin{equation*}
    -\gamma D_{\gamma}(\gamma^{2-\imath\alpha}\widetilde{\cF}u) = \gamma
    ^{2-\imath\alpha}(-\gamma D_{\gamma} + 2\imath +
    \alpha)\widetilde{\cF}u \in \CI(Y; S^{(s-1-\im\alpha)}),
  \end{equation*}
  and thus
  \begin{equation}\label{crap1}
    D_{\gamma} (\gamma ^{2-\imath\alpha}\widetilde{\cF}u)|_{[1,\infty)} \in \CI(Y;
    S^{(s-\im \alpha)}[1,\infty)).
  \end{equation}
  Note that due to the presence of $\epsilon > 0$ in the statement of
  the lemma, we may assume that $s - \im \alpha \neq 1/2$; this
  simplifies some formulae below (otherwise one would have logarithmic
  terms).  

  Now, if $b\in S^{(\ell)}([1,\infty))$, $\ell < 1/2$, then the
  indefinite integral of $b$ given by
  \begin{equation*}
    c(\gamma) = \int _{1}^{\gamma}b(\eta)\,d\eta,
  \end{equation*}
  satisfies (by Cauchy-Schwarz)
  \begin{equation*}\begin{aligned}
    |c(\gamma)|&\leq \left(
      \int_{1}^{\gamma}|\eta|^{2\ell}|b(\eta)|^{2}\,d\eta\right)^{1/2}\left(
      \int_{1}^{\gamma}|\eta|^{-2\ell}\,d\eta\right)^{1/2} \\
&\leq
    C\Norm[L^{2,\ell}]{b}\left(1 + |\gamma|^{\frac{1}{2}-\ell}\right).
  \end{aligned}\end{equation*}
  Thus $c\in L^{2,\ell-1-\epsilon}$ for all $\epsilon > 0$, and as
  $D_{\gamma}c = b$, $c\in S^{(\ell-1-\epsilon)}$.  (Note that
  constants are in $S^{(\ell-1-\epsilon)}$ since $\ell<1/2.$)

  Returning now to $u$ described by \eqref{crap1} above and setting
  $\ell=s-\Im \alpha,$ we see that
 \begin{equation*}
    \gamma ^{2-\imath\alpha}\widetilde{\cF}u=\tilde{w}\in S^{(s-\im \alpha
    - 1 - \epsilon)}
  \end{equation*}
provided $\ell<1/2.$
  On the other hand, if $\ell=s-\Im \alpha > 1/2$, then
  $S^{(\ell)}\subset L^{1}$, and if we define the indefinite integral
  as
  \begin{equation*}
    c(\gamma ) = -\int_{\gamma}^{\infty}b(\eta)\,d\eta,
  \end{equation*}
  then, by Cauchy--Schwarz,
  \begin{equation*}
    |c(\gamma)|\leq \left( \int
      _{\gamma}^{\infty}|\eta|^{2\ell}|b(\eta)|^{2}\,d\eta\right)^{1/2}
    \left( \int_{\gamma}^{\infty}|\eta|^{-2\ell}\,d\eta\right)^{1/2}
    \leq C\Norm[L^{2,\ell}]{b}|\gamma|^{\frac{1}{2}-\ell},
  \end{equation*}
  so $c\in S^{(\ell-1-\epsilon)}$.  Then, writing
  \begin{equation*}
    \gamma^{2-\imath\alpha}\widetilde{\cF}u-\widetilde{\cF}u|_{\gamma=1} =
    \int_{1}^{\infty}\pa_{\gamma}(\gamma^{2-\imath\alpha}\widetilde{\cF}u) - \int_{\gamma}^{\infty}\pa_{\gamma}(\gamma^{2-\imath\alpha}\widetilde{\cF}u) ,
  \end{equation*}
  we deduce that
  \begin{equation*}
    \gamma^{2-\imath\alpha}\widetilde{\cF}u|_{\gamma>1} = g_{+} + \tilde{w},\ g_{+}\in
    C^{\infty}(Y),\ \tilde{w}\in S^{(s-\im\alpha-1-\epsilon)},
  \end{equation*}
  and thus
  \begin{equation*}
    \widetilde{\cF}u|_{\gamma > 1} = g_{+}\gamma^{-2+\imath\alpha} +
    w_{+},\ w_{+}\in S^{(s+1-\epsilon)}.
  \end{equation*}
  A similar calculation applies to $\widetilde{\cF}u|_{\gamma < -1}$, yielding
  \begin{equation*}
    \widetilde{\cF}u|_{\gamma < -1} = g_{-}(-\gamma)^{-2+\imath\alpha} +
    w_{-} ,\ w_{-}\in S^{(s+1-\epsilon)}.
  \end{equation*}
  In summary, if $\psi_{+}$ is supported in $(1,\infty)$, identically
  $1$ on $[2,\infty)$, and $\psi_{-}(\gamma) = \psi_{+}(-\gamma)$,
  then
  \begin{equation*}
    \widetilde{\cF}u = g_{+}\psi_{+}\gamma^{-2+\imath\alpha} +
    g_{-}\psi_{-}(-\gamma)^{-2+\imath\alpha} + w,\ w\in S^{(s+1-\epsilon)}.
  \end{equation*}
  Now, the inverse partial Fourier transform of $w$ is in
  $I^{(s+1-\epsilon)}$, so it remains to deal with the other terms.
  Changing these by a compactly supported distribution does not affect
  their singularities, so we can replace these by the homogeneous
  distributions $\gamma_{\pm}^{-2+\imath\alpha}$ for a local
  description of the inverse partial Fourier transform.  But the
  inverse Fourier transforms of the latter are $v_{\pm\imath 0}^{1-\imath\alpha}$, so we conclude that
  \begin{equation*}
    u = g_{+}v_{+\imath 0}^{1-\imath\alpha} + g_{-}v_{-\imath 0}^{1-\imath\alpha} + \tilde{u}, \ \tilde{u} \in I^{(s+1-\epsilon)},
  \end{equation*}
  as claimed.
\end{proof}

Although the following corollary follows directly from the results of
\cite{Vasy-Dyatlov:Microlocal-Kerr},
we give a proof using
Lemma~\ref{lemma:asymp-leading}.
\begin{corollary}
  \label{cor:no-low-sing}
  If $Pu=f\in \CI(X)$, $u\in I^{(s_{0})}(N^{*}Y)$, $s_{0} > 3/2 +
  \im \alpha$, then $u\in \CI(X)$.  
\end{corollary}

\begin{proof}
  Let $\tilde{s}_{0} = \sup \{ s: u\in I^{(s)}(N^{*}Y)\}$, so
  $\tilde{s}_{0} > 3/2 + \im \alpha$ (possibly
  $\tilde{s}_{0}=+\infty$); if $\tilde{s}_{0} = +\infty$, then we are
  done as $\bigcap _{s\in\RR}I^{(s)} = \CI(X)$.  Thus, $u\in
  I^{(\tilde{s}_{0}-\epsilon)}(N^{*}Y)$ for all $\epsilon > 0$.  By
  Lemma~\ref{lemma:asymp-leading},
  \begin{equation*}
    u = g_{+}v_{+\imath 0}^{1-\imath\alpha} + g_{-}v_{-\imath 0}^{1-\imath\alpha} + \tilde{u},
  \end{equation*}
  with $g_{\pm} \in \CI(Y)$ (pulled back via a local product
  decomposition) and $\tilde{u}\in
  I^{(\tilde{s}_{0}+1-\epsilon)}(N^{*}Y)$ for all $\epsilon > 0$. 
For all $\epsilon > 0$, $\tilde{u}\in
  I^{(\tilde{s}_{0}+1-\epsilon)}(N^{*}Y)$, which is a subset of
  $I^{(3/2+\im\alpha)}(N^{*}Y)$ for sufficiently small $\epsilon > 0.$
On the other hand the
  sum of the first two terms is not in $I^{(3/2+\im \alpha)}(N^{*}Y)$
  unless $g_{\pm}$ vanish. Since $u \in I^{(3/2+\im \alpha)}(N^{*}Y),$
  $g_{\pm}$ must vanish, and thus $u=\tilde{u}\in
  I^{(\tilde{s}_{0}+1-\epsilon)}(N^{*}Y)$ for all $\epsilon >0$,
  contradicting the definition of $\tilde{s}_{0}$.  Thus,
  $\tilde{s}_{0} = +\infty.$
\end{proof}

Next, under the assumption that $\alpha$ is constant, we show that
distributions such as those in the first two terms on the right hand
side of the equation~\eqref{eq:conormal-expansion1} can be modified to
elements of the nullspace of $P$ modulo $\CI(X)$.

\begin{lemma}
  \label{lemma:approx-soln-expn}
  Suppose $\alpha \in \CC$ is a constant, $1-\imath\alpha$ is not an
  integer, and $g\in \CI(Y)$.  Then there exist $u_{\pm} = g_{\pm}
  v_{\pm\imath 0}^{1-\imath\alpha}\in \bigcap _{\epsilon > 0} I^{(3/2 +
    \im \alpha - \epsilon)}(N^{*}Y)$, with $g_{\pm}\in C^{\infty}(X)$
  such that $g_{\pm}|_{Y} = g$ and $Pu_{\pm}\in C^{\infty}(X)$.
\end{lemma}

\begin{remark}
  \label{rem:precise-expansion}
  If $1 - \imath\alpha$ is an integer, the proof below still proves a
  slightly different result: logarithmic terms appear.  Indeed, if
  $1-\imath\alpha$ is a nonnegative integer, then logarithmic terms
  appear from the definition of $v_{\pm\imath 0}^{1-\imath\alpha}$.
  If it is a negative integer, say, $1-\imath\alpha = -r \leq -1$, then an
  additional logarithmic term is incurred at the $r$-th step in the
  expansion.  

  It is more straightforward to state it as follows: $u_{\pm}$ is a
  classical conormal distribution of the appropriate order, with
  principal symbol the same as that of $gv_{\pm \imath
    0}^{1-\imath\alpha}$.
\end{remark}

\begin{remark}
  A similar expansion can be obtained in general, without assuming
  that $\alpha$ is a constant.  This is similar to the treatment of
  generalized Coulomb type spherical waves in \cite{Vasy:Geometric}.
\end{remark}

\begin{proof}
  We suppose first that $1-\imath\alpha$ is not an integer.
  As the indicial roots associated to the ordinary differential
  operator $vD_v^2+\alpha D_v$ are $0$ and $1-\imath\alpha,$ for $h\in
  \CI(X)$,
  \begin{equation*}
    Pv^{k}v_{\pm\imath 0}^{1-\imath\alpha} h = v^{k}v_{\pm \imath 0}^{-\imath\alpha} w,\ w\in \CI(X), w|_{Y} = c(k)h|_{Y},
  \end{equation*}
  with $c(0) = 0$, $c(k)\neq 0$ for $k\neq 0$.  (We suppress the
  dependence of $c(k)$ on $\alpha$.)  Correspondingly, given $g$,
  consider first $h_{\pm , 0}\in \CI(X)$ with $h_{\pm,0}|_{Y} = g$.  Then
  \begin{equation*}
    Pv^{0}v_{\pm\imath 0}^{1-\imath\alpha}h_{\pm,0} = v^{0}v_{\pm
      \imath 0}^{-\imath\alpha} w
  \end{equation*}
  with $w|_{Y} = 0$, so in fact
  \begin{equation*}
    Pv^{0}v_{\pm\imath 0}^{1-\imath\alpha}h_{\pm,0} =
    v^{1}v_{\pm\imath 0}^{-\imath\alpha} \tilde{w}_{\pm,1}
  \end{equation*}
for some $\tilde{w}_{\pm,1} \in \CI(X).$
  Now, in general, for $k\neq 0$, given $\tilde{w}_{\pm,k}\in \CI(X)$,
  one can let $h_{\pm,k} = -c(k)^{-1}\tilde{w}_{\pm,k}$, and then
  \begin{equation*}
    Pv^{k}v_{\pm\imath 0}^{1-\imath\alpha} h_{\pm ,k} + v^{k}v_{\pm
      \imath 0}^{-\imath\alpha}\tilde{w}_{\pm,k} = v^{k}v_{\pm \imath 0}^{-\imath\alpha}w_{\pm,k}
  \end{equation*}
  with $w_{\pm,k}|_{Y}=0$, thus the right hand side is of the form
  $v^{k+1}v_\pm^{-\imath\alpha}\tilde{w}_{\pm,k+1}$.
  Correspondingly, one can proceed inductively and construct
  $\tilde{h}_{\pm,k}$ with
  \begin{equation*}
    Pv_{\pm \imath 0}^{1-\imath\alpha}\tilde{h}_{\pm ,k} =
    v^{k+1}v_{\pm \imath 0}^{-\imath\alpha}\tilde{w}_{\pm,k+1}
  \end{equation*}
  with $\tilde{w}_{\pm,k+1}\in \CI(X),$ e.g.\ by taking $\tilde{h}_{\pm,k}
  = \sum_{j=0}^{k}v^{j}h_{\pm,j}.$  More generally, one can
  asymptotically sum the series $\sum_{j=0}^{\infty}v^{j}h_{\pm,j}$, i.e.,
  construct a function $h_{\pm}$ which differs from
  $\sum_{j=0}^{k}v^{j}h_{\pm,j}$ by a $\CI$ function vanishing to order
  $k+1$; then
  \begin{equation*}
    Pv_{\pm\imath 0}^{1-\imath\alpha}h_{\pm} = v^{k+1}v_{\pm\imath 0}^{-\imath\alpha}W_{\pm,k+1}
  \end{equation*}
  for every $k$ for some $W_{\pm,k+1}\in \CI(X)$, thus the right
  hand side is $\CI$, completing the proof.

  If $1-\imath\alpha$ is a non-negative integer, then the iterative
  construction requires including another logarithmic term owing to
  the logarithmic term in $v_{\pm}^{1-\imath\alpha}$.  If
  $1-\imath\alpha =r\leq -1$ is a negative integer, then the iterative
  construction breaks down when finding the coefficient of
  $v^{r}v_{\pm\imath 0}^{1-\imath\alpha}$, as in this setting $c(r) = 0$.  The
  proof goes through nearly as stated once we also include terms of
  the form $v^{r+k}v_{\pm\imath 0}^{1-\imath\alpha}\log(v\pm \imath 0)$
  for $k \geq 0$.
\end{proof}

In addition to knowing that we may formally parametrize elements in the
approximate nullspace by functions on $Y,$ we will need to know how to
formally solve certain inhomogeneous equations with specified conormal
right-hand sides.  For the following lemma, we assume that $Q$ is a
\emph{differential} operator in the module $\mathcal{M}_{\pa}^{2}$,
although it holds (with a slightly more complicated proof) if $Q$ is
pseudodifferential.  Note that for our operator $P_{\sigma}$, $Q$ is
in fact differential.  
\begin{lemma}\label{lemma:conormalinhomog}
  Suppose $\alpha \in \CC$ is constant and $Q$ is a
  \emph{differential} operator in $\mathcal{M}_{\pa}^{2}$.  Let $h \in
  \CI(X)$ and let $m$ be a nonnegative integer.  If $\imath\alpha$ is not a
  strictly positive integer, then there exist $g_\pm^0,\dots
  g_\pm^{m+1} \in \CI(X)$ such that the functions
  $$
  u_\pm =\sum_{m'=0}^{m+1} g_\pm^{m'} v_{\pm\imath 0}^{1-\imath\alpha} \log (v\pm \imath 0)^{m'}
  $$
  solve
  $$
  P u_\pm = h v_{\pm\imath 0}^{-\imath\alpha} \log (v\pm \imath 0)^m+\tilde{u}
  $$
  with $\tilde{u} \in \CI$ and $P$ as in \eqref{modelform}.  If $1-\imath\alpha=-k\leq 0$ is a non-positive
  integer, then the same statement is true with $u_{\pm}$ replaced by
  $u_{\pm }+ w_{\pm}$, where $w_{\pm}$ has the form
  \begin{equation*}
    w_{\pm} = g^{m+2}_{\pm}v^{k-1} v_{\pm\imath
      0}^{1-\imath\alpha}\log(v\pm \imath 0)^{m+2},
  \end{equation*}
 with
  $g^{m+2}_{\pm}$ a smooth function.
  (If $1-\imath\alpha$ is a non-negative integer, then there is an additional
  log term implicit in the formula owing to the definition of $v_{\pm
    \imath 0}^{1-\imath\alpha}$.)
\end{lemma}

\begin{proof}
  The proof is similar to that of Lemma~\ref{lemma:approx-soln-expn}.
  Indeed, as the indicial roots of $vD_{v}^{2} + \alpha D_{v}$ are $0$
  and $1-\imath\alpha$, for $g\in C^{\infty}(X)$,
  \begin{align}
    \label{eq:application-expansion}
    P\left( v^{k}v_{\pm\imath 0}^{1-\imath\alpha}(\log (v\pm \imath
      0))^{m'}g\right) &= \sum_{\ell =0}^{\min(m',2)}
    v^{k}v_{\pm\imath 0}^{-\imath \alpha} (\log (v\pm \imath
    0))^{m'-\ell}w^{(\ell)}, \\
    w^{(\ell)}\in \CI (X), &\quad w^{(\ell)}|_{Y} = c(k,m',\ell)
    g|_{Y}. \notag
  \end{align}
  We may calculate the $c(k,m',\ell)$ explicitly.  If $\imath\alpha$
  is not a negative integer, then
  \begin{align*}
    c(k,m',0) &= -k (k + 1 - \imath\alpha), \\
    c(k,m',1) &= -m' (2k + 1 - \imath\alpha), \\
    c(k,m',2) &= - m' (m'-1).
  \end{align*}
  If $\imath\alpha$ is a negative integer, then
  \begin{align*}
    c(k,m',0) &= -k (k+1-\imath\alpha), \\
    c(k,m',1) &= -(m'+1) (2k+1-\imath\alpha), \\
    c(k,m',2) &= -(m'+1)m' .
  \end{align*}
  In particular, we always have $c(0,m',0)=0$.  When $\imath\alpha
  \neq 1$, then $c(0,m',1)\neq 0$, while for $\imath\alpha = 1$,
  $c(0,m',1)=0$ as well, though $c(0,m',2)\neq 0$.  If $\imath\alpha$ is
  not a positive integer, then $c(k,m',0) \neq 0$ for all $k$.  If
  $\imath\alpha = r\neq 1$ is a positive integer, then $c(r-1,m',0)=0$ as
  well.  

  We start by assuming that $\imath\alpha$ is not a positive integer,
  so that $c(0,m',1) \neq 0$ and $c(k,m',0) \neq 0$ for all $k\geq
  1$.  We find the Taylor series for $g_{\pm}^{m'}$ iteratively; at
  each step we find the largest remaining term, starting with the
  leading term of $g^{m+1}_{\pm}$, then the leading term of
  $g^{m}_{\pm}$, and so on.  We then find the first Taylor coefficient
  of $g^{m+1}_{\pm}$, then $g^{m}_{\pm}$, and continue in this
  fashion.  

  We start by finding the leading term in the expansion.  For any
  function $g_{\pm,0}^{m+1} \in C^{\infty}(X)$, we have
  \begin{align*}
    &Pv^{0}v_{\pm\imath 0}^{1-\imath\alpha} (\log (v\pm \imath
    0))^{m+1}g_{\pm,0}^{m+1} \\
    &\quad = \sum_{\ell = 0}^{\min (m+1,2)}v^{0}v_{\pm\imath
      0}^{-\imath\alpha}(\log (v\pm\imath 0))^{m+1-\ell}
    w^{(\ell)}_{\pm, 0, m+1},
  \end{align*}
  with $w^{(\ell)}_{\pm,0,m+1} = c(0,m+1,\ell) g_{\pm,0}^{m+1}$.  As
  $c(0,m+1,0) = 0$, and $c(0,m+1,1)\neq 0$, we choose
  $g_{\pm,0}^{m+1}$ so that 
  \begin{equation*}
    g_{\pm,0}^{m+1}|_{Y} = \frac{1}{c(0,m+1,1)}h|_{Y}.
  \end{equation*}
  
  For a function $g^{m'}_{\pm,0}\in C^{\infty}$, we have
  \begin{align*}
    &Pv^{0}v_{\pm\imath 0}^{1-\imath\alpha}(\log(v \pm \imath
    0))^{m'}g^{m'}_{\pm,0} \\
    &= \sum_{\ell =0}^{\min (m',2)} v^{0}v_{\pm\imath 0}^{-\imath\alpha}(\log (v\pm \imath 0))^{m'-\ell} w^{(\ell)}_{\pm,
    0, m'}.
  \end{align*}
  
  Having now found $g_{\pm,0}^{m'+1}$, we then choose $g^{m'}_{\pm,0}$
  so that
  \begin{equation*}
    g_{\pm,0}^{m'}|_{Y} = - \frac{c(0,m'+1,2)}{c(0,m',1)} g^{m'+1}_{\pm,0}|_{Y}.
  \end{equation*}
  Observe that because $v_{\pm\imath 0}^{1-\imath\alpha}$ is in the
  approximate kernel of $P$, we may choose $g^{0}_{\pm,0}|_{Y}$ freely.

  Applying $P$, all terms other than the leading one cancel at $Y$ and
  so we have
  \begin{align*}
    &P\sum_{m' = 0}^{m+1}v^{0}v_{\pm\imath 0}^{1-\imath\alpha} (\log
    (v\pm \imath 0))^{m'}g^{m'}_{\pm,0}  - v_{\pm \imath 0}^{-\imath\alpha}(\log (v\pm \imath 0))^{m}h \\
    &\quad = \sum_{m'=0}^{m+1}v^{1}v_{\pm\imath 0}^{-\imath \alpha}
    (\log (v\pm \imath 0))^{m'} \tilde{w}_{\pm, 1, m'}.
  \end{align*}

  Now, in general, for $k\neq 0$, given $\tilde{w}_{\pm, 1,m'}\in
  \CI(X)$ for $0 \leq m' \leq m+1$, we set
  \begin{equation*}
    g^{m'}_{\pm,k} = -
    \frac{1}{c(k,m',0)}\left(\tilde{w}_{\pm,k,m'} +
      c(k,m'+1,1)w^{(1)}_{\pm,k,m'+1}
      +c(k,m'+2,2)w^{(2)}_{\pm,k,m'+2}\right), 
  \end{equation*}
  where $w^{(\ell)}_{\pm,k,m'}$ are the coefficients in
  equation~\eqref{eq:application-expansion} with applied to $g =
  g^{m'}_{\pm, k}$.  Applying $P$, we see
  \begin{align*}
    &P \sum_{m'=0}^{m+1}v^{k}v_{\pm \imath 0}^{1-\imath\alpha}(\log (v\pm \imath 0))^{m'}g^{m'}_{\pm, k} +
    \sum_{m'=0}^{m+1}v^{k}v_{\pm\imath 0}^{-\imath\alpha}\tilde{w}_{\pm,k,m'} \\
    &\quad\quad\quad\quad - v_{\pm\imath 0}^{-\imath \alpha}(\log (v\pm \imath
    0))^{m}h \\
    &\quad = \sum_{m'=0}^{m+1}v^{k}v_{\pm\imath 0}^{-\imath\alpha}w^{(0)}_{\pm,k,m'},
  \end{align*}
  where $w^{(0)}_{\pm,k,m'}|_{Y} = 0$, and so the right hand side is
  of the form
  \begin{equation*}
    \sum_{m'=0}^{m+1}v^{k+1}v_{\pm\imath 0}^{-\imath\alpha}\tilde{w}_{\pm,k+1,m'}.
  \end{equation*}
  We can thus proceed inductively and construct
  $\tilde{g}^{m'}_{\pm,k}$ with
  \begin{align*}
    &P \sum_{m'=0}^{m+1}v_{\pm\imath 0}^{1-\imath\alpha}(\log (v\pm \imath 0))^{m'}\tilde{g}^{m'}_{\pm,
      k}  - v_{\pm\imath 0}^{-\imath\alpha}(\log (v\pm \imath
    0))^{m}h \\
    &\quad = \sum_{m'=0}^{m+1}v^{k+1}v_{\pm\imath 0}^{-\imath\alpha}\tilde{w}_{\pm,k+1,m'},
  \end{align*}
  with $\tilde{w}_{\pm,k+1,m'}\in \CI(X)$.  (Namely,
  $\tilde{g}^{m'}_{\pm,k} = \sum_{j=0}^{k}v^{j}g^{m'}_{\pm,j}$ works.)  

  We now asymptotically sum the series
  $\sum_{j=0}^{\infty}v^{j}g^{m'}_{\pm,j}$ to construct a function
  $g^{m'}_{\pm}$ differing from each
  $\sum_{j=0}^{k}v^{j}g^{m'}_{\pm,j}$ by a smooth function vanishing
  to order $k+1$, and then
  \begin{align*}
    &P \sum_{m'=0} ^{m+1}v_{\pm\imath 0}^{1-\imath\alpha}(\log (v\pm
    \imath 0))^{m'} g^{m'}_{\pm} - v_{\pm\imath 0}^{-\imath
      \alpha}(\log (v\pm \imath 0))^{m}h \in \CI(X),
  \end{align*}
  completing the proof.

  If $1-\imath\alpha = -k$ is a non-positive integer, the iteration proceeds
  nearly as before, but at the expense of an additional log term at the
  $k$-th coefficient.  (For example, if $k=0$, then $c(0,m',1)$ also
  vanishes and so an additional log term is needed to find the first
  coefficient.)
\end{proof}

We now combine Lemma~\ref{lemma:asymp-leading} and
Lemma~\ref{lemma:approx-soln-expn} to obtain a complete asymptotic
expansion of elements of the nullspace of $P$ modulo $\CI(X)$.

\begin{proposition}\label{prop:conormal-inhomog}
  Suppose $\alpha \in \CC$ is a constant.  If $u\in I^{(-\infty)}(N^{*}Y)$ and $Pu\in \CI(X)$, then there exist $g_{\pm}\in
  \CI(X)$ and $\tilde{u}\in \CI(X)$ such that
  \begin{equation*}
    u = g_{+}v_{+\imath 0}^{1-\imath\alpha} + g_{-}v_{-\imath 0}^{1-\imath\alpha} + \tilde{u}.
  \end{equation*}
  See Remark~\ref{rem:precise-expansion} if $1-\imath\alpha$ is an
  integer.

If instead $Pu\in I^{(\tilde s)}(N^*Y)$, the same conclusion holds with
$\tilde{u}\in \CI(X)$  replaced by $\tilde u\in \bigcap_{\delta>0} I^{(\tilde s+1-
\delta)}(N^*Y)$.
\end{proposition}

\begin{proof}
First suppose $Pu\in\CI(X)$.
  Let $s_{0} = \sup \{ s: u\in I^{(s)}(N^{*}Y)\}$ (the set on the
  right is non-empty by hypothesis), so either $s_{0} = +\infty$, and
  then $\bigcap _{s\in\RR} I^{(s)} = \CI(X)$ shows that the
  conclusion holds with $g_{\pm}=0$, or $s_{0}\in \RR$ is finite, and
  then $u\in I^{(s_{0}-\epsilon)}$ for all $\epsilon > 0$.  By
  Lemma~\ref{lemma:asymp-leading}, there exist $\tilde{g}_{\pm}\in
  \CI(Y)$ so that
  \begin{equation*}
    u = \tilde{g}_{+}v_{+\imath 0}^{1-\imath\alpha} +
    \tilde{g}_{-}v_{-\imath 0}^{1-\imath\alpha} + u',
  \end{equation*}
  with $u' \in I^{(s_{0}+1-\delta)}(N^{*}Y)$ for all $\delta >0$.
  Here the first two terms are in $\bigcap _{\delta > 0} I^{(3/2 +
    \im\alpha-\delta)}(N^{*}Y)$ but not in $I^{(3/2 + \im
    \alpha)}(N^{*}Y)$ unless $\tilde{g}_{\pm}$ vanish; by the
  assumption that $s_{0}<\infty$, we find $3/2 + \im \alpha = s_{0}$ and
  $\tilde{g}_{\pm}$ cannot both vanish.  Let $g_{\pm}\in
  C^{\infty}(X)$, $u_{\pm} \in \bigcap_{\delta > 0}I^{(3/2 +
    \im\alpha - \delta)}(N^{*}Y)$ be given by
  Lemma~\ref{lemma:approx-soln-expn} with $\tilde{g}_{\pm}$ in place
  of $g$.  Thus, $Pu_{\pm}\in \CI(X)$, hence
  $P(u-u_{+}-u_{-})\in \CI(X)$.  Further,
  \begin{equation*}
    u - u_{+} - u_{-} = (g_{+}-\tilde{g}_{+})v_{+\imath
      0}^{1-\imath\alpha} + (g_{-}-\tilde{g}_{-})v_{-\imath 0}^{1-\imath\alpha} + u',
  \end{equation*}
  and
  \begin{equation*}
    (g_{\pm}-\tilde{g}_{\pm})v_{\pm\imath 0}^{1-\imath\alpha} =
    O(v) v_{\pm\imath 0}^{1-\imath\alpha} \in \bigcap _{\delta
      > 0}I^{(5/2 + \im \alpha - \delta)}(N^{*}Y).
  \end{equation*}
  Thus, $u-u_{+}-u_{-}\in \bigcap _{\delta > 0}I^{(5/2 + \im\alpha -
    \delta)}(N^{*}Y)$.  By Corollary~\ref{cor:no-low-sing},
  $u-u_{+}-u_{-}\in C^{\infty}(X)$, completing the proof of the
  proposition in the first case.

If $Pu\in I^{(\tilde s)}(N^*Y)$ instead, then defining $s_0$ as above,
$s_0\geq \tilde s+1$ means that we are done, so assume $s_0<\tilde
s+1$. Proceeding as above, we obtain
$\tilde{g}_{\pm}\in
  \CI(Y)$ so that
  \begin{equation*}
    u = \tilde{g}_{+}v_{+\imath 0}^{1-\imath\alpha} +
    \tilde{g}_{-}v_{-\imath 0}^{1-\imath\alpha} + u',
  \end{equation*}
  with $u' \in I^{(\min(s_{0},\tilde s)+1-\delta)}(N^{*}Y)$ for all $\delta >0$. Now
the first two terms are in $\bigcap _{\delta > 0} I^{(3/2 +
    \im\alpha-\delta)}(N^{*}Y)$ but not in $I^{(3/2 + \im
    \alpha)}(N^{*}Y)$ unless $\tilde{g}_{\pm}$ vanish; by the
  assumption that $s_{0}<\tilde s+1$ (so $(\min(s_{0},\tilde
  s)+1-\delta)>s_0$ for some $\delta>0$), we find $3/2 + \im \alpha = s_{0}$ and
  $\tilde{g}_{\pm}$ cannot both vanish. Proceeding as above,
  $P(u-u_{+}-u_{-})\in I^{(\tilde s)}(N^*Y)$, $u-u_{+}-u_{-}\in
  \bigcap _{\delta > 0}I^{(\min(3/2 + \im\alpha,\tilde s)+1 -
    \delta)}(N^{*}Y)$. If $\min(3/2 + \im\alpha,\tilde s)=\tilde s$,
  we are done, otherwise $3/2+\im\alpha<\tilde s$, $u-u_{+}-u_{-}\in
  \bigcap _{\delta > 0}I^{(5/2 + \im\alpha -
    \delta)}(N^{*}Y)$ so by \cite[Theorem~6.3]{Haber-Vasy:Radial}
  $u-u_{+}-u_{-}\in \bigcap_{\delta>0}I^{(\tilde s+1-\delta)}(N^*Y)$,
 completing the proof.
\end{proof}

In our setting, where by equation~\eqref{eq:Psigma} $$\alpha = \sigma -
\imath,$$ this gives:
\begin{corollary}
  \label{cor:Mellin-wave-asymp}
  If $u\in I^{(-\infty)}(N^{*}Y)$ and $P_{\sigma}u \in \CI(X)$, then there exist $g_{\pm}\in \CI(X)$ and $\tilde{u}\in
  C^{\infty}(X)$ such that
  \begin{equation}
    \label{eq:resonantstate}
    u=g_{+}v_{+\imath 0}^{-\imath\sigma} + g_{-} v_{-\imath 0}^{-\imath\sigma} + \tilde{u}.
  \end{equation}
  Again, see Remark~\ref{rem:precise-expansion} if $-\imath\sigma$ is
  an integer.

If instead $P_{\sigma} u\in I^{(\tilde s)}(N^*Y)$, the same conclusion holds with
$\tilde{u}\in \CI(X)$  replaced by $\tilde u\in \bigcap_{\delta>0} I^{(\tilde s+1-
\delta)}(N^*Y)$.
\end{corollary}

Note that $u$ as in the corollary lies in $H^{1/2 + \im\sigma -
  \epsilon}$ for all $\epsilon > 0$, but not in $H^{1/2 + \im\sigma}$
unless $g_{\pm }|_{Y}$ vanish.  Thus, for $s$ and $\sigma$
corresponding to the region~\eqref{eq:P-sigma-holomorphic}, this lies
in $H^{s}$, the target space of $(P_{\sigma})^{-1}_{\tow}$, as
expected---and this containment is sharp insofar as it would fail
whenever $g_{\pm} |_{Y}$ do not vanish if the inequality
in~\eqref{eq:P-sigma-holomorphic} is
replaced by equality.

Finally, we now use Lemma~\ref{lemma:conormalinhomog} to deduce the
structure of solutions to certain inhomogeneous equations with conormal
right hand side:
\begin{proposition}\label{prop:inverseofpsigma}
  If $u\in I^{(-\infty)}(N^{*}Y)$ and $$P_{\sigma}u \in v_{\pm \imath
    0}^{-\imath\sigma-1} \log (v\pm \imath 0)^m\CI(X),$$ then there
  exist $g_{\pm}^{m'}\in \CI(X)$ (for $m'=0,\dots m+1$) and
  $\tilde{u}\in C^{\infty}(X)$ such that
  \begin{equation}
    \label{eq:resonantstate2}
    u=\sum_{m'=0}^{m+1} g_\pm^{m'} v_{\pm\imath 0}^{-\imath\sigma} \log (v\pm \imath 0)^{m'}
 + \tilde{u}.
  \end{equation}
  See Remark~\ref{rem:precise-expansion} and
  Proposition~\ref{lemma:conormalinhomog} if $-\imath\sigma$ is
  an integer.
\end{proposition}

\begin{proof}
  By Lemma~\ref{lemma:conormalinhomog}, we may find a function
  $w$ of the form~\eqref{eq:resonantstate2} so that
  \begin{equation*}
    P_{\sigma}w - P_{\sigma}u \in \CI(X),  
  \end{equation*}
  with the leading term having the claimed form.  As the function $w$
  is also conormal, $w-u$ is conormal, and so we may apply
  Corollary~\ref{cor:Mellin-wave-asymp} to finish the proof.
\end{proof}

\section{The connection between $P_\sigma$ and asymptotically
  hyperbolic and de Sitter spaces}
\label{sec:ident-p_sigma-1}

In this section we clarify the action of $P_{\sigma}^{-1}$ on the caps
$C_{\pm}$ and in the equatorial region $C_0$ as in \cite[Sections 3.3
and 4]{Vasy-Dyatlov:Microlocal-Kerr}.  Recall that $P_{\sigma}^{-1}$
propagates regularity from $S_{-}$ to $S_{+}$; by contrast the
behavior at $C_-$ and $C_0$ is what is studied in detail in
\cite{Vasy-Dyatlov:Microlocal-Kerr}, with the behavior of
$P_{\sigma}^{-1}$ at $C_0$ and $C_+$ corresponding to the adjoint
operator in that paper.

On $C_\pm$ we consider the operators
  \begin{equation*}
\begin{aligned}
    L_{\sigma,\pm} &=
v^{\frac{1}{2}}v^{\frac{n}{4}\pm\frac{\imath\sigma}{2}}P_{\pm\sigma}v^{-\frac{n}{4}
    \mp \frac{\imath\sigma}{2}}v^{\frac{1}{2}}\\
&= - \Lap_{k_{\pm}} + \left( \sigma^{2} +
      \frac{(n-2)^{2}}{4}\right) + v\mathcal{X}(\pm\sigma) +
    vV(\pm\sigma),
\end{aligned}
  \end{equation*}
from \eqref{eq:relat-with-asymp-hyp}  (note the sign switch in $\sigma$
relative to \eqref{eq:relat-with-asymp-hyp} to keep the behavior for
$L_{\sigma,+}$ and $L_{\sigma,-}$ similar in terms of $\im\sigma>0$
being the physical half-plane).  Here
 $k_{\pm}$ are asymptotically hyperbolic metrics, $V$ a smooth potential
  and $\mathcal{X}$ a vector field tangent to $v=0.$  On $C_0$ we
  likewise consider
\begin{equation*}
 L_{\sigma,0}= \Box_{k_{0}} - \left(
    \sigma^{2} + \frac{(n-2)^{2}}{4}\right) + v \mathcal{X}(\sigma) + vV(\sigma)
\end{equation*}
from \eqref{eq:connection-with-dS}, with $V$, $\mathcal{X}$ as above
($|v|=-v$ being a defining function for $\overline{C_0}$).

Since $L_{\sigma,0}$ is an asymptotically de Sitter operator as in \cite{Vasy:asymp-dS}, it has a
forward solution operator $\cR_{C_0}(\sigma)$ propagating towards
$S_+$, i.e.\ if $f\in\CI_c(C_0)$,  $u=\cR_{C_0}(\sigma)f$ is the
unique solution of $L_{\sigma,0}u=f$ with $L_0$ vanishing near
$S_-$. On the other hand, $L_{\sigma, \pm}$ are non-self-adjoint
perturbations of  the asymptotically hyperbolic operator $-\Lap_{k_{\pm}} + ( \sigma^{2} +
      \frac{(n-2)^{2}}{4})$, as in \cite{Mazzeo-Melrose},
with the perturbation being non-trapping in the high energy sense.
In particular, if we let $H^2_0$ denote the \emph{$0$-Sobolev space}
associated to the $0$-calculus of \cite{Mazzeo-Melrose}, then
$L_{\sigma,\pm}:H^2_0(\overline{C_\pm})\to L^2(\overline{C_\pm})$,
$\im\sigma>0$, is an analytic Fredholm family.  A priori we do
not have automatic invertibility for such perturbations without appeal
to the large parameter behavior, which is only understood from the
perspective of the extended operator; we prove this below. Note that
if $ v \mathcal{X}(\sigma) + vV(\sigma)$ vanishes then the
invertibility of $L_{\sigma}$ is automatic when $\im\sigma>0$,
  $\im\sigma^2\neq 0$ as $\Lap_{k_{\pm}}$ is self-adjoint.

We begin by recording a result on supports that follows from the proof
of Proposition~3.9 of \cite{Vasy-Dyatlov:Microlocal-Kerr} (with the
complex absorption hypotheses employed there irrelevant
here).
\begin{lemma}\label{lemma:uppertriangular}
Suppose $P_\sigma^{-1}:\cY^{s_{\tow}-1}\to\cX^{s_{\tow}}$ is regular at
$\sigma\in \CC$ with $\im \sigma>0.$
If $f\in \cY^{s_{\tow}-1}$ and $\supp f \subset \overline{C_+\cup C_0}$ then $\supp P_\sigma^{-1} f \subset \overline{C_+\cup C_0}.$
\end{lemma}

Although the proof of Lemma~\ref{lemma:uppertriangular} is essentially
contained in \cite[Proposition~3.9]{Vasy-Dyatlov:Microlocal-Kerr}, we
include a sketch of the proof here.
\begin{proof}
  The microlocal argument proving the Fredholm estimates in
  Section~\ref{sec:results-citevasy} in fact yield a microlocal
  version of the same estimates.  In particular, let
  $\mathbf{t}$ be a global function that is time-like on $C_{0}$, with
  $\mathbf{t}^{-1}(T_{0}, T_{1})\subset C_{0}$.  As we already know
  $P_{\sigma}^{-1}f$ is smooth near $S_{-}$, we may estimate
  \begin{equation*}
    \Norm[\cX^{s_{\tow}}(\mathbf{t}^{-1}(-\infty, T_{0}))]{u} \leq C
    \left( |\Re \sigma |
      \Norm[\cY^{s_{\tow}-1}(\mathbf{t}^{-1}(-\infty,
      T_{1}))]{P_{\sigma}u} + |\Re \sigma|^{-1}
      \Norm[H^{-N}(\mathbf{t}^{-1}(-\infty, T_{1}))]{u}\right) .
  \end{equation*}

  As $P_{\sigma}$ is hyperbolic in $C_{0}$, energy estimates allow
  allow us to estimate $\Norm[\cX^{s_{\tow}}]{u}$ on
  $\mathbf{t}^{-1}(-\infty, T_{1})$ in terms of the same right hand
  side.  For $|\Re \sigma|$ large enough, we may then absorb the
  second term on the right side into the left side.  In particular, if
  $T_{1}$ is such that $\mathbf{t}^{-1}(-\infty, T_{1})$ is disjoint
  from $\supp f$, then $u$ must vanish on $\mathbf{t}^{-1}(-\infty,
  T_{1})$.  

  Meromorphic continuation then shows that the same support property
  holds at all $\sigma$ so that $P_{\sigma}^{-1}$ is regular.
\end{proof}

We now prove the following lemma relating invertibility of $P_\sigma$
and $L_{\sigma,\pm}:$
\begin{lemma}\label{lemma:Lpm}
Suppose that $P_\sigma:\cX^{s_{\tow}}\to\cY^{s_{\tow}-1}$ and
$P_\sigma^*:\cX^{s_{\away}^*}\to\cY^{s_{\away}^*-1}$
are invertible at
$\sigma\in\CC$ with $\im\sigma>0$.
Then $L_{\sigma,\pm}:H^2_0(\overline{C_\pm})\to L^2(\overline{C_\pm})$
is invertible.
\end{lemma}

We let $\mathcal{R}_{C_\pm}(\sigma)$ denote the inverse of
$L_{\sigma,\pm}$ thus obtained.

\begin{remark}
While we handle the invertibility within our framework, an alternative
would be the complex absorption framework used in
\cite{Vasy:Hyperbolic}; the absorption would be placed in $v<-\ep$
for some $\ep>0$.
\end{remark}

\begin{remark}
  We may verify the invertibility hypothesis above by employing
  Proposition~\ref{prop:nontrappingestimates}.
\end{remark}

\begin{proof}
As already remarked, $L_{\sigma,\pm}:H^2_0(\overline{C_\pm})\to
L^2(\overline{C_\pm})$
is Fredholm, so we only need to show that $\Ker L_{\sigma,\pm}$ and
$\Ker L_{\sigma,\pm}^*$ are trivial.
By the results of
\cite{Mazzeo-Melrose}, first any element of $\Ker
L_{\sigma,\pm}$ is in $H^\infty_0(\overline{C_\pm})$ by elliptic
regularity in the 0-calculus, and indeed using the parametrix construction,
of \cite{Mazzeo-Melrose} they are in $v^{-\frac{1}{2}+\frac{n}{4} +
    \frac{\imath\sigma}{2}}\CI(\overline{C_\pm})$, while any
  element of $\Ker
L_{\sigma,\pm}^*$ is in $v^{-\frac{1}{2}+\frac{n}{4} -
    \frac{\imath\overline{\sigma}}{2}}\CI(\overline{C_\pm})$. In
  particular, for any element $u_-$ of
  the kernel of $L_{\sigma,-},$ we can extend $v^{\frac{1}{2}-\frac{n}{4} -
    \frac{\imath\overline{\sigma}_{0}}{2}} u_-$ to an element $\tilde
  u$ of $\CI(X)$. Then $f=P_{\sigma} \tilde u$ is supported in
  $\overline{C_+\cup C_0}$ by \eqref{eq:relat-with-asymp-hyp}, hence
  by Lemma~\ref{lemma:uppertriangular}
  $P_{\sigma}^{-1}f$ is also supported in $\overline{C_+\cup C_0}$, so
  $u=\tilde u-P_\sigma^{-1}f$ solves $P_{\sigma} u=0$ and
  $u|_{C_-}=u_-$. Since $\Ker P_{\sigma}$ is trivial by assumption,
  $u$, and thus $u_-$, vanish.
A similar argument applies to elements of $\Ker L_{\sigma,+}^*$
as $P_\sigma^*:\cX^{s_{\away}^*}\to\cY^{s_{\away}^*-1}$ is also invertible;
in that case for an element $u_+$ of the kernel $v^{\frac{1}{2}-\frac{n}{4} +
    \frac{\imath\overline{\sigma}_{0}}{2}}u_+$ to an element $u$ of
  $\CI(X)$ and apply $(P_{\sigma}^*)^{-1}$ to the result. Finally,
  for $\Ker L_{\sigma,+}$ and $\Ker L^*_{\sigma,-}$ we switch the
  direction of propagation for the inverse $P_{\sigma}^{-1}$, i.e.\ we consider
$$
P_\sigma:\cX^{s_{\away}}\to\cY^{s_{\away}-1},\ P_\sigma^*:\cX^{s_{\tow}^*}\to\cY^{s_{\tow}^*-1},
$$
and then completely analogous arguments apply as the roles of $C_+$
and $C_-$ are simply reversed.
\end{proof}

We now make the connection between $P_{\sigma}^{-1}$ and the operators
on the $C_{\pm}$ and $C_{0}$.  Let $-\sigma_0,$ be a point with $\im \sigma_0<0$ 
at which the family $P_{\sigma}^{-1}$ is regular.  Then by  Lemma~\ref{lemma:Lpm} 
$\mathcal{R}_{C_{-}}(.)$ is likewise regular at $-\sigma_0.$
If $f\in
C^{\infty}_{c}(C_{-})\subset\CI(X)$, then $P_{\sigma_{0}}^{-1}f$ is smooth on
$\overline{C_-}$ by Proposition~\ref{proposition:fredholm}.  By
\eqref{eq:relat-with-asymp-hyp}, on $C_{-}$, 
$$
u=v^{-\frac{1}{2}+\frac{n}{4} +
    \frac{\imath\sigma_{0}}{2}} \big(P_{\sigma_{0}}^{-1}(v^{\frac{n}{4} + \frac{\imath\sigma_{0}}{2} + \frac{1}{2}}f)\big)|_{C_-} 
$$
solves $L_{-\sigma_0,-}u = f$, and  $u\in v^{-\frac{1}{2}+\frac{n}{4} +
    \frac{\imath\sigma_{0}}{2}}\CI(\overline{C_-})\subset L^2_0(\overline{C_-})$ (with
  $L^2_0(\overline{C_-})$ being the asymptotically hyperbolic $L^2$
  space as described in \cite{Mazzeo-Melrose}), with the inclusion holding as
$\re (\imath\sigma_{0})=-\im\sigma_0>0$. Thus,
\begin{equation*}
  v^{-\frac{1}{2}+\frac{n}{4} +
    \frac{\imath\sigma_{0}}{2}} \big(P_{\sigma_{0}}^{-1}(v^{\frac{n}{4} + \frac{\imath\sigma_{0}}{2} + \frac{1}{2}}f)\big)|_{C_-} = \mathcal{R}_{C_{-}}(-\sigma_{0})f,
\end{equation*}
since $\mathcal{R}_{C_{-}}(-\sigma_{0})f$ is the unique $L^2$
(relative to the asymptotically hyperbolic metric) solution of
$L_{-\sigma_0,-}u = f$. By the meromorphy
  of both sides, the formula is then valid at all $\sigma$ (regardless of
  the sign of $\im\sigma$) at which $P_{\sigma}^{-1}$ is regular.
Indeed, by the same argument, for $f\in\CI(C_-)$ such that
$v^{\frac{n}{4} + \frac{\imath\sigma_{0}}{2} + \frac{1}{2}}
f\in\CI(\overline{C_-})$, and thus has an extension $\tilde f=E(v^{\frac{n}{4} + \frac{\imath\sigma_{0}}{2} + \frac{1}{2}}
f)$ to an
element of $\CI(X)$, one still has
\begin{equation*}
  v^{-\frac{1}{2}+\frac{n}{4} +
    \frac{\imath\sigma_{0}}{2}} \big(P_{\sigma_{0}}^{-1}E(v^{\frac{n}{4} + \frac{\imath\sigma_{0}}{2} + \frac{1}{2}}f)\big)|_{C_-} = \mathcal{R}_{C_{-}}(-\sigma_{0})f.
\end{equation*}

Recall that if $f\in \cY^{s_{\tow}-1}$ and $\supp f \subset
\overline{C_+\cup C_0}$ then $\supp P_\sigma^{-1} f \subset
\overline{C_+\cup C_0}$ by Lemma~\ref{lemma:uppertriangular}.  Turning to the region $C_{0}$, the
Carleman-type estimates in \cite[Proposition 5.3]{Vasy:asymp-dS} (see
also \cite[Section 4]{Vasy-Dyatlov:Microlocal-Kerr}) imply that
$P_{\sigma_{0}}^{-1}f$ must vanish in a neighborhood of $S_{-}$.  In
particular, $P_{\sigma_{0}}^{-1}$ (applied to such $f$, with the
result restricted to $C_0$) must be a conjugate of the forward
fundamental solution of the operator in
equation~\eqref{eq:connection-with-dS} (applied to $f|_{C_0}$),
denoted $\mathcal{R}_{C_{0}}(\sigma)$ above. Indeed, again a simple
generalization shows the same conclusion when one merely has $f\in
\dCI(\overline{C_{0}}\cup\overline{C_+})$, with the dot denoting
infinite order vanishing at the boundary of this set, namely $S_-$.

Finally, if $f\in C^{\infty}_{c}(C_{+})$, or indeed $f\in\dCI(C_+)$ then the above discussion
implies that $P_{\sigma_{0}}^{-1}f$ vanishes in $\overline{C}_{-}$ and
$C_{0}$.  Moreover, the expansion of Corollary~\ref{cor:Mellin-wave-asymp} implies that in fact
$P_{\sigma_{0}}^{-1}$ is a conjugate of $\mathcal{R}_{C_{+}}(\sigma)$.

In particular, the above discussion proves the following proposition:
\begin{proposition}
  \label{prop:block-structure}
  Let $\varpi_\pm=1/2\pm n/4\pm\imath\sigma/2.$
If $P_{\sigma}^{-1}$ is regular at $\sigma$, then it
  has the following ``block structure'' (here the rows and columns
  correspond to support in $C_{+}$, $C_{0}$, and $C_{-}$):
  \begin{equation*}
\begin{pmatrix}
          |v|^{\varpi_-}
\mathcal{R}_{C_{+}}(\sigma) |v|^{\varpi_+}
 & * & * \\ 0 &
          |v|^{\varpi_-}
\mathcal{R}_{C_{0}}(\sigma) |v|^{\varpi_+}
& * \\ 0 & 0 &     |v|^{\varpi_-}
\mathcal{R}_{C_{-}}(-\sigma) |v|^{\varpi_+}
    \end{pmatrix}
 \end{equation*}
in the strong sense that if $P_\sigma^{-1}$ is applied to a $\CI$
function on $X$, the restriction of the result to $C_-$ is given by
the lower right block,
$$|v|^{\frac{1}{2}-\frac{n}{4}-\frac{\imath\sigma}{2}}
\mathcal{R}_{C_{-}}(-\sigma)
|v|^{\frac{1}{2}+\frac{n}{4}+\frac{\imath\sigma}{2}},
$$
if
$P_\sigma^{-1}$ is applied to a $\CI$ function supported in
$\overline{C_0\cup C_+}$,  the restriction of the result to $C_0$ is
given by $|v|^{\frac{1}{2}-\frac{n}{4}-\frac{\imath\sigma}{2}}
\mathcal{R}_{C_{0}}(\sigma)
|v|^{\frac{1}{2}+\frac{n}{4}+\frac{\imath\sigma}{2}}$ (and this result
vanishes in $\overline{C_-}$), while finally
if $P_\sigma^{-1}$ is applied to a $\CI$ function supported in
$\overline{C_+}$,  the restriction of the result to $C_+$ is
given by $|v|^{\frac{1}{2}-\frac{n}{4}-\frac{\imath\sigma}{2}}
\mathcal{R}_{C_{+}}(\sigma)
|v|^{\frac{1}{2}+\frac{n}{4}+\frac{\imath\sigma}{2}}$ (and the result
vanishes in $\overline{C_-}\cup C_0$).
\end{proposition}

By our non-trapping assumption on the null-geodesics of $g$, $-\Lap
_{k_{\pm}} + v\mathcal{X}(\sigma) + vV(\sigma) + \sigma ^{2} + (n-2)^{2}/4$ is a
semi-classically non-trapping operator and thus the following
proposition (which follows from, e.g., the work of
Vasy~\cite{Vasy-Dyatlov:Microlocal-Kerr}) applies.

\begin{proposition}[cf.\ {\cite[Theorem 4.7]{Vasy:Hyperbolic}}]
  \label{prop:meromorphic-hyperbolic}
  Consider the operators
  \begin{equation*}
    L_{\sigma,\pm} = - \Lap_{k_{\pm}} + \left( \sigma^{2} +
      \frac{(n-2)^{2}}{4}\right) + v\mathcal{X}(\pm\sigma) + vV(\pm\sigma),
  \end{equation*}
  with $k_{\pm}$ asymptotically hyperbolic metrics, $V$ a smooth potential
  and $\mathcal{X}$ a vector field tangent to $v=0$.  If $L_{\sigma,\pm}$
  is semiclassically non-trapping, then it has a meromorphic inverse
  $\mathcal{R}_{C_\pm}(\sigma)$ with finite rank poles, is holomorphic for
  $\Im \sigma \gg 0$, and has only finitely many
  poles in each strip $C_{1} \leq \Im \sigma \leq C_{2}$.  Moreover,
  non-trapping estimates hold in each strip $\Im \sigma > -C$ for
  large $\Re \sigma$:
  \begin{equation*}
    \Norm[H^{s}_{|\sigma|^{-1}}]{\mathcal{R}_{C_\pm}(\sigma)f} \leq C \Norm[H^{s-1}_{|\sigma|^{-1}}]{f}
  \end{equation*}
  
  Moreover, if $L_{\sigma, \pm}$ has no $L^{2}$ kernel (with
  respect to the metric $k_{\pm}$) then all poles $\sigma_{0}$ of
  $\mathcal{R}_{C_\pm}(\sigma)$ have $\Im \sigma_{0} \leq 0$.
\end{proposition} 

\begin{proof}
The bounded $\sigma$ properties were already explained above.
The high energy estimates
    then follow from those for $P_\sigma^{-1}$. Since $P_\sigma$ has
    index zero, its invertibility amounts to
    having a trivial kernel. Since an element of $\Ker P_\sigma$
    restricts to a $\CI$ function on $\overline{C_-}\cup C_0$, thus
    $v^{\frac{n}{4}-\frac{1}{2}+\imath\frac{\sigma}{2}}$ times the
    restriction to $C_-$ is an element of
    $v^{\frac{n}{4}-\frac{1}{2}+\imath\frac{\sigma}{2}}\CI(\overline{C_-})$,
    in view of the asymptotically hyperbolic metric on $C_-$
    this gives an element of $L^2$ if $\im\sigma<0$.
Thus under the assumption of no $L^2$ ``eigenvalues'' all poles $\sigma_0$ of
$\mathcal{R}_{C_-}(\sigma)$ indeed
have $\Im \sigma_{0} \leq 0$. On the other hand, by Section~\ref{sec:conormal-reg}, an
element of $\Ker P_\sigma$ whose support is disjoint from $C_-$ is
supported in $\overline{C_+}$, and restricted to $C_+$ it has the asymptotic form
$v^{-\imath\sigma}\CI(\overline{C_+})$, and thus $v^{\frac{n}{4}-\frac{1}{2}+\imath\frac{\sigma}{2}}$ times the
    restriction to $C_+$ is an element of
    $v^{\frac{n}{4}-\frac{1}{2}-\imath\frac{\sigma}{2}}\CI(\overline{C_-})$,
    and thus is in $L^2$ if $\im\sigma>0$, under the assumption of no $L^2$ ``eigenvalues'' all poles $\sigma_0$ of
$\mathcal{R}_{C_+}(\sigma)$ indeed
have $\Im \sigma_{0} \geq 0$.
\end{proof}

\begin{remark}
  \label{rem:extraneous-poles}
  Proposition~\ref{prop:block-structure} implies that the poles of
  $\mathcal{R}_{C_{+}}(\sigma)$ yield poles of $P_{\sigma}^{-1}$.  A
  partial converse is true as well.  If $\sigma_{0}$ is a pole of
  $P_{\sigma}^{-1}$ so that the corresponding resonant dual state has
  support intersecting $X\setminus \overline{C}_{-}$, then either
  $\sigma_{0}$ is a pole of $\mathcal{R}_{C_{+}}(\sigma)$ or the
  corresponding resonant state is supported at $S_{+}$ (see
  \cite[Remark 4.6]{Vasy-Dyatlov:Microlocal-Kerr} for more details).
  Such poles may occur only for $\sigma_{0}$ a pure imaginary negative
  integer.  In other words, the relevant poles of $P_{\sigma}^{-1}$
  are either poles of $\mathcal{R}_{C_{+}}(\sigma)$ or have state
  supported at $S_{+}$ (and hence are differentiated delta functions
  in $v$).  We remark that such states occur in even-dimensional
  Minkowski space, where $-\imath$ is a pole of $P_{\sigma}^{-1}$ in
  $2$- and $4$-dimensions.
\end{remark}

\section{Structure of the poles of $P_\sigma^{-1}$}
\label{sec:logs}

While the results in the previous section fully address the structure
of nullspace of $P_\sigma,$
knowledge of nullspace alone is clearly not sufficient to deal with the structure of the poles of
$P_\sigma^{-1}.$  Even for a spectral family of the form
$$
(P_0-\sigma \Id)^{-1},
$$ 
with $P_0$ as in \eqref{modelform},
the poles may of course be multiple owing to generalized eigenspaces;
thus knowing that the nullspace of $P_0$ has a particular conormal form
$v^\gamma$ would in general permit the range of the polar part of the
resolvent to have log terms.  Here the situation is further
complicated by the fact that our family $P_\sigma$ is not of the form
$P_0-\sigma \Id$ but rather has nontrivial dependence on $\sigma,$ so
that we cannot even employ the usual machinery of Jordan
decomposition.  A careful analysis of the log terms will, however, be
essential in order to see that excess log terms in our asymptotic
expansion \eqref{asymp} do not spoil the restriction of the rescaled
solution to the front face of the radiation field blowup, which we
know a priori must be smooth if we impose an additional support
hypothesis in $s=v/\rho$ (cf.\ Section~\ref{sec:blowup}).
In this section we demonstrate (among other things) that
the top-order terms with $\log \rho$ are balanced by terms containing
$\log v$ in such a way as to permit the solution to be smooth across
the front face.  (In light of the smoothness of the solution
  across the front face, we expect all such $\log$ terms to be
  balanced in this manner, but we consider only the top-order terms,
  i.e., the terms affecting the radiation field, here.)   We should
emphasize that these log terms are typically \emph{not} vanishing, and
are still a relevant part of the expansion away from the interior of
the front face.  In particular, we prove the following proposition,
which is an extension of Corollary~\ref{cor:Mellin-wave-asymp}:
\begin{proposition}\label{prop:boundarylogterms}
Let $\sigma_0$ be a pole of order $k$ of the operator family
$$
P_\sigma^{-1}: \cY^{s_{\tow}-1}\to \cX^{s_{\tow}},
$$
and let
$$
 (\sigma - \sigma_{0})^{-k}A_{k} +
    (\sigma-\sigma_{0})^{-k+1}A_{k-1} + \ldots +
    (\sigma-\sigma_{0})^{-1}A_{1} + A_{0}
$$
denote the Laurent expansion near $\sigma_0,$ with $A_0$ (locally)
holomorphic.  
If $f \in \cY^{s_{\tow}-1}$ vanishes in a neighborhood of
$\overline{C}_{-}$, there are smooth functions $\phi _{\pm,1},
  \ldots \phi_{\pm,k}$ so that for $0 \leq \ell \leq k-1$, 
 $A_{k-\ell}f$
  has an asymptotic expansion near $S_+:$
  \begin{align*}
    A_{k-\ell} f &= v_{+\imath 0}^{-\imath \sigma_{0}}
    \left[ \sum_{j = 0}^{\ell} \frac{(-\imath)^{j}}{j!}(\log (v+\imath
      0))^{j}\phi_{+,k-(\ell-j)}\right] \\
    &\quad\quad +  v_{-\imath 0}^{-\imath \sigma_{0}}
    \left[ \sum_{j = 0}^{\ell} \frac{(-\imath)^{j}}{j!}(\log (v-\imath
      0))^{j}\phi_{-,k-(\ell-j)}\right]
     + O\left(v^{-\imath\sigma_{0}+1}(\log v)^{\ell}\right).
  \end{align*}
 (Although we use the notation
    $O(v^{\gamma}(\log v)^{\kappa})$ here, the term in fact has a
    polyhomogeneous expansion with index sets shifted from the
    ``base'' ones.) 

  If $-\imath\sigma_{0}$ is a non-negative integer, then there are
  smooth functions $\phi_{1}, \ldots, \phi_{k}$ so that $A_{k-\ell}f$
  has a similar expansion in terms of the distributions
  $v_{+}^{-\imath\sigma_{0}} = H(v)v^{-\imath\sigma_{0}}$:
  \begin{align*}
    A_{k-\ell}f = v_{+}^{-\imath\sigma_{0}} \sum_{j=0}^{\ell}
    \frac{(-\imath )^{j}}{j!} (\log |v|)^{j} \phi_{k-(\ell -j)} +
    O(v^{-\imath\sigma_{0} +1}(\log |v|)^{\ell}).
  \end{align*}
\end{proposition}

\begin{remark}
  \label{rem:purposes-of-prop}
  This proposition serves two purposes.  The first is to show that
  Laurent coefficients have asymptotic expansions at $v=0$, while the
  second is to show that the leading terms in this expansion have a
  specific form.  This form is later used to show that the terms of
  the form $\log \rho$ cancel at the radiation field face so that the
  radiation field may be defined.  

  The additional logarithmic terms occurring at imaginary integers in
  Proposition~\ref{lemma:conormalinhomog} would in general disrupt the
  form of this expansion, but we use the support of the states to
  conclude that in fact it has the desired form.

  One could also write the entire expansion in terms of
  $H(v)v^{-\imath\sigma_{0}}$ even if $-\imath\sigma_{0}$ is not a
  positive integer.  To do this, we would have to include derivatives
  of delta functions if $-\imath\sigma_{0}$ is a negative integer.
  (As noted in Remark~\ref{rem:extraneous-poles} these occur even in
  the case of even dimensional Minkowski space.)
\end{remark}

The proof requires the following lemma:
\begin{lemma}
  \label{lem:support-of-polar-part}
  If $f$ vanishes in a neighborhood of $\overline{C}_{-}$ then
  $A_{k-\ell}f$ is supported in $\overline{C}_{+}$ for $\ell = 0, 1,
  \ldots, k-1$.
\end{lemma}

\begin{remark}
  This lemma implies that there are two types of ``resonant'' states.
  If the state is given by $\phi \langle \psi, \cdot \rangle$, then
  either $\phi$ is supported in $\overline{C}_{+}$ or $\psi$ is
  supported in $\overline{C}_{-}$. See
  \cite[Section~4.9]{Vasy-Dyatlov:Microlocal-Kerr}, especially
  Remark~4.6, for more details.
\end{remark}

\begin{proof}
  Near a pole $\sigma_{0}$ of $P_{\sigma}^{-1}$, we may write
  \begin{equation}
    \label{eq:P-near-pole}
    P_{\sigma} = P_{0} + (\sigma - \sigma_{0})P_{1} + (\sigma - \sigma_{0})^{2}P_{2},
  \end{equation}
  where $P_{0} = P_{\sigma_{0}}$, $P_{1} = D_{v} + E$, and $P_{2}$ is
  a smooth function.  Here $E \in \mathcal{M}_{\pa}$ is a first order
  differential operator characteristic on $N^{*}\YS$.  The proof relies
  on the following relationships between $P_{i}$ and $A_{j}$, which
  holds because $P_{\sigma}P_{\sigma}^{-1}= \id$:
  \begin{align}
    \label{eq:relationship-P-and-A}
    &P_{0}A_{k} = 0\\
    &P_{1}A_{k} + P_{0}A_{k-1} = 0 \notag \\ 
    &P_{2}A_{k-i} + P_{1}A_{k-(i-1)} + P_{0}A_{k-(i-2)} = 0, \quad i = 0,
    \ldots, k-3 \notag
  \end{align}
  
  We first observe that $A_{k-\ell}f$ vanishes near
  $\overline{C}_{-}$.  Indeed, for $\ell=0$ this follows from the
  Cauchy integral formula applied to $(\sigma -
  \sigma_{0})^{k-1}P_{\sigma}^{-1}f$ and
  Proposition~\ref{prop:block-structure}, while for $\ell > 0$, it
  follows inductively from Proposition~\ref{prop:block-structure},
  \eqref{eq:relationship-P-and-A}, and the
  Cauchy integral formula applied to
  \begin{equation*}
    (\sigma-\sigma_{0})^{k-\ell-1}P_{\sigma}^{-1}f.
  \end{equation*}

  To observe that $A_{k-\ell}f$ vanishes in $C_{0}$, we again proceed
  inductively.  For $\ell=0$, as $P_{0}A_{k}f = 0$, the proof of
  Proposition~\ref{prop:block-structure} implies that it vanishes in a
  neighborhood of $S_{-}$ and hence in all of $C_{0}$.  (The
  Proposition does not apply as stated only because we are not at a
  regular point of $P_\sigma^{-1},$ but the Carleman and energy
  estimates---the latter being used in the proof of
  Lemma~\ref{lemma:uppertriangular}---employed nonetheless apply here
  as well.)  If $\ell > 0$, the
  relationship~\eqref{eq:relationship-P-and-A} implies that
  $P_{0}A_{k-\ell}f$ vanishes in $C_{0} \cup \overline{C}_{-}$ and so
  $A_{k-\ell}f$ also vanishes in $C_{0} \cup \overline{C}_{-}$.
\end{proof}

\begin{proof}[Proof of Proposition~\ref{prop:boundarylogterms}]
  We rely on the structure of $P_{\sigma}$ near $S_{+}$.  Indeed,
  recall from above that $P_{\sigma} = D_{v}(vD_{v} + \sigma) + Q$,
  where $Q\in M_{\pa}^{2}$ (in the notation of
  Section~\ref{sec:conormal-reg}) is a differential operator.  We rely
  on the form~\eqref{eq:P-near-pole} of $P_{\sigma}$ near a pole
  $\sigma_{0}$ of $P_{\sigma}^{-1}$ as well as the
  relationships~\eqref{eq:relationship-P-and-A} between $P_{i}$ and
  $A_{j}$.

  We start by assuming that $-\imath \sigma_{0}$ is not an integer and
  proceed by induction on $\ell$.  As $f$ vanishes near
  $\overline{C}_{-}$, Lemma~\ref{lem:support-of-polar-part} implies
  that $A_{k}f$ is supported in $\overline{C}_{+}$, while
  Proposition~\ref{proposition:fredholm} (or, indeed, elliptic
  regularity) implies it is smooth away from the radial set
  $\Lambda^{+}$.  We may thus apply a theorem of
  Haber--Vasy~\cite[Theorem 6.3]{Haber-Vasy:Radial} to conclude that
  in fact $A_{k}f\in I^{(-\infty)}(\Lambda^{+}) =
  I^{(-\infty)}(N^{*}\YS)$.  In particular, then,
  Corollary~\ref{cor:Mellin-wave-asymp} implies that there are smooth
  functions $\phi_{\pm,k}$ and $\psi$ so that
  \begin{equation*}
    A_{k}f = v_{+\imath 0}^{-\imath\sigma_{0}}\phi_{+,k} + v_{-\imath
      0}^{-\imath \sigma_{0}} \phi_{-,k} + \psi .
  \end{equation*}
  By Lemma~\ref{lem:support-of-polar-part}, $A_{k}f$ is supported in
  $\overline{C}_{+}$, so $\psi$ vanishes to infinite order at $S_{+}$
  and may be absorbed into the other terms, i.e.,
  \begin{equation*}
    A_{k}f = v_{+\imath 0}^{-\imath\sigma_{0}}\phi_{+,k} + v_{-\imath
      0}^{-\imath \sigma_{0}} \phi_{-,k}.
  \end{equation*}
  
  Now suppose that the statement is true for $0 \leq \ell' \leq \ell
  -1$.  As $P_{0}A_{k-\ell}f = - P_{1}A_{k-\ell + 1}f -
  P_{2}A_{k-\ell+2}f$, we have
  \begin{align*}
    P_{0}A_{k-\ell}f &= 
      \sum_{\pm}\sum_{j=0}^{\ell-1} v^{-\imath\sigma_{0}
        -1}_{\pm\imath 0} (\log (v\pm \imath 0))^{j}
      \frac{(-\imath)^{j}}{j!}\left[\sigma_{0} \phi_{\pm, k -
          (\ell-1-j)} + \phi_{\pm, k - (\ell -2 -j)} \right] \\
      &\quad\quad\quad+
      O(v^{-\imath \sigma_{0}} (\log v)^{\ell}),
  \end{align*}
  where the $O(v^{-\imath\sigma_{0}}(\log v)^{\ell})$ in fact has an
  asymptotic expansion of a similar form.  Observe that the right hand
  side is an element of $I^{(-\infty)}(N^{*}\YS)$, so again
  Haber--Vasy \cite{Haber-Vasy:Radial}
  implies that $A_{k-\ell}f\in I^{(-\infty)}(N^{*}\YS)$.
  Proposition~\ref{prop:inverseofpsigma} then implies that
  $A_{k-\ell}f$ has a similar expansion, say (suppressing dependence
  of coefficients on $k,\ell$)
  \begin{align*}
    A_{k-\ell}f = \sum_{\pm}\sum_{j=0}^{\ell}
    v_{\pm\imath 0}^{-\imath\sigma_{0}}(\log (v\pm \imath 0))^{j} a_{\pm, j}
  \end{align*}
  To determine the leading coefficients in the expansion, we calculate
  \begin{align*}
    P_{0}A_{k-\ell}f &= -\sum_{\pm}\sum_{j=0}^{\ell-1} v_{\pm \imath 0}^{-\imath
      \sigma_{0} -1}(\log (v\pm\imath 0))^{j}(j+1)\left[ (-\imath
      \sigma_{0}) a_{\pm,j+1} + (j+2)a_{\pm,j+2}\right] \\
    &\quad \quad \quad + O(v^{-\imath\sigma_{0}}(\log v)^{\ell}),
  \end{align*}
  where again the last term has an expansion.  We now simply equate
  coefficients, starting with the largest one.  If $j=\ell-1$, we must
  have
  \begin{equation*}
    \imath\sigma_{0}\ell a_{\pm, \ell} = \frac{(-\imath)^{\ell-1}}{(\ell-1)!} \sigma_{0} \phi_{\pm, k},
  \end{equation*}
  i.e., $$a_{\pm, \ell} = \frac{(-\imath)^{\ell}}{\ell!} \phi_{\pm, k}.$$  Now
  for $j<\ell-1$, we have
  \begin{align*}
    &\imath \sigma_{0}(j+1)a_{\pm, j+1} -
    (j+2)(j+1)\frac{(-\imath)^{j+2}}{(j+2)!}\phi_{\pm,k-(\ell-2-j)} \\
    & \quad\quad=
    \frac{(-\imath)^{j}}{j!}\left( \sigma_{0}\phi_{\pm,k-(\ell-1-j)} + \phi_{\pm,k-(\ell-2-j)}\right),
  \end{align*}
  i.e., $$a_{\pm,j+1} = \frac{(-\imath)^{j+1}}{(j+1)!} \phi_{\pm,
    k-(\ell-1-j)}.$$  This determines $a_{\pm,2}, \ldots,
  a_{\pm,\ell}$, while $a_{\pm,1}$ are given by
  Corollary~\ref{cor:Mellin-wave-asymp} and are denoted $\phi_{\pm, k-\ell}$.
  
  We now consider when $-\imath\sigma_{0}$ is an integer, in which
  case additional logarithmic terms appear in
  our application of Proposition~\ref{prop:inverseofpsigma}.  If $-\imath\sigma_{0}< 0$,
  these additional logarithms are not in the leading order terms and
  so the results above still hold.  For $-\imath \sigma_{0} \geq 0$ an
  integer, however, we must be a bit more careful and rely on
  Lemma~\ref{lem:support-of-polar-part} as follows.

  Let us assume for now that $-\imath\sigma_{0} \neq 0$.  Indeed, we
  again proceed inductively.  Consider first $A_{k}f$.  The same
  arguments as above imply that $A_{k}f$ has an expansion of the form
  \begin{equation*}
    A_{k}f = v_{+\imath 0}^{-\imath\sigma_{0}}\phi_{+} + v_{-\imath
      0}^{-\imath\sigma_{0}} \phi_{-} +  \phi.
  \end{equation*}
  As $A_{k}f$ is supported in $\overline{C}_{+}$ and $v_{\pm\imath
    0}^{-\imath\sigma_{0}} = v^{-\imath\sigma_{0}} \log (v\pm\imath 0)$
  in this case, given $\phi_{+}$, the behavior of $\phi_{-}$ and
  $\phi$ at $\YS$ is determined by the support condition.  Indeed, we
  must have $\phi_{-} = -\phi_{+}$ and $\phi = -2\pi \imath
  v^{-\imath\sigma_{0}}\phi_{+}$.  In other words, there is a smooth
  function $\phi_{k}$ so that
  \begin{equation*}
    A_{k}f = v^{-\imath\sigma_{0}}H(v)\phi_{k}
  \end{equation*}
with $H$ denoting the Heaviside function.
  
  Now suppose that the statement holds for $A_{k-\ell'}f$ for $0 \leq
  \ell' \leq \ell -1$.  Then $P_{0}A_{k-\ell}f$ must satisfy
  \begin{align*}
    P_{0}A_{k-\ell}f &= - \sum_{j=0}^{\ell-1}
    v^{-\imath\sigma_{0}-1}H(v)(\log|v|)^{j}
    \frac{(-\imath)^{j}}{j!}\left( \sigma_{0}\phi_{k-(\ell-1-j)} +
      \phi_{k-(\ell-2-j)}\right) \\
    &\quad\quad\quad+ O(v^{-\imath\sigma_{0}}(\log|v|)^{\ell}),
  \end{align*}
  where again the $O(v^{-\imath\sigma_{0}}(\log|v|)^{\ell})$ term has
  an expansion of a similar form.  The theorem of Haber--Vasy and
  Lemma~\ref{lem:support-of-polar-part} then imply that $A_{k-\ell}f$
  has an expansion
  \begin{align*}
    A_{k-\ell}f = \sum_{j=0}^{\ell} H(v) v^{-\imath\sigma_{0}}(\log|v|)^{j}a_{j}.
  \end{align*}
  Applying $P_{0}$ and equating coefficients finishes the proof in
  this case.

  Finally, if $-\imath \sigma_{0} = 0$, the same argument as in the
  case of $-\imath\sigma_{0} > 0$ still works, but differentiating the
  $j=0$ term yields a $\delta(v)$ term.  This term is no problem, as
  we still simply solve for its coefficient.  This process yields an
  identical result.
\end{proof}

\section{An asymptotic expansion}
\label{sec:an-asympt-expans}

In this section we detail the iteration scheme used to obtain a
preliminary asymptotic expansion for (smooth) solutions $w$ of
$\Box_{g}w \in \dCI (M)$ that vanish in a neighborhood of
$\overline{C_{-}}$.  

We start with a tempered solution $w$ of $\Box_{g} w = f \in
\dCI (M)$ vanishing in a neighborhood of
$\overline{C_{-}}$.  We immediately replace $w$ by $\chi w$, where
$\chi$ is a cut-off function supported near the boundary and vanishing
identically near $\overline{C_{-}}$.  By choosing $\chi$
appropriately, we guarantee that $\chi w = w$ near the boundary and
that the replacement is supported in a collar neighborhood of the
boundary and still solves an equation of the same form.

Because $w$ is tempered, we know that $w \in \rho
^{\gamma}H^{s_{0}}_{b}(M) = H^{s_{0},\gamma}_{b}(M)$ for some $s_{0}$
and $\gamma$.  \emph{We decrease $s_{0}$ so that $\gamma + s_{0} <
  1/2$.}  Corollary~\ref{cor:radial-b-estimate} and our non-trapping
hypothesis then imply that $w$ has module regularity at $\Lambda^{+} =
N^{*}S_{+}$ relative to $H^{s_{0},\gamma}_{b}(M)$.

As before, write
\begin{equation*}
  L = \rho ^{-(n-2)/2}\rho^{-2}\Box_{g} \rho^{(n-2)/2}\in \Diff_{b}^{2}(M),
\end{equation*}
so that setting
\begin{equation*}
  u = \rho ^{-(n-2)/2}w
\end{equation*}
we have
\begin{equation*}
  Lu = g \in \dCI (M),
\end{equation*}
with $u,g$ vanishing near $\overline{C_{-}}$.  Now let $N(L)$ denote
the normal operator of $L$ and set $E= L-N(L)$; $E$ thus measures the
failure of $L$ to be dilation-invariant in $\rho$.  Thus,
\begin{equation*}
  E \in \rho \Diff^{2}_{b}(M).
\end{equation*}

By the form of $G^{-1}$ given by
equation~\eqref{eq:inverse-perturbation}(3.3), we note that the
coefficient of $D_{v}^{2}$ in $E$ is of the form $O(\rho^{2}) + O(\rho
v)$ and hence $Ew \in \rho^{\gamma + 1}H^{s_{0}-1}_{b}(M)$.  

At this juncture, we discuss the mapping properties of the Mellin
conjugate of $E$.  To begin, we let $R_{\sigma}$ be the family of
operators satisfying
\begin{equation*}
  \mathcal{M}\circ E = R_{\sigma} \circ \mathcal{M};
\end{equation*}
thus $R_{\sigma}$ is an operator on meromorphic families in $\sigma$
in which $\rho D_{\rho}$ is replaced by $\sigma$ and multiplication by
$\rho$ translates the imaginary part.  Since, as remarked above, the
coefficient of $D_{v}^{2}$ in $E$ is of the form $O(\rho^{2}) + O(\rho
v)$, i.e., is a sum of terms having better decay either in the sense
of $v$ or $\rho$ than the rest of the operator, we have the following
result on the mapping properties of $R_{\sigma}$ (cf.\ Section~\ref{section:mellin} above for
notation):
\begin{lemma}
  \label{lemma:Rmappingproperties}
  For each $\nu, k, \ell, s$ the operator family $R_{\sigma}$ enjoys
  the following mapping properties:
  \begin{enumerate}
  \item $R_{\sigma}$ enlarges the region of holomorphy at the cost of
    regularity at $\Lambda^{+}$:
    \begin{align}
      \label{Rmapping}
      R_{\sigma} : &\holcon[-k]{\nu}{s} \\ 
      &\to 
      \holcon[-k+2]{\nu + 1}{s-1}  \notag \\
      &\quad\quad +
        \holcon[-k+2]{\nu + 2}{s-2} \notag
    \end{align}
  \item If $f_{\sigma}$ vanishes near $\overline{C_{-}}$ for $\Im
    \sigma \geq -\nu$, then $R_{\sigma}f_{\sigma}$ also vanishes near
    $\overline{C_{-}}$ for $\Im \sigma \geq - \nu - 1$.
  \item If
    \begin{equation*}
      \phi \in \holcon{\nu}{\infty}
    \end{equation*}
    then
    \begin{align*}
      &R_{\sigma} \left( (\sigma - \sigma_{0})^{-k-1}v_{\pm\imath
          0}^{-\imath\sigma_{0}}(\log (v\pm \imath 0))^{k}\phi\right)
      \\
      &\quad = (\sigma - (\sigma_{0}-\imath))^{-k-1}v_{\pm\imath
        0}^{-\imath\sigma_{0}-1}\sum_{j=0}^{k}(\log(v\pm\imath
      0))^{j}\tilde{\phi}_{j,1} \\
      &\quad \quad + (\sigma - (\sigma_{0}-2\imath))^{-k-1}v_{\pm\imath
        0}^{-\imath\sigma_{0}-2}\sum_{j=0}^{k}(\log(v\pm \imath
      0))^{j}\tilde{\phi}_{j,2},
    \end{align*}
    where $\tilde{\phi}_{j,i}$ enjoy the same same properties and are
    holomorphic on $\im \sigma \geq -\nu-1$.
  \item If
    \begin{equation*}
      \phi \in \holcon{\nu}{\infty}
    \end{equation*}
    then
    \begin{align*}
      R_{\sigma} \left( v_{\pm\imath 0}^{-\imath\sigma} \phi \right)
      &\in v^{-\imath\sigma -1}_{\pm\imath 0} \holcon{\nu + 1}{\infty} \\
      &\quad \quad + v^{\imath\sigma -2}_{\pm\imath 0}\holcon{\nu + 2}{\infty}
    \end{align*}
  \end{enumerate}
\end{lemma}

The proof is a simple application of the properties of the Mellin
transform discussed in Section~\ref{section:mellin}.  Note that
in the first term the Sobolev order has decreased by $1$ arising from
the action of the $O(\rho)vD_{v}^{2}$ term in $\Box_{g}$ (rather than
by $2$ as would be the effect of $O(\rho)D_{v}^{2}$ term).  In the
second term, we see the action of $O(\rho^{2})D_{v}^{2}$ terms, which
give a family holomorphic in an even larger strip, at the cost of
further worsening of Sobolev regularity.  We also lose at high
frequency owing to the $(\rho D_{\rho})^{2}$ error term in the
rescaled $\Box_{g}$, which Mellin transforms to an $O(\sigma^{2})$.
(We further note that a sharper result holds, keeping precise
accounts of tradeoffs between $\sigma$ powers and Sobolev orders in the boundary, 
but this refinement
will not be needed for our argument.) 

We now apply the Mellin transform to the identity $Lu = g$, splitting
up $L = N(L) + E$ to obtain
\begin{equation}
  \label{mellinbox}
  P_{\sigma} \tu_{\sigma} = \tg _{\sigma} - R_{\sigma}\tu_{\sigma},
\end{equation}
where, as above, $P_{\sigma} = \widehat{N}(L)$.  As $g\in \dCI (M)$,
we have
\begin{equation*}
  \tg_{\sigma} \in \holcon{C}{s} \text{ for all }C,s.
\end{equation*}

Because $w = \rho^{(n-2)/2}u \in H^{s_{0},\gamma}_{b}(M)$,
Lemma~\ref{lemma:mellinspaces} implies that
\begin{equation}
  \label{initialspace}
  \tu_{\sigma} \in  \holcon[\max (0,-s_{0})]{\varsigma_{0}}{s_{0}},
\end{equation}
where
\begin{equation*}
  \varsigma_{0} = \gamma - (n-2)/2.
\end{equation*}
(Recall that we have already reduced $s_{0}$ so that $s_{0} + \gamma
< 1/2$ and Corollary~\ref{cor:radial-b-estimate} applies.)

As $u$ vanishes near $\overline{C_{-}}$, $\tu_{\sigma}$ also vanishes
there.  Thus, in the notation of
Section~\ref{sec:results-citevasy}, the right-hand-side of
equation~\eqref{mellinbox} is in $\mathcal{Y}^{s_{\tow}-1}$ and
$\tu_{\sigma} \in \mathcal{X}^{s_{\tow}}$.  Here we may choose
$s_{\tow}$ to be constant on the singular support of $\tu_{\sigma}$ as
$\tu_{\sigma}$ is smooth near $\overline{C_{-}}$; in fact, we may take
it to be constant except in a small neighborhood of $\overline{C_{-}}$
where $\tu_{\sigma}$ vanishes.  We take $s_{\tow}$ equal to $s_{0}$
outside this neighborhood.

Because $w$ is conormal with respect to $N^{*}S_{+} = \{ \rho = v =
0\}$, Lemma~\ref{lemma:mellinspaces} implies that
\begin{equation*}
  \tu_{\sigma} \in \holcon{\varsigma_{0}}{-\infty}.
\end{equation*}
Thus, by interpolation with equation~\eqref{initialspace},
\begin{equation*}
  \tu_{\sigma} \in \holcon{\varsigma_{0}}{s_{0}-0}.
\end{equation*}

By Lemma~\ref{Rmapping}, then,
\begin{align}
  \label{resolvent1}
  R_{\sigma}\tu_{\sigma} &\in \holcon{\varsigma_{0}+1}{s_{0}-1-0} \\
  &\quad \quad +
  \holcon{\varsigma_{0}+2}{s_{0}-2-0}, \notag
\end{align}
and hence $P_{\sigma}\tu_{\sigma}$ lies in this space as well.  

Because $P_{\sigma} \tu_{\sigma}$ is now known to be holomorphic in a
larger half-plane, we can now invert $P_{\sigma}$ to obtain meromorphy
of $\tu_{\sigma}$ on this larger space: by Propositions~\ref{proposition:fredholm}
and \ref{prop:nontrappingestimates}, $P_{\sigma}$ is Fredholm as
a map
\begin{equation*}
  \mathcal{X}^{s_{\tow}} \to \mathcal{Y}^{s_{\tow}-1}
\end{equation*}
and $P_{\sigma}^{-1}$ has finitely many poles in any horizontal strip
$\Im z \in [a,b]$, and satisfies polynomial growth estimates as
$\smallabs{\Re z} \to \infty$.  

Here we recall from Section~\ref{sec:results-citevasy} that given
any $\varsigma'$, in order for $P_{\sigma}$ to be Fredholm for $\sigma
\in \CC_{\varsigma'}$, the (constant) value $s(S_{+})$ assumed by the
variable Sobolev order $s_{\tow}$ near $S_{+}$ must satisfy $s(S_{+})
< 1/2 - \varsigma'$.  In other words, as one enlarges the domain of
meromorphy for $\tu_{\sigma}$, one must relax control of the
derivatives.  

We then conclude that $\tu_{\sigma}$ is the sum of two terms arising
from the two terms on the right side of equation~\eqref{resolvent1}.
For any $\delta > 0$, the first term is meromorphic in
$\CC_{\varsigma_{0}+1}$ with values in
\begin{equation}
  \label{rhs1}
  \langle\sigma \rangle^{-\infty} L^{2}(\reals , H^{\min (s_{0}, 1/2 -
  (\varsigma_{0}+1) - \delta)})
\end{equation}
with finitely many poles in this strip arising from the poles of
$P_{\sigma}^{-1}$.  (Note the improvement in the Sobolev orders: by
applying $P_{\sigma}^{-1}$ we win back the derivative we lost from
applying $R_{\sigma}$, but only up to the threshold value.)  Likewise, the second term is meromorphic in
$\CC_{\varsigma_{0}+2}$ with values in
\begin{equation}
  \label{rhs2}
  \langle\sigma\rangle^{-\infty} L^{2}(\reals, H^{\min (s_{0}-1, 1/2 -
  (\varsigma_{0} + 2) - \delta)}),
\end{equation}
again with finitely many poles in the strip.  (Here and below we are
ignoring the distinction between $\mathcal{X}^{s_{\tow}}$ and $H^{s}$
as $\tu_{\sigma}$ is trivial by hypothesis on the set where the
regularity in the variable-order Sobolev space differs from $H^{s}$.)

We now refine the description of the terms in
equations~\eqref{rhs1} and \eqref{rhs2} in two steps using what we know about the
conormal structure of solutions to inhomogeneous equations involving
$P_{\sigma}$.  To begin, since $P_{\sigma}$ maps the terms in question
to conormal spaces, they must in fact lie in the conormal spaces
\begin{align*}
  &\langle\sigma\rangle^{-\infty} L^{2}(\reals , I^{(\min (s_{0} -0,
    1/2 - (\varsigma_{0}+1) - 0))}), \text{ resp.,} \\
  &\langle\sigma\rangle^{-\infty} L^{2}(\reals , I^{(\min (s_{0}-1-0,
    1/2 - (\varsigma_{0}+2)-0))}).
\end{align*}
This improvement follows by
Proposition~\ref{proposition:bpropagation} (propagation of
singularities away from the radial points) and the first case of
Theorem~6.3 of \cite{Haber-Vasy:Radial}.

\begin{remark}
  Here Theorem~6.3 of \cite{Haber-Vasy:Radial} is applied pointwise in
  $\sigma$.  The result there is not stated in terms of bounds (just
  as a membership in the claimed set), but, just as in the case of
  Proposition~\ref{prop:radial-b-estimate} here, estimates can be
  recovered from the statement of Theorem~6.3 by the closed graph
  theorem or alternatively recovered from examination of the proof
  (which proceeds via such estimates).
\end{remark}

Finally, since we have now established conormality, we may use
Corollary~\ref{cor:Mellin-wave-asymp} to determine the detailed
structure of the conormal singularities.  We find that for any $\delta
' > 0$, the two terms in question lie in
\begin{align*}
  &\langle\sigma\rangle^{-\infty}L^{2}(\reals, I^{(s_{0}-\delta')}) + v_{+\imath 0}^{-\imath
    \sigma}\langle\sigma\rangle^{-\infty}L^{2}(\reals , \CI) +
  v_{-\imath
    0}^{-\imath\sigma}\langle\sigma\rangle^{-\infty}L^{2}(\reals ,
  \CI), \text{ resp., } \\
  &\langle\sigma\rangle^{-\infty}L^{2}(\reals, I^{(s_{0}-1-\delta')}) + v_{+\imath 0}^{-\imath
    \sigma}\langle\sigma\rangle^{-\infty}L^{2}(\reals , \CI) +
  v_{-\imath
    0}^{-\imath\sigma}\langle\sigma\rangle^{-\infty}L^{2}(\reals ,
  \CI).
\end{align*}
Again, Corollary~\ref{cor:Mellin-wave-asymp} does not state the
desired estimates explicitly, but these follow either directly from
the proofs or by an application of the closed graph theorem.

Consequently, we have now established 
\begin{align}
  \label{eq:tu-step-1}
  \tu_{\sigma} &\in \holcon{\varsigma_{0}+1}{s_{0}-0} +
  \holcon{\varsigma_{0}+2}{s_{0}-1-0} \\
  \notag&\quad \notag + v_{+\imath
    0}^{-\imath\sigma}\holcon{\varsigma_{0}+1}{\infty}+ v_{-\imath
    0}^{-\imath\sigma}\holcon{\varsigma_{0}+1}{\infty} \\
  \notag&\quad \notag + \sum _{\Im
    \sigma _{j} > - (\varsigma_{0}+2)}(\sigma - \sigma_{j})^{-m_{j}} a_{j},
\end{align}
where
\begin{equation*}
  a_{j} \in \holcon{\varsigma_{0}+1}{\Im\sigma_{j} + 1/2 - 0}.
\end{equation*}
(Note that we have dropped terms $v_{\pm\imath
    0}^{-\imath\sigma}\holcon{\varsigma_{0}+2}{\infty}$ by absorbing
  them in the terms of the same form, holomorphic in the smaller
  half-space; also recall that $I^{(\infty)}=\CI.$)
As the $a_{j}$ are given by the polar parts of $P_{\sigma}^{-1}$ at
value $\sigma_{j}$ lying in a strip in $\CC$, the coefficients of the
polar part of the sum are described by
Proposition~\ref{prop:boundarylogterms}:
\begin{align}
  \label{eq:expn-coeffs}
  a_{j} &= \sum _{\kappa = 0}^{m_{j}-1}(\sigma- \sigma_{j})^{\kappa}
  \sum _{\ell = 0}^{\kappa} \left(v_{+\imath 0}^{-\imath\sigma_{j}}(\log
    (v + \imath 0))^{\ell}a_{j\kappa\ell+} + v_{-\imath 0}^{-\imath
      \sigma_{j}}(\log (v-\imath 0))^{\ell}a_{j\kappa\ell
      -}\right) \notag \\
  &\quad\quad + O((\sigma - \sigma_{j})^{m_{j}}),
\end{align}
with
\begin{equation*}
  a_{j\kappa\ell\pm} = \frac{(-\imath)^{\ell}}{\ell
    !}\phi_{j,m_{j}-(\kappa - \ell), \pm}
\end{equation*}
and $\phi_{j,\alpha,\pm}\in \CI (X)$.  (If $-\imath\sigma_{j}$ is a
positive integer, then the term is of the form in
Proposition~\ref{prop:boundarylogterms}.)  We may further arrange
that $a_{j\kappa\ell\pm}$ are holomorphic and rapidly decaying in
strips.  Also, $\supp a_{j} \cap \overline{C_{-}} = \emptyset$ by
Proposition~\ref{prop:block-structure}.  

Finally, we remark that
\begin{align}
  \label{eq:smooth-terms-to-conormal}
  &v^{-\imath\sigma}_{+\imath 0}\left(
    \holcon{\varsigma_{0}+1}{\infty}\right) \\
  \notag &\quad + v^{-\imath\sigma}_{-\imath 0}\left(
    \holcon{\varsigma_{0}+1}{\infty}\right) \\
  \notag &\quad \subset \holcon{\varsigma_{0}+1}{1/2 + \Im \sigma -0} \\
  \notag &\quad \subset
  \holcon{\varsigma_{0}+1}{1/2 - (\varsigma_{0}+1)-0},
\end{align}
so that we have established
\begin{align}
  \label{eq:tu-step-1A}
  \tu_{\sigma} \in &\holcon{\varsigma_{0}+1}{\min (s_{0}-0, 1/2 -
    (\varsigma_{0}+1) - 0)} \\
  \notag&+ \holcon{\varsigma_{0}+2}{s_{0}-1-0} \\
  \notag&+ \sum_{\Im \sigma_{j} > - (\varsigma_{0}+2)}(\sigma - \sigma_{j})^{-m_{j}}a_{j}
\end{align}
with $a_{j}$ given by
equation~\eqref{eq:expn-coeffs}.  (Through careful accounting
  it is possible to keep track of the $v^{-\imath\sigma}_{\pm\imath
    0}$ terms.  Indeed, these terms contribute both to the radiation field
  and to other terms in the expansion after the blow-up, but become
  rapidly decaying along the front face after taking the Mellin
  transform.) 
  
Now we iterate this argument: by equation~\eqref{eq:smooth-terms-to-conormal},
Lemma~\ref{lemma:Rmappingproperties}, and
Lemma~\ref{lemma:approx-soln-expn},
\begin{align*}
  R_{\sigma}\tu_{\sigma} &\in \holcon{\varsigma_{0}+2}{\min(s_{0}-1-0,
    1/2 - (\varsigma_{0}+2) - 0)} \\
  &\quad +  \holcon{\varsigma_{0}+3}{\min (s_{0}-2-0, 1/2 - (\varsigma_{0}+3)-0)} \\
 &\quad + \holcon{\varsigma_{0}+4}{s_{0}-3-0} \\
  &\quad + \sum_{\Im\sigma_{j} > - (\varsigma_{0} + 2)}(\sigma -
  (\sigma_{j}-\imath))^{-m_{j}}b_{j}' \\
  &\quad + \sum_{\Im\sigma_{j} > - (\varsigma_{0}+2)}(\sigma -
  (\sigma_{j} - 2\imath))^{-m_{j}}b_{j}''
\end{align*}
where the polar parts lie in the spaces
\begin{align*}
  b_{j}' &\in \holcon{\varsigma_{0}+2}{\Im\sigma_{j} - 1/2 - 0} \\
  b_{j}'' &\in \holcon{\varsigma_{0} + 2}{\Im \sigma_{j} - 3/2 - 0},
\end{align*}
and have the forms
\begin{align*}
  b_{j}' &= \sum_{\kappa = 0}^{m_{j}}\sum_{\ell = 0}^{\kappa}(\sigma -
  (\sigma _{j} - \imath))^{\kappa}\Big( v_{+\imath
    0}^{-\imath\sigma_{j}-1}(\log (v+\imath 0))^{\ell}b_{j\kappa\ell
    +}'\\
&\qquad\qquad\qquad\qquad\qquad\qquad +
    v_{-\imath 0}^{-\imath \sigma_{j}-1}(\log
    (v- \imath 0))^{\ell}b_{j\kappa\ell
      -}'\Big) \notag \\
  &\qquad\qquad+ O((\sigma - (\sigma_{j}-\imath))^{m_{j}})
\end{align*}
and
\begin{align*}
   {b}_{j}'' = &\sum_{\kappa =0}^{m_{j}-1}\sum_{\ell
    =0}^{\kappa}(\sigma-(\sigma_{j}-2\imath))^{\kappa}\Big(
    v_{+\imath 0}^{-\imath\sigma_{j}-2}(\log (v+\imath
    0))^{\ell}b_{j\kappa\ell +}''\\
&\qquad\qquad\qquad\qquad\qquad\qquad +
    v_{-\imath 0}^{-\imath \sigma_{j}-2}(\log
    (v-\imath 0))^{\ell}b_{j\kappa\ell
      -}''\Big) \notag \\
  &\qquad\qquad+ O((\sigma - (\sigma_{j}-2\imath))^{m_{j}}).
\end{align*}
Moreover, the $b_{j\kappa\ell\pm}'$ and $b_{j\kappa\ell\pm}''$ are
smooth and supported away from $\overline{C_{-}}$.  

Again inverting $P_{\sigma}$ and employing
Proposition~\ref{prop:boundarylogterms} and
Corollary~\ref{cor:Mellin-wave-asymp} yields
\begin{align*}
  \tu_{\sigma} &\in \holcon{\varsigma_{0}+2}{\min (s_{0}-0, 1/2 -
    (\varsigma_{0}+1)-0)} \\
  &\quad + \holcon{\varsigma_{0}+3}{\min (s_{0} - 1 - 0, 1/2 -
    (\varsigma_{0}+2)-0)} \\
  &\quad + \holcon{\varsigma_{0}+4}{s_{0}-2} \\
  &\quad + v_{+\imath 0 }^{-\imath\sigma}
  \holcon{\varsigma_{0}+2}{\infty} \\
  &\quad + v_{-\imath 0}^{-\imath\sigma}
  \holcon{\varsigma_{0}+2}{\infty} \\
  &\quad + \sum_{\Im\sigma_{j} > -(\varsigma_{0}+4)}(\sigma -
  \sigma_{j})^{-m_{j}}a_{j} \\
  &\quad + \sum_{\Im \sigma_{j}>-(\varsigma_{0}+2)}(\sigma -
  (\sigma_{j}- \imath))^{-\tilde{m}_{j}}\tilde{a}_{j1} \\
  &\quad + \sum_{\Im\sigma_{j}> - (\varsigma_{0}+2)}(\sigma -
  (\sigma_{j} - 2\imath))^{-\tilde{m}_{j}} \tilde{a}_{j2},
\end{align*}
where again the coefficients of the poles have expansions as in
equation~\eqref{eq:expn-coeffs} (although the expansion for
$\tilde{a}_{j2}$ begins at $v^{-\imath\sigma_{j}-1}$) and support away
from $\overline{C_{-}}$.  Here we may have $\tilde{m_{j}} > m_{J}$ if
there are integer coincidences among the poles of $P_{\sigma}^{-1}$,
i.e., if $\sigma_{j}$ and $\sigma_{j}-\imath$ or $\sigma_{j} -
2\imath$ are both poles.  As before we may use the inclusions of
$v_\pm^{-\imath\sigma}$ in conormal spaces to rewrite this (for
purposes of the next step of our iteration) as
\begin{align*}
  \tu_{\sigma} &\in \holcon{\varsigma_{0}+2}{\min (s_{0}-0, 1/2 -
    (\varsigma_{0}+2)-0)} \\
  &\quad + \holcon{\varsigma_{0}+3}{\min (s_{0} - 1 - 0, 1/2 -
    (\varsigma_{0}+2)-0)} \\
  &\quad + \holcon{\varsigma_{0}+4}{s_{0}-2} \\
 &\quad + \sum_{\Im\sigma_{j} > -(\varsigma_{0}+4)}(\sigma -
  \sigma_{j})^{-m_{j}}a_{j} \\
  &\quad + \sum_{\Im \sigma_{j}>-(\varsigma_{0}+2)}(\sigma -
  (\sigma_{j}- \imath))^{-\tilde{m}_{j}}\tilde{a}_{j1} \\
  &\quad + \sum_{\Im\sigma_{j}> - (\varsigma_{0}+2)}(\sigma -
  (\sigma_{j} - 2\imath))^{-\tilde{m}_{j}} \tilde{a}_{j2},
\end{align*}

Iterating in this fashion, we obtain after $N$ such steps (and
slightly weakening our Sobolev exponents over the foregoing for the sake of
simplicity):
\begin{align}
  \label{foo:2.28}
  \tu_{\sigma} &\in \holcon{\varsigma_{0}+N}{\min (s_{0}-0, 1/2 -
    (\varsigma_{0}+N-1)-0)} \\
  \notag &+ \dots + \holcon{\varsigma_{0}+2N}{\min (s_{0}-N-0, 1/2 -
    (\varsigma_{0} + 2N-1)-0)} \\
  \notag &+ v_{+\imath 0}^{-\imath\sigma} \left(
    \holcon{\varsigma_{0}+N}{\infty}\right) \\
  \notag &+ v_{-\imath 0}^{-\imath\sigma} \left(
    \holcon{\varsigma_{0}+N}{\infty} \right) \\
  \notag &+ \sum_{\Im\sigma_{j} > -(\varsigma_{0}+2N)}(\sigma -
  \sigma_{j})^{-m_{j}}a_{j} \\
  \notag &+ \sum_{\ell=1}^{N}\sum_{\Im\sigma_{j} > - (\varsigma_{0} +
    2N)}(\sigma - (\sigma_{j} - \imath\ell))^{-\tilde{m}_{j}}\tilde{a}_{j\ell}.
\end{align}
Here again $\tilde{m}_{j}$ may exceed $m_{j}$ in the case of
integer coincidences among poles of $P_{\sigma}^{-1}$.  Moreover,
while $a_{j}$ is described by equation~\eqref{eq:expn-coeffs}, we also
have
\begin{align*}
  \tilde{a}_{j\ell} = \sum_{\kappa = 0}^{\tilde{m}_{j}-1}(\sigma -
  (\sigma_{j} - \imath \ell))^{\kappa} \sum_{k = 0}^{\ell -1} \sum_{p
    = 0}^{P(j,\ell, \kappa , k)} \Big( &v_{+\imath
      0}^{-\imath\sigma_{j}-k} (\log (v+\imath 0))^{p}a_{j\ell\kappa
      kp+} \\
    & + v_{-\imath 0}^{-\imath \sigma_{j}-k}(\log (v-\imath
    0))^{p}a_{j\ell\kappa kp-}\Big) \\
    & + O((\sigma - (\sigma_{j} - \imath \ell))^{\tilde{m}_{j}}).
\end{align*}

As the inverse Mellin transform of $(\sigma - \sigma_{0})^{-m}$ is
\begin{equation*}
  \frac{\imath ^{m}}{(m-1)!} \rho^{\imath\sigma_{0}}(\log \rho)^{m-1},
\end{equation*}
under inverse Mellin transform with a contour deformation to the line
$\reals - \imath(\varsigma_{0} + N)$, equation~\eqref{eq:expn-coeffs}
shows that the poles in the sum $\sum (\sigma -
\sigma_{j})^{-m_{j}}a_{j}$ yield the residues
\begin{equation*}
  \sum_{j} \sum_{\kappa = 0}^{m_{J}}\sum_{\ell=0}^{\kappa}\frac{\imath
    ^{m_{j} - \kappa}(-\imath)^{\ell}}{(m_{j}-\kappa - 1)!\ell !} \rho
  ^{\imath\sigma_{j}}v^{-\imath\sigma_{j}} (\log
  \rho)^{m-\kappa-1}(\log v)^{\ell}\phi_{j,m-(\kappa -\ell)},
\end{equation*}
i.e., the main terms in our asymptotic expansion.  (Strictly speaking,
this is an expansion in powers of $v_{\pm\imath 0}$ rather than in
powers of $v;$ however, we are primarily concerned with asymptotics in
the regime $v/\rho\gg 0.$ Likewise, we have chosen to write
$\phi_{\bullet}=\phi_{-,\bullet}+\phi_{+,\bullet},$ but could carry
along both terms separately if desired.)  Rearranging this sum yields
\begin{equation*}
  \sum_{j}\sum_{k=0}^{m_{j}-1}\frac{\imath^{-k+1}}{k!}\rho^{\imath\sigma_{j}}v^{-\imath\sigma_{j}}
  (\log v - \log \rho)^{k}\phi_{j,k+1},
\end{equation*}
i.e., the only logarithmic terms in this sum are powers of $\log v -
\log \rho$.

The terms in the last sum in equation~\eqref{foo:2.28} become
\begin{equation*}
  \sum_{j}\sum_{\ell=1}^{N}\sum_{\kappa + \alpha < \tilde{m}_{j\ell}}
  a_{j\ell\kappa\alpha}\rho^{\imath\sigma_{j}+\ell}v^{-\imath\sigma_{j}-\ell
  + 1} \smallabs{\log\rho}^{\kappa}\smallabs{\log v}^{\alpha}.
\end{equation*}
The other, ``remainder'' terms in equation~\eqref{foo:2.28} lie in
\begin{align}
  \label{eq:remainder-of-foo:2.28}
  \sum_{j=0}^{N}&\holcon{\varsigma_{0}+N+j}{\min(s_{0}-j-0, 1/2 -
    (\varsigma_{0}+N+j-1)-0)} \\
  \notag & + v_{+\imath 0}^{-\imath\sigma}
  \left( \holcon{\varsigma_{0}+N}{\infty}\right) \\
  \notag & + v_{-\imath 0}^{-\imath\sigma}\left( \holcon{\varsigma_{0}
      + N}{\infty}\right).
\end{align}
Thus by Lemma~\ref{lemma:mellinspaces}, after applying the
inverse Mellin transform the terms in the first sum of
equation~\eqref{eq:remainder-of-foo:2.28} become
\begin{equation*}
  O\left( \rho^{\varsigma_{0}+N+j-0}v^{\min (s_{0}, 1/2 -
      (\varsigma_{0} + N -1)) - 1/2 - j - 0}\right).
\end{equation*}
In particular, the power of $v$ appearing is at least one larger than
the power of $\rho^{-1}$.  The last two terms of equation~\eqref{eq:remainder-of-foo:2.28}
Mellin transform to Schwartz functions of $v/\rho$.  

Thus, returning to the solution $w$ of $\Box_{g}w = f \in \dCI (M)$,
and taking $N$ large to simplify the remainder, we find that near $\{
\rho = v = 0\}$, $w$ has an asymptotic expansion of the form
\begin{align}
  \label{theexpansion}
  \notag w = &\rho^{(n-2)/2}\sum_{j}\sum_{\kappa \leq
    m_{j}}\rho^{\imath\sigma_{j}}v^{-\imath\sigma_{j}}(\log v - \log
  \rho)^{\kappa}a_{j\kappa} \\
  &+ \rho^{(n-2)/2}\sum_{j}\sum_{\ell=1}^{N}\sum_{\kappa + \alpha <
    \tilde{m}_{j\ell}}
  \tilde{a}_{j\ell\kappa\alpha}\rho^{\imath\sigma_{j}+\ell}
  v^{-\imath\sigma_{j}- \ell + 1}\smallabs{\log
    \rho}^{\kappa}\smallabs{\log v}^{\alpha} + w'
\end{align}
with
\begin{equation*}
  w' \in \sum_{j=0}^{N}\rho^{(n-2)/2+\varsigma_{0} + N + j - 0}v^{-\varsigma_{0}-N-j+1-0}L^{\infty},
\end{equation*}
where $\sigma_{j}$ are the poles of the meromorphic inverse
$P_{\sigma}^{-1}$, and the coefficients are the corresponding
resonance states.  Here $v^{-\imath\sigma_{j}}a_{j\kappa}$ (and its
counterpart in the second sum) is understood to mean a sum of the two
$v_{\pm\imath 0}^{-\imath\sigma_{j}}$ terms (which we write out fully
below).

\begin{remark}
Note that the only log terms in the
$\rho^{\imath\sigma_{j}}v^{-\imath\sigma_{j}}$ term occur as powers of
$\log v - \log \rho$.  Because $\log v - \log \rho = \log s$ in the
radiation field blow-up, this will imply that $\rho^{-(n-2)/2}w$ has a
restriction to the front face of the radiation field blowup.
\end{remark}

\section{The asymptotics of the radiation field}
\label{sec:asympt-radi-field}

We now lift the expansion \eqref{theexpansion} to the radiation field blowup.
We thus introduce the ``radiation field'' coordinates
$\rho,y,s=v/\rho;$ note that these constitute a coordinate system on
the blown up space described in Section~\ref{sec:blowup}, and note
that $\pd[s]$ is well-defined as a vector field on the fibers of $\mathrm{ff}$.  In these
coordinates, then, homogeneity yields the expansion
\begin{align*}
  &\sum_j \sum_{\alpha+\kappa< m_j} (\log s)^{\alpha}\big(
  a'_{j\kappa\alpha,+} s_{+\imath
    0}^{-\imath\sigma_{j}}+a'_{j\kappa\alpha,-} s_{-\imath 0}^{-\imath\sigma_{j}}\big) \\
  &+ \sum_{j}\sum_{\ell =1}^{N}\sum_{\kappa +\alpha <
    \tilde{m}_{j\ell}}\rho^{\ell}\smallabs{\log\rho}^{\kappa}\left(
    \log\rho + \log s\right)^{\alpha} \left( a_{j\ell\kappa\alpha,
      +}s_{+\imath 0}^{-\imath\sigma_{j}} + a_{j\ell\kappa\alpha,
      -}s_{-\imath 0}^{-\imath\sigma_{j}}\right)\\
  &+u'
\end{align*}
for  $u=\rho^{-\frac{n-2}{2}}w$.
Consequently, restricting terms of the expansion to
$\rho=0$ yields an expansion
$$
\sum_j \big( a_{jk,+} s_{+\imath 0}^{-\imath\sigma_{j}}+a_{jk,-}
s_{-\imath 0}^{-\imath\sigma_{j}}+\tilde \tu_{jk}\big) 
$$
with a remainder term $u'$.  Notice that the presence of $\log
\rho$ factors in the $\rho^0$ ($\ell=0$) terms would prevent the restriction of $u$ to the front face of
the blow-up, but in Section~\ref{sec:an-asympt-expans} we showed that
in fact (at top-order) those terms possessing a logarithmic factor
cancel. We can now define the radiation field as in Section~\ref{sec:blowup}:
\begin{definition}
  If $w$ is a solution of $\Box_{g}w = f$, $f\in\dCI(M)$,
  $w$ vanishing near $\overline{C_-}$, we define the
  (forward) \emph{radiation field} of $w$ by
  \begin{equation*}
    \mathcal{R}_{+}[w](s, y) = \pd[s] u(0, s, y),\qquad u=\rho^{-\frac{n-2}{2}}w.
  \end{equation*}
\end{definition}
The rest
of Theorem~\ref{thm:main-theorem} follows immediately.  As identified in
Section~\ref{sec:ident-p_sigma-1}, the exponents $\sigma_{j}$ are the
poles of $\mathcal{R}_{C_{+}}(\sigma)$, i.e., the resonances of the
asymptotically hyperbolic problem on the cap $C_{+}$, while the terms
supported at $S_{+}$ do not contribute to the expansion as $s\to
\infty$.  

\begin{remark}
  \label{rem:dont-see-deltas}
  While it may seem that the coefficients in the expansion are
  singular at $s=0$, this is an artifact of the basis chosen.  The
  b-regularity established in Section~\ref{sec:b-regularity} (see, in
  particular, Remark~\ref{rem:blow-up-conormal}) implies that the
  solution is conormal to the front face of the radiation field
  blow-up and hence the coefficients may be taken to be smooth.
\end{remark}

\subsection{Asymptotically Minkowski space}
\label{sec:minkowski-space-1}

We now consider the special case of asymptotically Minkowski space
(i.e., ``normally very short range'' perturbations of Minkowski space).  Here
we are assuming that the metric takes the form \eqref{Minkowskimetric}
modulo
$$
O(\rho^{2}) \frac{d\rho^{2}}{\rho^{4}} + O(\rho) \frac{d\rho
  d(v,y)}{\rho^{3}} + O(\rho) \left( \frac{d(v,y)}{\rho}\right)^{2}.
$$
Then the induced metric on $C_+$ (which is diffeomorphic to a ball) is
the metric on $(n-1)$-dimensional hyperbolic space; $P_{\sigma}$ is a
conjugate of the spectral family on hyperbolic space.  (See Section~5
of \cite{Vasy-Dyatlov:Microlocal-Kerr}.)  In particular, the relevant
poles of $P_{\sigma}^{-1}$ (i.e., those of $L_{\sigma,+}^{-1}$ from
Section~\ref{sec:ident-p_sigma-1}) are given by the poles of the
meromorphic expansion of $\left( \Lap_{\mathbb{H}^{n-1}} - \sigma^{2}
  - \frac{(n-2)^{2}}{4}\right)^{-1}$.  These poles can be calculated
explicitly: when $n$ is even (and hence the spatial dimension is odd),
there are no poles, while if $n$ is odd, the poles are given by
$\sigma _{j} = - \imath \frac{n-2}{2} - \imath j$ for $j \in \NN$.  In
particular, $\mathcal{R}_{+}[w]$ has an asymptotic expansion of the
following form:
\begin{equation*}
  \mathcal{R}_{+}[w](s,\omega) \sim
  \begin{cases}
    O(s^{-\infty}) & n \text{ even} \\
    \sum_{j=0}^{\infty} \sum_{\kappa \leq j} s^{-\frac{n}{2}-j}(\log
    s)^{\kappa} a_{j\kappa} & n \text{ odd.}
  \end{cases}
\end{equation*}
(Recall that one differentiates the restriction of $u$ in $s$ to obtain
$\mathcal{R}_+.$)  In the special case when the metric is in fact
\emph{exactly} Minkowski in a neighborhood of $C_+$ in $M,$ we remark that
the whole iterative apparatus of Section~\ref{sec:an-asympt-expans} can be
dispensed with, in favor of a single application of $P_\sigma^{-1}$ to the
Mellin-transformed inhomogeneity, with the result that the the $\log$ terms
in the expansion do not appear in that case.

The stability of $P_{\sigma}^{-1}$ under perturbations implies that
for small ``normally short range'' perturbations of Minkowski space,
the radiation field still decays.  In this setting, however, poles of
$P_{\sigma}^{-1}$ that are not poles of $L_{\sigma,+}^{-1}$ (and hence
do not affect the decay of the radiation field) may become relevant
under perturbations.  As discussed earlier, such poles must occur at
purely imaginary negative integers and the corresponding states must be
supported exactly at $S_{+}$.  Such a state occurs even in
4-dimensional Minkowski space at $\sigma = -\imath$.  Under small
``normally short range'' perturbations, then, the first pole occurs
close to $\sigma = -\imath$ and so we conclude that the radiation field
is $O(s^{-2+\epsilon})$ as $s\to \infty$. 

\appendix

\section{Variable order Sobolev spaces}\label{section:variableorder}

First recall that (uniform) symbols $a\in S^m_{\rho,\delta}$ on
$\RR^n\times\RR^n$ of type $(\rho,\delta)$ of order $r\in\RR$ are $\CI$
functions on $\RR^n_z\times\RR^n_\zeta$ such that
$$
|D^\alpha_z D^\beta_\zeta a|\leq C\langle \zeta\rangle^{r+\delta|\alpha|-\rho|\beta|}.
$$
For various applications, the natural type is $\rho=1-\delta$,
$\delta\in[0,1/2)$, with $\delta=0$ corresponding to the standard
symbol class. We assume these restrictions from now on; for us the
relevant regime will be $\delta>0$ but arbitrarily small. Note that
$S^{-\infty}=\bigcap_r S^r_{1-\delta,\delta}$ is independent of $\delta$. There is a
symbol calculus within this class $S^{r}_{1-\delta,\delta}$, which works modulo
$S^{r-1+2\delta}_{1-\delta,\delta}$; the principal symbol of the composition of two operators
is the product of the two principal symbols in this sense. Further,
one has the full symbol expansion of the composition modulo
$\Psi^{-\infty}$; namely if $(Au)(z)=(2\pi)^{-n}\int
e^{i(z-z')\cdot\zeta}a(z,\zeta)u(z')\,dz'$ is the left quantization of
$a\in S^r_{1-\delta,\delta}$, and $B$ is the left quantization of
$b\in S^{r'}_{1-\delta,\delta}$ then $AB$ is the left quantization of
$$
c\sim\sum_\alpha\frac{\imath^{|\alpha|}}{\alpha!}D_\zeta^\alpha a D_z^\alpha b.
$$
As usual these can be transferred to manifolds by local coordinates,
allowing the addition of globally $\CI$ kernels as well.

We can now turn to variable order operators. Suppose that $s$ is a
real-valued function on $S^*\RR^n=\RR^n\times\sphere^{n-1}
=\RR^n\times(\RR^n\setminus\{0\})/\RR^+$, which we assume is constant
outside a compact set since we are interested only in transferring the result to
manifolds via local coordinates -- one could assume instead uniform
bounds on derivatives on $\RR^n\times\sphere^{n-1}$.
On $\RR^n\times\RR^n$, we say that $a$ is a (variable order) symbol of order
$s$, written $a\in S^s_{1-\delta,\delta}$, $\delta\in(0,1/2)$ if
\begin{equation}\label{eq:weight-times-symbol}
a=\langle\zeta\rangle^s a_0,\qquad a_0\in S^0_{1-\delta,\delta}(T^*X).
\end{equation}
So $S^s_{1-\delta,\delta}\subset S^{s_0}_{1-\delta,\delta}$ with
$s_0=\sup s$. Thus, one can quantize these symbols, with the result,
$\Psi^s_{1-\delta,\delta}$ being a subset of
$\Psi^{s_0}_{1-\delta,\delta}$. One calls the equivalence class of $a$
in $S^s_{1-\delta,\delta}/S^{s-1+2\delta}_{1-\delta,\delta}$ the
principal symbol of the left quantization $A$ of $a.$ We could of course just
as well used another choice of quantization such as right- or
Weyl-quantization.  Note, though, that the condition
$\delta>0$ is crucial for making the
different choices of quantizations equivalent since the right reduction formula is
$$
\sim\sum\frac{(-\imath)^{|\alpha|}}{\alpha!} D_z^\alpha D_\zeta^\alpha
  a,
$$
and the derivatives falling on the exponent of
  $\langle\zeta\rangle$ give logarithmic terms, which do not have the
  full $S_{1,0}$ type gain.

The full asymptotic expansion
for composition shows that if $s,s'$
are real valued functions on $S^*\RR^n$ then
$$
A\in\Psi^s_{1-\delta,\delta},\ B\in\Psi^{s'}_{1-\delta,\delta}\Longrightarrow
AB\in\Psi^{s+s'}_{1-\delta,\delta},
$$
and modulo
$\Psi^{s+s'-1+2\delta}_{1-\delta,\delta}$ it is given by a quantization
of the product of the principal symbols; again $\delta>0$ is
important. The commutator $[A,B]$ is then in
$\Psi^{s+s'-1+2\delta}_{1-\delta,\delta}$, and its principal symbol
(modulo $S^{s+s'-2+4\delta}_{1-\delta,\delta}$) is
$\frac{1}{\imath}\{a,b\}$, where $\{.,.\}$ is the Poisson bracket, and
$a,b$ are the respective principal symbols. Defining $a\in S^s_{1-\delta,\delta}$ to
be elliptic if there exists $c,R>0$ such that
$|a|\geq c\langle\zeta\rangle^s$ for
$\langle\zeta\rangle\geq R$, i.e.\ if $a_0$ is elliptic in
\eqref{eq:weight-times-symbol} in the analogous standard sense,
the (microlocal)
elliptic parametrix construction works, i.e.\ if $A\in\Psi^s_{1-\delta,\delta}$ has elliptic
principal symbol then there is $G\in\Psi^{-s}_{1-\delta,\delta}$ such
that $GA-\Id,AG-\Id\in\Psi^{-\infty}$.
We can transfer these operators to manifolds $X$ via localization and
adding $\CI$ Schwartz kernels to the space; here we may assume that
$X$ is compact.
In this manner, for $s$ a real-valued
function on $S^*X=(T^*X\setminus o)/\RR^+$ with
$s_0=\sup s$, we define $\Psi^s_{1-\delta,\delta}(\tilde
X)\subset\Psi^{s_0}_{1-\delta,\delta}(X)$.
The principal symbol of $A\in\Psi^s_{1-\delta,\delta}(X)$ is
a well-defined element of $S^s_{1-\delta,\delta}(T^*X)/S^{s-1+2\delta}_{1-\delta,\delta}(T^*X)$.

We can now define Sobolev spaces: fix $A\in\Psi^s(X)$ elliptic,
$s_1=\inf s$.
We write
$$
H^s=\{\ubdry\in H^{s_1}:\ A\ubdry\in L^2\},\ \|\ubdry\|_{H^s}^2=\|\ubdry\|_{H^{s_1}}^2+\|A\ubdry\|^2_{L^2};
$$
this is a Hilbert space and all the standard mapping properties of ps.d.o's
apply.
Different elliptic choices $A,B\in\Psi^s_{1-\delta,\delta}$ defining $H^s$ give the same
space, since if $G$ is a parametrix for $A$, then $B\ubdry=BGA\ubdry+E\ubdry$, where
$E\in\Psi^{-\infty}$, so $BG\in\Psi^0_{1-\delta,\delta}$, $A\ubdry\in L^2$
shows $B\ubdry\in L^2$ by the standard $L^2$-boundedness of $\Psi^0_{1-\delta,\delta}$,
and also shows the equivalence of the norms.
Further,  if
$s,s'\in\CI(S^*X)$ and $B$ is order $s$ then
$$
B:H^{s'}\to H^{s'-s}
$$
is continuous; taking $\Lambda_s$ elliptic of order $s$, then this is
equivalent to
$$
\Lambda^{s'-s}B\Lambda^{-s}:L^2\to L^2
$$
bounded, but the left hand side is in $\Psi^0_{1-\delta,\delta}$, so
this is again the standard $L^2$ boundedness.

Since the elliptic parametrix construction works, elliptic estimates hold without conditions on $s$ in this
setting. In our considerations, near the radial sets $s$ will be taken constant, so the
previous results apply microlocally there. However, one needs new real principal type
estimates; these hold if $s$ is non-increasing along the direction of
the $\hamvf_p$-flow in which we want to propagate the estimates.

\begin{proposition}\label{proposition:variablepropagation}
Suppose that $P\in\Psi^m(X)$ has real principal symbol.
Suppose that $s\in\CI(S^*X)$ is non-increasing along $\hamvf_p$ on a
neighborhood $O$ of $q\in S^*X$. Let 
$B,G,R\in\Psi^0$, with the property that $\WF'(B)\subset\Ell(G)$ and
such that if $\alpha\in\WF'(B)\cap\Sigma$ then the backward (null-)bicharacteristic
of $p$ from $\alpha$ reaches $\Ell(R)$ while remaining in $\Ell(G)\cap O.$
Then for all $N$ there is $C>0$ such that
$$
\|B\ubdry\|_{H^s}\leq C(\|GP \ubdry\|_{H^{s-m+1}}+\|R\ubdry\|_{H^s}+\|\ubdry\|_{H^{-N}}).
$$

A similar result holds if $s$ is non-decreasing along $\hamvf_p$ and
``backward'' is replaced by ``forward.''
\end{proposition}

Related results appear in \cite{Unterberger:Resolution}, but there the
weights arise from the base space $X$, and logarithmic weights are
used as well, which would require some definiteness of the derivative
of $s$ along $\hamvf_p$ that we do not have here.

\begin{proof}
As the result states nothing about radial points, one may assume that
$\hamvf_p$ is non-radial on $O$.
This then reduces to a microlocal result, namely that there is a neighborhood of a
point $q$ in which the analogous property holds. This can be proved by
a positive commutator estimate as in \cite{Hormander:Existence}.
Let $|\xi|$ be a positive homogeneous degree $1$ elliptic function on
$T^*X$; since we are working microlocally, we may take $|\xi|$
to be the function $|\zeta|$ in local coordinates. With
$\hamvf_{p,m}=|\xi|^{-m+1}\hamvf_p$ denoting the rescaled Hamilton vector
field, which is homogeneous of degree zero, thus a vector field on
$S^*X$, one can introduce local coordinates
$q_1,\ldots,q_{2n-1}$ on $S^*X$ centered at $\alpha$ such that $\hamvf_{p,m}=\frac{\pa}{\pa
  q_1}$; one writes $q'=(q_2,\ldots,q_{2n-1})$. Then one fixes $t_2<t_1<0<t_0$ and a neighborhood $U$ of $0$ in
$\RR^{2n-2}_{q'}$ such that $[t_2,t_0]_{q_1}\times \overline{U_{q'}}\subset O$
and such that one has a priori regularity near $[t_2,t_1]_{q_1}\times
\overline{U_{q'}}$, i.e.\ $R$ in the notation of the proposition is
elliptic there. For $r\in [0,1]$ (the regularization parameter) one considers
$$
a_r=|\xi|^{s-(m-1)/2}\chi(q_1)\phi(q')\psi_r(|\xi|),
$$
where $\phi\in\CI_c(\RR^{2n-2})$ is supported in $U$,
$$
\psi_r(t)=(1+rt)^{-1},
$$
and
$$
\chi(q_1)=\chi_0(q_1)\chi_1(q_1),
$$
with $\chi_0(t)=e^{-\digamma/(t_0-t)}$, $t<t_0$, $\chi_0(t)=0$ for
$t\geq t_0$ and $\chi_1(t)\equiv 1$ near $[t_1,\infty)$, $0$ near
$(-\infty,t_2]$; here $\digamma>0$ will be taken sufficiently
large. Taking $\delta\in(0,1/2)$ arbitrary (i.e.\ $\delta$ can be very
small), $\psi_r$ reduces the order of $a_r$ for $r>0$, so $a_r\in
S^{s-(m-1)/2-1}_{1-\delta,\delta}$ for $r>0$, and for $r\in[0,1]$, $a_r$ is
uniformly bounded in $S^{s-(m-1)/2}_{1-\delta,\delta}$, converging to $a_0$ in
$S^{s-(m-1)/2+\ep}_{1-\delta,\delta}$ for $\ep>0$.
Then as $\hamvf_{p,m}q_1=1$ and $\psi_r'=r\psi_r^2$,
\begin{equation*}\begin{aligned}
\hamvf_p a_r^2
&=2|\xi|^{2s}\phi(q')^2\psi_r(|\xi|)^2\chi_1(q_1)\chi_0(q')\\
&\qquad\qquad\times\Big(\chi'_0(q_1)\chi_1(q_1)+\chi_0(q_1)\chi'_1(q_1)\\
&\qquad\qquad\qquad+(s-(m-1)/2+r|\xi|\psi_r)|\xi|^{-1}(\hamvf_{p,m}|\xi|)\chi_0(q_1)\chi_1(q_1)\\
&\qquad\qquad\qquad+(\log|\xi|)(\hamvf_{p,m}s) \chi_0(q_1)\chi_1(q_1)\Big).
\end{aligned}\end{equation*}
Now, $\chi_0'\leq 0$, giving rise to the main ``good'' term, while the $\chi_1'$ term is supported in $(t_2,t_1)_{q_1}\times
U_{q'}$, where we have a priori regularity and estimates. Further, by making
$\digamma$ large, taking into account that $r|\xi|\psi_r$ is bounded,
we can dominate the $|\xi|^{-1}\hamvf_{p,m}|\xi|$ term
since $\chi_0$ can be bounded by a small multiple of $\chi_0'$ for
$\digamma>0$ large, and $\hamvf_{p,m}s\leq 0$, i.e.\ has the same sign as
the $\chi_0'$ term. The imaginary
(or skew-adjoint in the non-scalar setting)
part of the subprincipal symbol also gives a
contribution that can be dealt with as the $|\xi|^{-1}\hamvf_{p,m}|\xi|$
term.
Thus, taking $A_r$ to have principal symbol $a_r$ and family wave
front set
$\WF'(\{A_r\})=\esssupp a$ (for instance a quantization of $a_r$), $B_r$ have
principal symbol
$$
b_r=|\xi|^s\phi(q')\chi_1(q_1)\sqrt{\chi_0(q')\chi_0'(q')}\psi_r(\xi),
$$
and similar
$\WF'$ one obtains an estimate of the desired kind, and by estimating the
$\chi_1'$ term (which is the only one having the wrong sign) by the
$R$ term, first by obtaining an estimate for $r>0$ and then letting
$r\to 0$ to obtain the result of the desired form.
Corresponding to the symbol class, this can give $1/2-\delta$ order of improvement (i.e.\
allows $-N=s-1/2+\delta$) for all $\delta>0$;
iterating gives the stated result.
\end{proof}

\def\cprime{$'$} \def\cprime{$'$}

\end{document}